\documentclass{amsart}

\title[Algebraic relations among
$\Gamma$-values in characteristic $p$]
{Determination of the algebraic relations
among special
$\Gamma$-values in positive characteristic}
\date{July 19, 2002}
\author[G.\ W.\ Anderson]{Greg W.\ Anderson}
\address{School of Mathematics \\ University of Minnesota \\
Minneapolis, MN 55455}
\email{gwanders@math.umn.edu}
\author[W.\ D.\ Brownawell]{W.\ Dale Brownawell}
\address{Department of Mathematics \\ Penn State University \\
University Park, PA 16802}
\email{wdb@math.psu.edu}
\author[M.\ A.\ Papanikolas]{Matthew A.\ Papanikolas}
\address{Department of Mathematics \\ Brown University \\
Providence, RI 02912}
\email{map@math.brown.edu}

\newcommand{\Kboldbar}{{\overline{\Kbold}}}

\newcommand{\Mbold}{{\mathbf M}}
\newcommand{\Rbold}{{\mathbf R}}
\newcommand{\Lbold}{{\mathbf L}}
\newcommand{\Lboldbar}{{\overline{\Lbold}}}
\newcommand{\Fbold}{{\mathbf F}}

\newcommand{\Kbold}{{\mathbf K}}

\DeclareMathOperator{\rank}{{\mathrm rk}}

\DeclareMathOperator{\Betti}{{\mathrm Betti}}

\newcommand{\hh}{{\mathbf h}}
\newcommand{\ggg}{{\mathbf g}}
\newcommand{\RRR}{{\mathcal R}}
\newcommand{\DDD}{{\mathcal D}}
\DeclareMathOperator{\weight}{{{\mathrm{wt}}}}

\DeclareMathOperator{\Hom}{{\rm Hom}}
\DeclareMathOperator{\Moore}{{\rm Moore}}
\newcommand{\abold}{{\mathbf a}}
\newcommand{\bbold}{{\mathbf b}}

\newcommand{\size}[1]{{\left\Vert#1\right\Vert}}

\newcommand{\diamondbracket}[1]{{\left\langle#1\right\rangle}}
\newcommand{\FF}{{\mathbb F}}
\newcommand{\QQ}{{\mathbb Q}}
\newcommand{\AAA}{{\mathcal A}}
\newcommand{\OO}{{\mathcal O}}
\newcommand{\iso}{\xrightarrow{\sim}}
\newcommand{\RR}{{\mathbb R}}
\newcommand{\CC}{{\mathbb C}}
\newcommand{\Tbar}{{\widetilde{T}}}
\newcommand{\ebold}{{\mathbf e}}
\newcommand{\ZZ}{{\mathbb Z}}
\DeclareMathOperator{\Gal}{{\rm Gal}}
\DeclareMathOperator{\Res}{{\rm Res}}
\DeclareMathOperator{\GL}{{\rm GL}}
\newtheorem{Lemma}[subsubsection]{Lemma}
\newtheorem{Corollary}[subsubsection]{Corollary}
\newtheorem{Proposition}[subsubsection]{Proposition}
\newtheorem{Theorem}[subsubsection]{Theorem}
\theoremstyle{remark}
{\newtheorem{Remark}[subsubsection]{Remark}}
\newcommand{\one}{{\mathbf 1}}
\newcommand{\card}{{\#}}
\newcommand{\kbar}{{\bar{k}}}

\DeclareMathOperator{\Mat}{{\mathrm Mat}}
\DeclareMathOperator{\coker}{{\rm coker}}

\DeclareMathOperator{\End}{{\mathrm End}}

\newcommand{\EE}{{\mathcal E}}

\newcommand{\Xbar}{{\overline{X}}}
\newcommand{\Ubar}{{\overline{U}}}

\begin{document}

\begin{abstract}
We devise a new criterion for linear independence over function
fields.  Using this tool in the setting of dual
$t$-motives, we find that all algebraic relations among special
values of the geometric
$\Gamma$-function over
$\FF_q[T]$ are explained by the standard
functional equations.
\end{abstract}

\maketitle
\tableofcontents

\section{Introduction} \label{section:Introduction}
\subsection{Background on special $\Gamma$-values}
\subsubsection{Notation}
Let ${\FF}_q$ be 
the
a field of $q$ elements.
where $q$ is a power of a prime $p$.  
Let $A := {\FF}_q[T]$ and $k :=
{\FF}_q(T)$, where $T$ is a variable. Let $A_{+} \subset A$ be the subset
of monic polynomials. Let
$|\cdot|_\infty$ be the unique valuation of $k$ for
which
$|T|_\infty = q$.
Let $k_\infty := {\FF}_q(\!(1/T)\!)$ be the $|\cdot|_\infty$-completion
of
$k$, let $\overline{k_\infty}$ be an algebraic closure of $k_\infty$,
let
$\CC_\infty$ be the
$|\cdot|_\infty$-completion of $\overline{k_\infty}$,
and let $\kbar$ be the algebraic closure of $k$ in $\CC_\infty$.

\subsubsection{The geometric $\Gamma$-function}
In \cite{ThakurGamma}, Thakur studied
the geometric $\Gamma$-function 
over $\FF_q[T]$,
\[
\Gamma(z) := \frac{1}{z} \prod_{n \in A_{+}} \left( 1 +
  \frac{z}{n} \right)^{-1} \quad (z \in \CC_\infty),
\]
which is a meromorphic function on $\CC_\infty$.
 Notably,
it satisfies several natural functional equations, which are the analogues
of the 
translation, 
reflection and Gauss multiplication
identities satisfied by the classical Euler $\Gamma$-function, and to 
which we refer as the
\emph{standard functional equations} (see
\S\ref{subsubsection:KeyPi}).

\subsubsection{Special $\Gamma$-values and the fundamental period of the
Carlitz module} We define the set of {\em special
$\Gamma$-values} to be
\[
\{\Gamma(z) \mid z \in k \setminus (-A_+ \cup \{ 0\}) \}\subset
k_\infty^\times.
\]
Up to factors in $k^\times$ a special
$\Gamma$-value
$\Gamma(z)$ depends only on $z$ modulo $A$.
In connection with special $\Gamma$-values it is natural also
to consider the number
\[
\varpi := T \sqrt[q-1]{-T} \prod_{i=1}^{\infty} \left( 1 - T^{1-q^i}
  \right)^{-1} \quad \in k_\infty\left(\sqrt[q-1]{-T}\right)^\times
\]
where $\sqrt[q-1]{-T}$ is a fixed $(q-1)^{st}$ root of $-T$
in $\CC_\infty$. 
The number $\varpi$ is the fundamental period of the Carlitz
module (see \S\ref{subsection:CarlitzExponential}) and hence
deserves to be regarded as the $\FF_q[T]$-analogue
of
$2\pi i$.
The number $\varpi$  is transcendental over
$k$ by \cite{Wade}. Our goal in this paper to determine all Laurent
polynomial relations with
coefficients in
$\kbar$ among special $\Gamma$-values and $\varpi$.

\subsubsection{Transcendence of special $\Gamma$-values}
For all $z\in A$ the value $\Gamma(z)$ belongs
to
$k
\cup
\{
\infty\}$.  However, it is known that for all $z \in k\setminus A$  the
value
$\Gamma(z)$ is transcendental over $k$. A short history of this
transcendence result is as follows. Isolated results on the transcendence
of special
$\Gamma$-values were obtained in \cite{ThakurGamma}.  The first
transcendence results for general classes of values of the
$\Gamma$-function were obtained in 
\cite{Sinha}.  Namely, Sinha showed that the values
$\Gamma(\frac{a}{f} + b)$ are transcendental over $k$
for all $a$, $f \in A_{+}$ and $b\in A$  such that $\deg a < \deg
f$. Sinha's results were obtained by representing the
$\Gamma$-values in question as periods of $t$-modules defined over
$\kbar$ and then invoking a transcendence criterion of Gelfond-Schneider
type from
\cite{YuMuchEarlier}.  
Subsequently all the values
$\Gamma(z)$ for $z \in k\setminus
A$  were represented in \cite{BrownPapa} as periods of $t$-modules
defined over $\kbar$ and thus proved transcendental.

\subsubsection{$\Gamma$-monomials and the diamond bracket criterion}
An element of the subgroup of $\CC_\infty^\times$ generated
by $\varpi$ and the special $\Gamma$-values will for brevity's sake be
called a {\em
$\Gamma$-monomial}. 
By adapting the Koblitz-Ogus
criterion 
given in \cite{DeligneCorvallis} to
the function field 
setting
along lines suggested in \cite{ThakurGamma},
we have at our disposal a {\em diamond bracket criterion}
(see Corollary~\ref{Corollary:KoblitzOgus}) capable of deciding in a
mechanical way whether between a given pair of
$\Gamma$-monomials there exists a $\kbar$-linear relation explained by the
standard functional equations.  We call the two-term $\kbar$-linear
dependencies thus arising {\em diamond bracket relations}.

\subsubsection{Cautionary example}
\label{subsubsection:DangerousBend}
A $\kbar$-linear relation between $\Gamma$-monomials
can be explained by the standard functional equations without 
strictly speaking being a
consequence of them. Consider the following
example concerning the classical
$\Gamma$-function taken from
\cite{Das}. The relation 
$$\frac{\Gamma\left(\frac{4}{15}\right)
\Gamma\left(\frac{1}{5}\right)}{\Gamma\left(\frac{1}{3}\right)
\Gamma\left(\frac{2}{15}\right)}
=
\sqrt{\left\{
\begin{array}{l}
\;\;\;
\frac{\Gamma\left(\frac{1}{5}\right)}
{\Gamma\left(\frac{1}{15}\right)
\Gamma\left(\frac{2}{5}\right)\Gamma\left(\frac{11}{15}\right)}\times
\frac{\Gamma\left(\frac{2}{5}\right)}
{\Gamma\left(\frac{2}{15}\right)
\Gamma\left(\frac{4}{5}\right)\Gamma\left(\frac{7}{15}\right)}\\\\
\times\frac{\Gamma\left(\frac{1}{15}\right)
\Gamma\left(\frac{4}{15}\right)\Gamma\left(\frac{7}{15}\right)
\Gamma\left(\frac{2}{3}\right)\Gamma\left(\frac{13}{15}\right)}
{\Gamma\left(\frac{1}{3}\right)}\\\\
\times
\frac{\Gamma\left(\frac{4}{15}\right)\Gamma\left(\frac{11}{15}\right)
}
{\Gamma\left(\frac{1}{3}\right)\Gamma\left(\frac{2}{3}\right)} \times
\frac{
\Gamma\left(\frac{1}{5}\right)\Gamma\left(\frac{4}{5}\right)}{\Gamma\left(\frac{2}{15}\right)\Gamma\left(\frac{13}{15}\right)}
\end{array}\right.}=3^{-\frac{1}{5}}5^{\frac{1}{12}}\sqrt{\frac{\sin\frac{\pi}{3}
\cdot\sin\frac{2\pi}{15}}{\sin\frac{4\pi}{15}\cdot
\sin\frac{\pi}{5}}}$$ 
confirms the Koblitz-Ogus criterion but decomposes
into instances of the standard functional
equations only after  the terms are squared.
The results of \cite{Kubert2}
imply the existence of such peculiar examples  in great
abundance. See
\cite{Das} for a method by which essentially all such examples
can be constructed explicitly. The analogous phenomena occur in
the function field situation. For a discussion of the latter, see
\cite{BaeEtAl}.
For a simple example in the case
$q=3$, which was in fact  discovered before all the others mentioned in
this paragraph, see
\cite[\S 4]{SinhaDR}.

\subsubsection{Linear independence}
It
was shown in \cite{BrownPapa} that the only relations
of
$\kbar$-linear dependence among $1$, $\varpi$, 
and special $\Gamma$-values are those following
from the diamond bracket relations. 
This result was obtained by carefully analyzing $t$-submodule
structures and then invoking Yu's powerful Theorem of the
$t$-Submodule \cite{Yu}.

\subsection{The main result}
We prove:
\begin{Theorem}[cf.\ Theorem
\ref{Theorem:LinearRelations}]\label{Theorem:LinearRelationsIntro}
 A set of $\Gamma$-monomials is
$\kbar$-linearly
  dependent exactly when some pair of $\Gamma$-monomials is.  
Pairwise  $\kbar$-linear (in)dependence of
$\Gamma$-monomials is entirely
decided
  by the diamond bracket criterion.
\end{Theorem}
In other words, all $\kbar$-linear relations among
$\Gamma$-monomials follow from the diamond bracket relations. 
The theorem has the following implication concerning transcendence
degrees:
\begin{Corollary}[cf.\ Corollary
\ref{Corollary:TranscendenceDegree}]\label{Corollary:TranscendenceDegreeIntro}
For all
$f\in A_+$ of positive degree, the extension of $\kbar$ generated
by the set
$\{\varpi\}\cup\left\{\Gamma(x)\left| x\in \frac{1}{f}A\setminus 
(\{0\}\cup-A_+)
\right.\right\}$
is of transcendence degree
$1+\frac{q-2}{q-1}\cdot\card(A/f)^\times$
over $\kbar$.
\end{Corollary}
In fact the corollary is equivalent to the theorem (see
Proposition~\ref{Proposition:Strategy}).

\subsection{Methods}
We outline the proof of Theorem~\ref{Theorem:LinearRelationsIntro},
emphasizing the new methods introduced here,
and compare our techniques to
those used previously.

\subsubsection{A new
linear
independence criterion}
We develop a new method for detecting $\kbar$-linear independence
of sets of numbers in $\overline{k_\infty}$, culminating in a quite
easily stated criterion. Let
$t$ be a variable 
independent of $T$.
Given $f=\sum_{i=0}^\infty
a_it^i\in \CC_\infty[\![t]\!]$ and $n\in \ZZ$, put
$f^{(n)}:=\sum_{i=0}^\infty a_i^{q^n}t^i$ and extend the
operation $f\mapsto f^{(n)}$ entrywise to matrices.  
Let $\EE \subset
\kbar[\![t]\!]$
be the subring 
consisting
of power series $\sum_{i=0}^\infty a_it^i$
such that $[k_\infty\left(\{a_i\}_{i=0}^\infty\right):k_\infty]<\infty$
and $\lim_{i\rightarrow\infty}\sqrt[i]{|a_i|_\infty}=0$. 
We now state our criterion (Theorem~\ref{Theorem:TranscendenceCriterion}
 is the verbatim repetition; see also
Proposition~\ref{Proposition:Rationale}):
\begin{Theorem}
\label{Theorem:TranscendenceCriterionIntro}
Fix a matrix
$\Phi=\Phi(t)\in \Mat_{\ell\times \ell}(\kbar[t])$
such that $\det \Phi$ is a polynomial in $t$ vanishing (if at all)
only at $t=T$.
Fix a (column) vector
$\psi=\psi(t)\in
\Mat_{\ell\times 1}(\EE)$
satisfying the functional equation
$\psi^{(-1)}= \Phi\psi$.
Evaluate $\psi$ at $t=T$, thus obtaining a vector
$\psi(T)\in\Mat_{\ell\times 1}\left(\overline{k_\infty}\right)$.
For every (row) vector
$\rho\in \Mat_{1\times \ell}(\kbar)$
such that
$\rho\psi(T) =0$
there exists a (row) vector
$P=P(t)\in \Mat_{1\times \ell}(\kbar[t])$
such that
$P(T)=\rho$ and $P \psi  =0$.
\end{Theorem}
Thus, in the situation of this theorem, every
$\kbar$-linear relation 
among entries of the specialization $\psi(T)$ is explained by a
$\kbar[t]$-linear relation 
among entries of $\psi$ itself.

\subsubsection{Dual $t$-motives}
The category of dual $t$-motives 
(see \S\ref{subsection:DualtMotives})
provides a natural setting in which we
can apply Theorem~\ref{Theorem:TranscendenceCriterionIntro}.
Like $t$-motives in \cite{AndersonMotive}, dual $t$-motives are
modules of a certain sort over a certain skew polynomial
ring.
 From a formal
algebraic perspective dual $t$-motives differ very little
from $t$-motives, and consequently
most 
$t$-motive concepts carry over
naturally to the dual
$t$-motive setting. In particular,
the concept of rigid analytic triviality carries over (see
\S\ref{subsection:DualtMotives}).
Crucially, to give a rigid analytic trivialization 
of a dual
$t$-motive 
is to give a square matrix with 
columns 
usable as input 
to Theorem~\ref{Theorem:TranscendenceCriterionIntro}
(see
Lemma \ref{Lemma:RATEquation}).

\subsubsection{Position of the new linear independence criterion with
respect to Yu's theorem of the
$t$-Submodule} We came upon
Theorem~\ref{Theorem:TranscendenceCriterionIntro} in the process of
searching for a
$t$-motivic translation of Yu's Theorem of the $t$-Submodule
\cite{Yu}.
Our discovery of a direct proof of
Theorem~\ref{Theorem:TranscendenceCriterionIntro} was a happy accident,
but it was one for which we were psychologically prepared by close study
of the proof of Yu's theorem. 

Roughly speaking, the points of view adopted in
the two theorems correspond as follows.
If $H=\Hom({\mathbb G}_a,E)$
is the dual
$t$-motive defined over $\kbar$ corresponding canonically to a
uniformizable abelian
$t$-module
$E$ defined over $\kbar$, and
$\Psi$ is a matrix describing a
rigid analytic trivialization of
$H$ as in Lemma~\ref{Lemma:RATEquation}, then it is possible
to express the periods of
$E$ in a natural way as
$\kbar$-linear combinations of entries of 
$\left(\Psi\vert_{t=T}\right)^{-1}$ and {\em vice
versa}. Thus it becomes at least plausible that
Theorem~\ref{Theorem:TranscendenceCriterionIntro}
and Yu's theorem provide similar information about $\kbar$-linear
independence.
A detailed comparison of
the two  theorems is not going to be presented here;
indeed, such has yet to be worked out. But
we are inclined to believe that at the end of the day the theorems differ
insignificantly in terms of  ability to detect
$\kbar$-linear independence.  

In any case, it is clear that both
theorems are  strong enough to handle the analysis of
$\kbar$-linear relations among $\Gamma$-monomials. 
Ultimately Theorem~\ref{Theorem:TranscendenceCriterionIntro} is our
tool of choice just because it is 
the easier to apply. Theorem~\ref{Theorem:TranscendenceCriterionIntro}
allows us to carry out our analysis entirely within the category  of dual
$t$-motives, which means that we can exclude $t$-modules
from the picture altogether at a considerable savings of labor in
comparison to
\cite{Sinha} and
\cite{BrownPapa}. 

\subsubsection{Linking $\Gamma$-monomials to dual
$t$-motives via Coleman functions}  In order
to generalize beautiful
examples in \cite{Coleman} and
\cite{ThakurGamma}, {\em solitons} over $\FF_q[T]$
were defined  and studied in \cite{AndersonStick}. In turn, 
in order to obtain various results on transcendence and
algebraicity of special $\Gamma$-values, variants of solitons called {\em
Coleman functions} were defined
 and studied in \cite{Sinha} and \cite{SinhaDR}.

We present in this paper a self-contained
elementary approach to Coleman functions producing new simple explicit
formulas  for them (see
\S\ref{section:Ragbag}, 
\S\ref{subsection:ColemanConstruction}). 
From the Coleman functions we then construct dual $t$-motives with
rigid analytic trivializations described by matrices with entries
specializing at
$t=T$ to $\kbar$-linear combinations of
$\Gamma$-monomials (see \S\ref{subsection:ColemanGCM}), thus putting
ourselves  in a position where
Theorem~\ref{Theorem:TranscendenceCriterionIntro} is 
at least potentially 
 applicable. 

Our method for attaching dual
$t$-motives to Coleman functions is straightforwardly adapted
from \cite{Sinha}. 
But our method for obtaining rigid
analytic trivializations is more elementary than that of
\cite{Sinha} because
the explicit formulas for Coleman functions at our disposal
obviate sophisticated apparatus from rigid analysis.

\subsubsection{Geometric complex multiplication and the end of the proof}
The dual $t$-motives engendered by Coleman functions are equipped with
extra endomorphisms and are examples of dual $t$-motives with
geometric complex multiplication, GCM for short (see
\S\ref{subsection:GCM}). We  
extend a technique developed in
\cite{BrownPapa} for analyzing $t$-modules with complex multiplication to
the setting of dual
$t$-motives with GCM, dubbing the generalized technique the {\em
Dedekind-Wedderburn trick} (see
\S\ref{subsection:DedekindWedderburnTrick}).
We determine that rigid analytically trivial dual
$t$-motives with GCM are semi-simple up to isogeny. In fact each
such object is
isogenous to a power of a simple dual $t$-motive.

\subsubsection{End of the proof}
Combining our general results on the structure of GCM dual $t$-motives
with our concrete results on the structure of the dual $t$-motives 
engendered by Coleman functions, we can finally  apply 
Theorem~\ref{Theorem:TranscendenceCriterionIntro} (in the guise
of Proposition~\ref{Proposition:Rationale}) to rule out  all
$\kbar$-linear relations among
$\Gamma$-monomials not following from the diamond bracket relations (see
\S\ref{subsection:LooseEnds}),
thus proving Theorem~\ref{Theorem:LinearRelationsIntro}.

\subsection{Comments on the classical case}
\label{subsubsection:ClassicalConjectures}
In the classical situation various people
have formed a clear picture about what algebraic relations should hold
among special $\Gamma$-values.
Those ideas stimulated our interest 
and guided our intuition
in
the function field setting. We discuss these ideas in more detail below.
\subsubsection{Temporary notation and terminology} Let
$\Gamma(s)$ be the classical
$\Gamma$-function, call 
$\{\Gamma(s)\mid s\in \QQ\setminus\ZZ_{\leq 0}\}$
the set of {\em special $\Gamma$-values}, 
and let a {\em
$\Gamma$-monomial} be any element of the subgroup of $\CC^\times$
generated by the special $\Gamma$-values and $2\pi i$.
\subsubsection{Rohrlich's
conjecture}\label{subsubsection:RohrlichConjecture} Rohrlich in the late
1970's made a conjecture which in rough form can be stated thus:
all multiplicative algebraic relations among
special
$\Gamma$-values and
$2\pi i$ are explained by the standard functional equations. See 
\cite[Appendix to \S 2, p.~66]{LangCyclotomic}
for a more precise formulation of the conjecture
in the language of distributions.  In language very similar to that
we have used above, Rohrlich's conjecture can be formulated 
as the assertion that the Koblitz-Ogus criterion for a
$\Gamma$-monomial to belong to $\overline{\QQ}^\times$ is not only
sufficient, but necessary as well.
\subsubsection{Lang's conjecture}
\label{subsubsection:LangConjecture}
Lang subsequently
strengthened Rohrlich's conjecture to a conjecture which in rough form
can be stated thus: all polynomial algebraic relations among special
$\Gamma$-values and $2\pi i$ with coefficients in $\overline{\QQ}$ are
explained by the standard functional equations. 
See \cite[loc.\ cit.]{LangCyclotomic} for a formulation of this conjecture
in the language of distributions. In language very similar
to that we have used above, Lang's conjecture can be formulated
as the assertion that all 
$\overline{\QQ}$-linear relations among
$\Gamma$-monomials follow from the
two-term relations provided by the Koblitz-Ogus criterion. 
Yet another formulation of Lang's conjecture 
is the assertion that for every integer $n>2$ the transcendence degree
of the extension of $\overline{\QQ}$ generated by the set
$\{2\pi i\}\cup\left\{\Gamma(x)\left| x\in
\frac{1}{n}\ZZ\setminus\ZZ_{\leq 0}\right.\right\}$
is equal to 
$1+\phi(n)/2$, where $\phi(n)$ is Euler's totient.
What we prove in this paper is the analogue 
of Lang's conjecture for the geometric $\Gamma$-function over $\FF_q[T]$.
\subsubsection{Evidence in the classical case}
There are very few integers $n>1$
such that all Laurent polynomial relations among elements of the set 
$\{2\pi
i\}\cup\left\{\Gamma\left(\frac{1}{n}\right),\dots,
\Gamma\left(\frac{n-1}{n}\right)\right\}$
with coefficients in $\overline{\QQ}$
can be ruled out save those following from the two-term relations
provided by the
Koblitz-Ogus criterion,
to wit:
\begin{itemize}
\item 
$n=2$ (Lindemann 1882, since $\Gamma(1/2) =
  \sqrt{\pi}$)
\item 
$n=3,4$  (Chudnovsky 1974, cf.\ \cite{Wald})
\end{itemize}
The only other
evidence 
known
for Lang's conjecture is 
indirect, and it is contained in a result of
\cite{WolWust}: all
$\overline{\QQ}$-linear relations 
among the special Beta values
$$B(a,b)=
\frac{\Gamma(a)\Gamma(b)}{\Gamma(a+b)}\;\;(a,b\in \QQ,\;\;\;
a,b,a+b\not\in \ZZ_{\leq 0})$$ and $2\pi
i$ follow from the two-term relations provided by the Koblitz-Ogus
criterion. 

\subsection{Acknowledgements}
Research of the second author was supported in part by an NSF grant.
The second  
and third authors would like to thank the Erwin Schr\"odinger
Institute for 
its hospitality during some of the final editorial
work.

\section{Notation and terminology}
\label{section:NotationAndTerminology}

\subsection{Table of special symbols}
\label{subsection:SymbolTable}
$$
\begin{array}{lcl}
T,t,z&:=&\mbox{independent
variables}\\
\FF_q&:=&\mbox{a field of
$q$ elements}\\
k&:=&\FF_q(T)\\
|\cdot|_\infty&:=&\mbox{the unique valuation of $k$ such
that $|T|_\infty=q$}\\
k_\infty&:=&\FF_q(\!(1/T)\!)=\mbox{the $|\cdot|_\infty$-completion of
$k$}\\
\overline{k_\infty}&:=&\mbox{an algebraic closure of $k_\infty$}\\
\CC_\infty&:=&\mbox{the $|\cdot|_\infty$-completion of
$\overline{k_\infty}$}\\
\kbar&:=&\mbox{the algebraic closure of $k$ in
$\CC_\infty$}\\
\Tbar&:=&\mbox{a fixed choice in $\kbar$ of a $(q-1)^{st}$ root of $-T$}\\
\CC_\infty\{t\}&:=&\mbox{the subring of the power series
ring $\CC_\infty[\![t]\!]$ consisting of  }\\
&&\mbox{power series convergent in the ``closed'' unit disc
$|t|_\infty\leq 1$}\\
\card S&:=&\mbox{the cardinality of a set $S$}\\
\Mat_{r\times s}(R)&:=&\mbox{the set of $r$ by $s$ matrices with entries
in a ring or module $R$}\\
R^\times&:=&\mbox{the group of units of a ring $R$ with unit}\\
\GL_n(R)&:=&\Mat_{n\times n}(R)^\times,\;\mbox{where $R$ is a ring
with unit}\\
\one_n&:=&\mbox{the $n$ by $n$ identity matrix}\\
 A&:=&\FF_q[T]\\
\deg&:=&\mbox{the function associating to each element of
$A$ its degree in $T$}\\
A_+&:=&\mbox{the set of elements of $A$
monic in
$T$}\\
D_N&:=&\prod_{i=0}^{N-1}(T^{q^N}-T^{q^i})\;\in A_+\\
\Res&:=&\left(\sum_i a_iT^i\mapsto
a_{-1}\right):k_\infty\rightarrow \FF_q\\
\end{array}
$$

\subsection{Twisting}
\label{subsection:Twisting} Fix $n\in \ZZ$.
Given a formal power series 
$f=\sum_{i=0}^\infty a_it^i\in \CC_\infty[\![t]\!]$
we define the {\em $n$-fold twist}
by the rule
$f^{(n)}:=\sum_{i=0}^\infty a_i^{q^n}t^i$.
The $n$-fold twisting operation is an automorphism of the power series
ring $\CC_\infty[\![t]\!]$ stabilizing various subrings, e.~g.,
$\kbar[\![t]\!]$, $\kbar[t]$,  and $\CC_\infty\{t\}$.
More
generally, for any matrix $F$ with entries in $\CC_\infty[\![t]\!]$
we define the $n$-fold twist $F^{(n)}$ by the rule
$\left(F^{(n)}\right)_{ij}:=(F_{ij})^{(n)}$.
In particular, for any matrix $X$ with entries in $\CC_\infty$ we
have $\left(X^{(n)}\right)_{ij}=(X_{ij})^{q^n}$. The
$n$-fold twisting operation commutes with matrix addition and
multiplication.

\subsection{Norms}
For any matrix $X$ with entries in $\CC_\infty$ we put
$|X|_\infty:=\max_{ij}|X_{ij}|_\infty$.
We have
$\left|X^{(n)}\right|_\infty=|X|_\infty^{q^n}$
for all $n\in \ZZ$.
We have 
$$|U+V|_\infty\leq \max(|U|_\infty,|V|_\infty),\;\;\;|XY|_\infty\leq
|X|_\infty\cdot|Y|_\infty$$ for all matrices $U$, $V$, $X$, $Y$ with
entries in $\CC_\infty$ such that $U+V$ and $XY$ are defined.

\subsection{The ring $\EE$}
\label{subsection:RingEE}
We define $\EE$ to be the ring consisting of formal power series
$$\sum_{n=0}^\infty a_nt^n\in \kbar[\![t]\!]$$
such that
$$\lim_{n\rightarrow\infty} \sqrt[n]{|a_n|_\infty}=0,\;\;\;\;
[k_\infty(a_0,a_1,a_2,\dots):k_\infty]<\infty.$$
The former condition guarantees that such a series has an infinite radius
of convergence with respect to the valuation
$|\cdot|_\infty$. The latter condition guarantees that for any
$t_0\in\overline{k_\infty}$ the  value of such a series at $t=t_0$
belongs again to
$\overline{k_\infty}$. Note that the ring
$\EE$ is stable under  the $n$-fold twisting operation
$f\mapsto f^{(n)}$ for all $n\in \ZZ$. 

\subsection{The Schwarz-Jensen formula}
Fix 
$f\in
\EE$ not vanishing identically. It is possible to enumerate the
zeroes of
$f$ in $\CC_\infty$ because there are only finitely many zeroes
in each disc of finite radius.
Put
$$\{\omega_i\}:=\mbox{an enumeration (with multiplicity)
of the zeroes of
$f$ in
$\CC_\infty$}$$
and
$$\lambda:=\mbox{the leading coefficient of the Maclaurin expansion
of $f$.}
$$
The {\em
Schwarz-Jensen} formula \\
$$\sup_{\begin{subarray}{l}
x\in \CC_\infty\\
|x|\leq r
\end{subarray}}|f(x)|_\infty=
|\lambda|_\infty \;\cdot\;r^{\card\{i\left|\omega_i=0\}\right.}
\;\cdot\;
\prod_{i:\;0<|\omega_i|_\infty<r}\frac{r}{|\omega_i|_\infty}
\;\;\;\;\;\;(r\in
\RR_{>0})$$ \\
relates the growth  of the modulus of $f$ to the
distribution of the zeroes of $f$. This fact is an easily deduced
corollary to the Weierstrass Preparation Theorem over a complete discrete
valuation ring.

\section{A linear independence
criterion}\label{section:TranscendenceCriterion}

\subsection{Formulation and discussion of the
criterion}
\begin{Theorem}\label{Theorem:TranscendenceCriterion}
Fix a matrix
$$\Phi=\Phi(t)\in \Mat_{\ell\times \ell}(\kbar[t]),$$
such that $\det \Phi$ is a polynomial in $t$ vanishing (if at all)
only at $t=T$.
Fix a (column) vector
$$\psi=\psi(t)\in
\Mat_{\ell\times 1}(\EE)$$
satisfying the functional equation
$$\psi^{(-1)}= \Phi\psi.$$
Evaluate $\psi$ at $t=T$, thus obtaining a vector
$$\psi(T)\in\Mat_{\ell\times 1}\left(\overline{k_\infty}\right).$$
For every (row) vector
$$\rho\in \Mat_{1\times \ell}(\kbar)$$
such that
$$\rho\psi(T) =0$$
there exists a (row) vector
$$P=P(t)\in \Mat_{1\times \ell}(\kbar[t])$$
such that
$$P(T)=\rho,\;\;\;P \psi  =0.$$
\end{Theorem}
The proof commences in \S\ref{section:SiegelLemmas} and
takes up the rest of
\S\ref{section:TranscendenceCriterion}.
We think of the
$\kbar[t]$-linear relation $P$ among the
entries of $\psi$ produced by the theorem as an ``explanation'' or
a ``lifting'' of the given $\kbar$-linear relation $\rho$ among the entries
of
$\psi(T)$.

\subsubsection{The basic application}
\label{subsubsection:TranscendenceIllustration}
Consider the power series
$$\Omega=\Omega(t):=\tilde{T}^{-q}\prod_{i=1}^\infty
\left(1-t/T^{q^i}\right)\in
k_\infty(\Tbar)[\![t]\!]\subset\CC_\infty[\![t]\!].$$
The power series
$\Omega(t)$ has an infinite radius of convergence and satisfies the
functional equation
$$\Omega^{(-1)}=(t-T)\cdot\Omega.$$
Consider the Maclaurin expansion
$$\Omega(t)=\sum_{i=0}^\infty a_it^i.$$
The functional equation satisfied by $\Omega$ implies the recursion
$$\sqrt[q]{a_{i}}+Ta_{i}=\left\{\begin{array}{ll}
a_{i-1}&\mbox{if $i>0$,}\\
0&\mbox{if $i=0$.}
\end{array}\right.
$$
Therefore $\Omega$ belongs to 
$\kbar[\![t]\!]$ and hence to $\EE$.
Suppose now that there exists a nontrivial $\kbar$-linear relation
$$\sum_{i=0}^n
\rho_i\Omega(T)^i=0\;\;\;(\rho_i\in \kbar,\;\;\;n>0,\;\;\;\rho_0\rho_n\neq
0)$$ among the powers of the number 
$$\Omega(T)=\tilde{T}^{-q}\prod_{i=1}^\infty
(1-T^{1-q^i})\in k_\infty(\Tbar).$$
Theorem~\ref{Theorem:TranscendenceCriterion} provides a
$\kbar[t]$-linear ``explanation''
$$\sum_{i=0}^n P_i\Omega^i=0\;\;\;(P_i\in
\kbar[t],\;\;P_i(T)=\rho_i).$$
But  
the polynomial $P_0$ must vanish at all the zeroes $t = T^{q^i}$ of the
function $\Omega$.  Thus $P_0$ vanishes identically, contrary to our
assumption that $\rho_0 = P_0(T) \ne 0$.  We conclude that $\Omega(T)$
is transcendental over $k$.

See \S\ref{subsection:CarlitzExponential} below for the interpretation
of $-1/\Omega(T)$ as the fundamental period of the Carlitz module. The
power series $\Omega(t)$ plays a key role in this paper.

\begin{Proposition}\label{Proposition:HypothesisChecker}
Suppose we are given 
$$\Phi\in \Mat_{\ell\times \ell}(\kbar[t]),\;\;\;
\psi\in
\Mat_{\ell\times 1}\left(\CC_\infty\{t\}\right)$$
such that 
$$\det \Phi(0)\neq 0,\;\;\;\psi^{(-1)}=\Phi\psi.$$
Then we necessarily have
$$\psi\in \Mat_{\ell\times 1}(\EE).$$
\end{Proposition}
The proposition simplifies the task of checking the
hypotheses of Theorem~\ref{Theorem:TranscendenceCriterion}.
\proof  
Write
$$\Phi=\sum_{i=0}^N
b_{(i)}t^i\;\;\;(b_{(i)}\in
\Mat_{\ell\times \ell}(\kbar),\;N:\;\mbox{positive integer}).$$
By hypothesis 
$b_{(0)}\in \GL_\ell(\kbar)$. 
By the theory
of Lang isogenies
\cite{LangAlgebraicGroups} there exists 
$U\in \GL_{\ell\times \ell}(\kbar)$
such that 
$$U^{(-1)}b_{(0)}U^{-1}=\one_\ell\;\;\;\left(\mbox{equivalently:}\;\;
b_{(0)}^{(1)}=U^{-1}U^{(1)}\right).$$
After making the replacements
$$\psi\leftarrow U\psi,\;\;\;\;\Phi\leftarrow U^{(-1)}\Phi
U^{-1},
$$ 
we may assume without loss of generality that
$b_{(0)}=\one_\ell$.
Now write
$$\psi=\sum_{i=0}^\infty a_{(i)}t^i\;\;\;(a_{(i)}\in \Mat_{\ell\times
1}(\CC_\infty)).$$
We have
$$a_{(n)}^{(-1)}-a_{(n)}=\sum_{i=1}^{\min(n,N)}
b_{(i)}a_{(n-i)},$$ and hence 
$$a_{(n)}\in \Mat_{\ell\times 1}(\kbar)$$
for all integers $n\geq 0$. By hypothesis   
$$\lim_{n\rightarrow\infty}|a_{(n)}|_\infty=0,$$
and hence the series 
$$\tilde{a}_{(n)}:=\sum_{\nu=1}^\infty
\left(\sum_{i=1}^{N}
b_{(i)}a_{(n-i)}\right)^{(\nu)}$$
converges
for all $n \gg 0$. Moreover, 
$$\lim_{n\rightarrow\infty}|\tilde{a}_{(n)}|_\infty=0.$$ 
Since we have
$$\left(\tilde{a}_{(n)}-a_{(n)}\right)^{(-1)}=
\left(\tilde{a}_{(n)}-a_{(n)}\right)$$
for $n\gg 0$, 
it follows that 
$\tilde{a}_{(n)}=a_{(n)}$
for $n\gg 0$ and hence that
the collection of entries of all the vectors $a_{(n)}$ generates
an extension of
$k_\infty$ of finite degree.
Now fix $C>1$ arbitrarily.
{From} the fact that $\tilde{a}_{(n)} = a_{(n)}$ for $n \gg 0$, we have inequalities
$$\begin{array}{rcl}
C^{n}|a_{(n)}|_\infty&\leq&\displaystyle \max_{i=1}^{N}
C^{n}|b_{(i)}|_\infty^q
|a_{(n-i)}|_\infty^q\\\\
&\leq&\displaystyle \left(\max_{i=1}^{N}
C^{i}|b_{(i)}|_\infty^q\right)
\cdot \left(\max_{i=1}^N |a_{(n-i)}|_\infty\right)^{q-1}
\cdot\left(\max_{i=1}^N C^{n-i}|a_{(n-i)}|_\infty\right)\\\\
&\leq&\displaystyle
\max_{i=0}^{n-1}C^{i}|a_{(i)}|_\infty
\end{array}$$
for  $n\gg 0$, and hence 
$$\sup_{n=0}^\infty C^n|a_{(n)}|_\infty<\infty.$$ 
Therefore the
radius of convergence of each entry of
$\psi$ is infinite.
\qed

\subsubsection{Remark}
Theorem~\ref{Theorem:TranscendenceCriterion} is in essence the (dual)
$t$-motivic translation of Yu's Theorem of the
$t$-Submodule
\cite[Thms.~3.3 and 3.4]{Yu}.
Once the setting is sufficiently developed, we expect that the class
of numbers about which Theorem~\ref{Theorem:TranscendenceCriterion}
provides $\kbar$-linear independence information is essentially the
same as that handled by Yu's theorem of the $t$-Submodule, and the
type of information provided is essentially the same, too.  We omit
discussion of the comparison.

\subsection{Specialized notation for making estimates}
\subsubsection{Degree in $t$}
Given a polynomial $f\in \kbar[t]$ let 
$\deg_t f$ denote its degree in $t$,
and, more generally, given a matrix
$F$ with entries in $\kbar[t]$
put 
$\deg_t F:=\max_{ij}\deg_t F_{ij}$.
We have 
$\deg_t F^{(n)}=\deg_t F$ for all $n\in \ZZ$.
We have 
$$\deg_t (D+E)\leq \max\left(\deg_t D,\deg_t E\right),\;\;\;\deg_t
(FG)\leq
\deg_t F+\deg_t  G$$ for all matrices $D$, $E$, $F$, $G$ with entries in
$\kbar[t]$ such that $D+E$ and $FG$ are defined.

\subsubsection{Size}
Given an algebraic number $x\in \kbar$ we set 
$\size{x}:=\max_{\tau}|\tau x|_\infty$,
where $\tau$ ranges over the automorphisms of $\kbar/k$, thereby
defining the {\em size} of $x$. More generally given a polynomial
$f=\sum_i a_it^i\in \kbar[t]$,
we define 
$\size{f}:=\max_i \size{a_i}$. Yet more generally, given a matrix
$F$ with entries in $\kbar[t]$
we define  $\size{F}:=\max_{ij}\size{F_{ij}}$. Then we have 
$\size{F^{(n)}}=\size{F}^{q^n}$
for all $n\in \ZZ$.
We have 
$$\size{D+E}\leq \max(\size{D},\size{E}),\;\;\;
\size{FG}\leq
\size{F}\cdot\size{G}$$ for all matrices $D$, $E$, $F$, $G$ with entries in
$\kbar[t]$ such that $D+E$ and $FG$ are defined.

\subsection{The basic estimates}\label{section:SiegelLemmas}
\subsubsection{The setting}
Throughout \S\ref{section:SiegelLemmas} we fix fields
$$k\subset K_0\subset K\subset\kbar$$
and rings
$$A\subset\OO_0\subset\OO\subset K$$
such that 
\begin{itemize}
\item $K_0/k$ is a finite separable extension,
\item $K$ is the closure of $K_0$ in $\kbar$  under the extraction of
$q^{th}$ roots, 
\item $\OO$ is the integral closure of $A$ in $K$, and
\item $\OO_0=\OO\cap K_0$.
\end{itemize} 
Note that $\kbar$ is the
union of all its subfields of the form $K$. 
\subsubsection{Lower bound from size}\label{subsubsection:LowerBound}
We claim that
$$\size{x}\geq 1,\;\;\;|x|_\infty\geq \size{x}^{1-[K_0:k]}$$
for all $0\neq x\in \OO$. Clearly these estimates hold in the case
$0\neq x\in
\OO_0$, because in that case $x$ has at most $[K_0:k]$ conjugates over
$k$ and the product of those conjugates is a nonzero element of
$A$; but then, since we have
$$\OO=\bigcup_{\nu=0}^\infty \OO_0^{q^{-\nu}},$$
the claim holds in general.

\begin{Lemma}[Liouville Inequality]
\label{Lemma:Clincher}
Fix a polynomial
$$f(z):=\sum_{i=0}^n a_iz^i\in \OO[z]$$
not vanishing identically.
For every nonzero root $\lambda\in \kbar$ of $f(z)$ of order $\nu$ we
have
\[
|\lambda|_\infty^{\nu}
\geq \left(\max_{i=0}^n\size{a_i}\right)^{-[K_0:k]}.
\]
\end{Lemma}
\proof We may of course assume that 
$|\lambda |_\infty < 1$, for otherwise the
claim is obvious.
After factoring out a power of $z$ we may also assume that $a_0\neq
0$. 
Write
$$f(z+\lambda)=
\sum_{i=\nu}^n
b_i z^i\;\;\;(b_i\in \OO[\lambda]),$$
noting that
$$|b_i|_\infty\leq \max_{j=i}^n |a_j|_\infty\leq \max_{i=0}^n\size{a_i}.$$
Evaluate  the displayed expression for $f(z+\lambda)$ at $z =-\lambda$,
thus obtaining an estimate
\[
|a_0|_{\infty} = |f(0)|_{\infty} \le
\max_{i = \nu}^n |b_i\lambda^i|_{\infty} \le
|\lambda|_{\infty}^{\nu} \max_{i = 0}^{n} \size{a_i}.
\]
Finally, apply the
fundamental lower bound of \S\ref{subsubsection:LowerBound} to 
$a_0$.
\qed

\begin{Lemma}\label{Lemma:AsymptoticRiemannRoch}
For all
constants
$C>1$ we have
$$\lim_{\nu\rightarrow\infty}\left(\card\left\{\left.x\in
\OO_0^{q^{-\nu}}\right|
\size{x}\leq C\right\}\right)^{\frac{1}{q^\nu\cdot
[K_0:k]}}=C.$$
\end{Lemma}
The normalization $|T|_\infty=q$ was imposed to make this formula hold.
\proof
We may assume without loss of generality that $C$ is of the form 
$$C=q^\delta\;\;\;\;\left(\delta\in \bigcup_{\nu=0}^\infty
q^{-\nu}\ZZ,\;\;\;\delta>0\right).$$ The Riemann-Roch theorem yields
constants
$n_0$ and $n_1$ such that
$$n> n_0\Rightarrow \card\{x\in \OO_0\vert \size{x}\leq |T|_\infty^n\}
=q^{[K_0:k]n+n_1}$$
for all $n\in \ZZ$. 
We then have $$
\card\left\{\left.x\in \OO_0^{q^{-\nu}}\right|
\size{x}\leq |T|_\infty^\delta\right\}=
\card\left\{x\in \OO_0\left| \size{x}\leq
|T|_\infty^{q^\nu\delta}\right.\right\}=q^{[K_0:k]q^\nu
\delta+n_1}$$
for all integers $\nu\gg 0$, whence the result.
\qed

\begin{Lemma}[Thue-Siegel Analogue]\label{Lemma:SiegelLemmaOne}
Fix parameters
$$C>1,\;\;\;0<r<s\;\;\;(C\in \RR,\;\;\;r,s\in \ZZ).$$
For each matrix
$$M\in \Mat_{r\times s}(\OO)$$
such that
$$\size{M}<C$$ 
there exists
$$x\in \Mat_{s\times 1}(\OO)$$
such that
$$x\neq 0,\;\;\;Mx=0,\;\;\;\size{x}<C^{\frac{r}{s-r}}.
$$
\end{Lemma}
\proof 
Choose $C'>1$ and 
$\epsilon>0$ such that
$$\size{M}<
C',\;\;\;(1+\epsilon)(C')^{\frac{r}{s-r}}<C^{\frac{r}{s-r}}.$$
 For all
$\nu\gg 0$ the cardinality of the set
$$\left\{\left.x\in \Mat_{s\times 1}\left(\OO_0^{q^{-\nu}}\right)
\right| \size{x}\leq
(1+\epsilon)(C')^{\frac{r}{s-r}}\right\}$$ exceeds the
cardinality of the set
$$
\left\{\left.x\in \Mat_{r\times 1}\left(\OO_0^{q^{-\nu}}\right)
\right|
\size{x}\leq (1+\epsilon)(C')^{\frac{s}{s-r}}\right\}
$$
by Lemma~\ref{Lemma:AsymptoticRiemannRoch}. Further, for all $\nu\gg 0$
multiplication by $M$ maps the former set to the latter. Therefore the
desired vector $x$  exists by the pigeonhole principle. 
\qed
\begin{Lemma}\label{Lemma:SiegelLemmaTwo}
Again fix parameters
$$C>1,\;\;\;0<r<s\;\;\;(C\in \RR,\;\;r,s\in \ZZ).$$
For each matrix
$$M\in \Mat_{r\times s}(\OO[t])$$
such that
$$\size{M}<C$$ 
there exists
$$x\in \Mat_{s\times 1}(\OO[t])$$
such that
$$x\neq 0,\;\;\;Mx=0,\;\;\;\size{x}<C^{\frac{r}{s-r}}.$$ 
\end{Lemma}
\proof Let $d$ and $e$ be nonnegative integers presently to be chosen
efficaciously large and put
$$\;r':=r(d+e+1),\;\;\;s':=s(e+1).$$
Choose $d$ large enough so that
$$\deg_t  M\leq d,$$
and then choose $e$ large enough so that
$$r'<s',\;\;\;\size{M}<C':=
C^{\frac{r}{s-r}/\frac{r'}{s'-r'}}.$$
Consider now the $\OO$-linear map  
$$\{x\in \Mat_{s\times 1}(\OO[t])\vert \deg_t  x\leq
e\}\rightarrow\{x\in
\Mat_{r\times 1}(\OO[t])\vert \deg_t  x\leq d+e\}$$ induced by
multiplication by $M$. With respect to the evident choice of bases
the map under consideration is represented by a matrix
$$M'\in \Mat_{r'\times s'}(\OO)$$
such that
$$\size{M'}<C'.$$ The existence of
$x\in \Mat_{s\times 1}(\OO[t])$ such that
$$x\neq 0,\;\;\;Mx=0,\;\;\;\deg_t  x\leq e,\;\;\;
\size{x}<(C')^{\frac{r'}{s'-r'}}=C^{\frac{r}{s-r}}$$ now follows by an
application of the preceding lemma with the triple of parameters
$(C',r',s')$ in place of the triple
$(C,r,s)$.
\qed

\subsection{Proof of the criterion}
\subsubsection{The case $\ell=1$}
Assume for the moment that $\ell=1$. In this case we
may assume without loss of generality that $\rho\neq 0$
and hence that
$$\psi(T)=0,$$
in which case our task is to show that $\psi$ vanishes identically.
For any integer
$\nu\geq 0$ we have
$$\left(\psi\left(T^{q^{-\nu}}\right)\right)^{q^{-1}}
=\psi^{(-1)}\left(T^{q^{-(\nu+1)}}\right)=
\Phi\left(T^{q^{-(\nu+1)}}\right)\psi\left(T^{q^{-(\nu+1)}}\right),$$
$$
\Phi\left(T^{q^{-(\nu+1)}}\right)\neq 0$$
and hence we have
$$\psi\left(T^{q^{-\nu}}\right)=0\;\;\;(\nu=0,1,2,\dots).$$
Since $\psi$ vanishes infinitely many times in the disc
$|t|_\infty\leq|T|_\infty$, necessarily 
$\psi$ vanishes identically. Thus the case
$\ell=1$ of Theorem~\ref{Theorem:TranscendenceCriterion} is
dispatched. 

\subsubsection{Reductions and further notation}
Assume now that
$\ell>1$. We may of course assume that
$$\rho\neq 0.$$
As in \S\ref{section:SiegelLemmas} let
$$k\subset K_0\subset K\subset \kbar$$
be fields such that $K_0/k$ is a finite separable extension and
$K$ is the closure of $K_0$ under the extraction of $q^{th}$ roots. Since
$\kbar$ is the union of fields of the form $K$ we may assume without loss
of generality that
$$\Phi\in \Mat_{\ell\times \ell}(K[t]),\;\;\;
 \rho\in \Mat_{1\times \ell}(K).$$
As in \S\ref{section:SiegelLemmas} let
$\OO$ be the integral closure of $A$ in $K$. 
After making  replacements
$$\Phi\leftarrow a^{q-1}\Phi,\;\;\psi\leftarrow a^{-q}\psi,\;\;\;
\rho\leftarrow b\rho$$
for suitably chosen
$$ a,b\in A,\;\;\;ab\neq 0,$$
we may assume without loss of generality that
$$\Phi\in \Mat_{\ell\times \ell}(\OO[t]),\;\;\;\rho\in
\Mat_{1\times
\ell}(\OO).$$
 Fix a matrix
$$\vartheta\in \Mat_{\ell\times (\ell-1)}(\OO)$$
of maximal rank such that
$$\rho\vartheta =0.
$$
Then the $K$-subspace of $\Mat_{1\times \ell}(K)$ annihilated by
right multiplication by $\vartheta$ is the $K$-span of $\rho$.
 Let
$$\Theta\in \Mat_{\ell\times \ell}(\OO[t])$$
be the transpose of the matrix of cofactors of $\Phi$. Then we have
$$\Phi\Theta=\Theta\Phi=\det
\Phi\cdot\one_\ell=c(t-T)^s\cdot\one_\ell$$  for some $0\neq c\in
\OO$ and integer $s\geq 0$.  
Let $N$ be a parameter taking values in the set of positive integers
divisible by $2\ell$. 

\subsubsection{Construction of the auxiliary function $E$}
We claim there exists
$$h=h(t)\in \Mat_{1\times \ell}(\OO[t])$$
depending on the parameter $N$ such that\\
\begin{itemize}
\item $\size{h}=O(1)$ as $N\rightarrow\infty$\\
\end{itemize}
and with the following properties
for each value of $N$:\\
\begin{itemize}
\item $h\neq 0$.\\
\item $\deg_t  h<\left(1-\frac{1}{2\ell}\right)N$.\\
\item $E\left(T^{q^{-(N+\nu)}}\right)=0$ for
$\nu=0,\dots,N-1$, where $E:=h\psi\in \EE$.\\ (We call $E$ the {\em
auxiliary function}.)\\
\end{itemize}
Before proving the claim, we note first that the auxiliary
function
$E$ figures in the following identity:
$$\begin{array}{cl}
&\displaystyle h\Theta^{(-0)}\cdots\Theta^{(-(N+\nu-1))}\psi^{(-(N+\nu))}
\\\\
=&\displaystyle
h\Theta^{(-0)}\cdots\Theta^{(-(N+\nu-1))}\Phi^{(-(N+\nu-1))}
\cdots\Phi^{(-0)}\psi\\\\
=&\displaystyle
c^{q^{-(N+\nu-1)}+\cdots+q^0}\left(t-T^{q^{-(N+\nu-1)}}\right)^s\cdots
\left(t-T^{q^{-0}}\right)^s E
\end{array}
$$
This identity is useful 
again below and so for convenient
reference we dub it the {\em key identity}. By the key identity,
the hypothesis 
$$\rho\psi(T)=0\;\;\;\left(\mbox{equivalently:}\;\;\;\;
\rho^{(-(N+\nu))}\psi^{(-(N+\nu))}\left(T^{q^{-(N+\nu)}}\right)=0
\right),$$
and the definition of $\vartheta$,
the following condition forces the
desired vanishing of
$E$:
\\
\begin{itemize}
\item
$\left.
h\Theta^{(0)}\cdots\Theta^{(-(N+\nu-1))}\vartheta^{(-(N+\nu))}\right|_{t=T^{q^{-(N+\nu)}}}=0$
for
$\nu=0,\dots,N-1$.\\
\end{itemize}
Put
$$r:=\left(\ell-1\right)N,\;\;\;
s:=\left(\ell-\frac{1}{2}\right)N.$$
With respect to the evident choices of bases, the
homogeneous system of linear equations that we need to solve
is described
by a matrix 
$$M\in \Mat_{r\times s}(\OO)$$ depending on $N$ such that
$$\size{M}\leq |T|_\infty^{q^{-N}((1-\frac{1}{2\ell})N+2N\cdot \deg_t
\Theta)}\cdot
\size{\Theta}^{\frac{q}{q-1}}\cdot
\size{\vartheta}=O(1)\;\;\;\mbox{as}\;\;\;N\rightarrow\infty,$$  and
the solution we need to find is described by a vector
$$x\in \Mat_{s\times 1}(\OO)$$
depending on $N$ such that
$$x\neq
0,\;\;\;Mx=0,\;\;\;\size{x}=O(1)\;\;\;\mbox{as}\;\;\;N\rightarrow\infty.$$
Lemma~\ref{Lemma:SiegelLemmaOne} now proves our claim.

\subsubsection{A functional equation
for $E$}
We claim
there exist polynomials
$$a_0,\dots,a_\ell\in \OO[t]$$
depending on the parameter $N$ such that
\begin{itemize}
\item $\displaystyle\max_{i=0}^\ell\size{a_i}=O(1)$ as
$N\rightarrow\infty$\\
\end{itemize}
and with the following properties for each value of
$N$:\\
\begin{itemize}
\item Not all the $a_i$ vanish identically.\\
\item $a_0E+a_1E^{(-1)}+\dots+a_\ell E^{(-\ell)}=0$.\\
\end{itemize}
Since
$$E^{(-\nu)}=h^{(-\nu)}\Phi^{(-(\nu-1))}\cdots\Phi^{(-0)}\psi$$
for any integer $\nu\geq 0$, the functional equation
we want $E$ to satisfy is implied by the following condition:\\
\begin{itemize}
\item $a_0h^{(0)}+a_1h^{(-1)}\Phi^{(-0)}
+\cdots+a_\ell h^{(-\ell)}\Phi^{(-(\ell-1))}\cdots\Phi^{(-0)}=0$.\\
\end{itemize}
The latter system of homogeneous linear equations for $a_0,
\dots, a_{\ell}$ is with respect to the evident choice of bases
described by a matrix
$$M\in \Mat_{\ell\times(\ell+1)}(\OO[t])$$
depending on $N$ such that
$$\size{M}=O(1)\;\;\;\mbox{as}\;\;N\rightarrow\infty,$$
and the solution we have to find is described by a vector
$$x\in \Mat_{(\ell+1)\times 1}(\OO[t])$$
depending on $N$ such that
$$x\neq 0,\;\;\;Mx=0,\;\;\;
\size{x}=O(1)\;\;\;\mbox{as}\;\;N\rightarrow\infty.$$
Lemma~\ref{Lemma:SiegelLemmaTwo} now proves our  claim.
After dividing out common factors of $t$ we may further
assume that for each value of $N$:
\begin{itemize}
\item Not all the constant terms $a_i(0)$ vanish.
\end{itemize}

\subsubsection{Vanishing of $E$}
We claim that $E$ vanishes identically for some $N$.
Suppose that this is not the case.
Let $\lambda$ be the leading coefficient
of the Maclaurin expansion of
$E$. We have
$$a_0(0)\lambda^{q^0}+\dots+a_\ell(0)\lambda^{q^{-\ell}}=0,$$
and hence
$$1/|\lambda|_\infty=O(1) \ \text{as $N\rightarrow\infty$}$$
by Lemma~\ref{Lemma:Clincher}.
But we also have
$$
|\lambda|_\infty\cdot
|T|_\infty^{N-\frac{q}{q-1}}\leq \sup_{\begin{subarray}{c}
x\in \CC_\infty\\
|x|_\infty\leq |T|_\infty
\end{subarray}}|E(x)|_\infty
\leq \sup_{\begin{subarray}{c}
x\in \CC_\infty\\
|x|_\infty\leq |T|_\infty
\end{subarray}}|\psi(x)|_\infty\cdot
\size{h}\cdot|T|_\infty^{N\left(1-\frac{1}{2\ell}\right)},
$$
for all $N$, the inequality on the left by the
Schwarz-Jensen formula, and hence
$$|\lambda|_\infty = O\left(|T|_\infty^{-\frac{N}{2\ell}}\right)\;\;\;
\text{as $N \to \infty$}.$$
These bounds for $|\lambda|_\infty$ are
contradictory for
$N\gg 0$. Our claim is proved.

\subsubsection{The case $E=0$}
Now fix a value of $N$ such that the auxiliary function $E$ vanishes
identically. Since the entries of the vector
$h$ are polynomials in $t$ of degree $<N$, not all vanishing
identically, there exists some $0\leq \nu<N$ such that
$$h^{(N+\nu)}(T)=h\left(T^{q^{-(N+\nu)}}\right)^{q^{N+\nu}}\neq
0.$$ Put
$$P=P(t):= h^{(N+\nu)}\Theta^{(N+\nu)}\cdots\Theta^{(1)}
\in
\Mat_{1\times \ell}(\OO[t]).$$ 
Since
$$\left.\det 
\left(\Theta^{(N+\nu)}\cdots \Theta^{(1)}\right)\right|_{t=T}\neq 0,$$
we have 
$$P(T)\neq 0.$$
We also have
$$P(T)\vartheta=\left(\left.h\Theta^{(-0)}\cdots\Theta^{(-(N+\nu-1))}
\vartheta^{(-(N+\nu))}
\right|_{t=T^{q^{-(N+\nu)}}}\right)^{q^{N+\nu}}=0,$$
and hence
$$P(T)\in (\mbox{$K$-span of $\rho$})\subset\Mat_{1\times \ell}(K).$$
Finally, we have 
$$
 P \psi=
c^{q+\dots+q^{N+\nu}}(t-T^q)^s\cdots\left(t-T^{q^{N+\nu}}\right)^s
E^{(N+\nu)}=0
$$ 
by the key identity.  Therefore (up to a nonzero correction factor in
$K$) the vector $P$ is the vector we want, and the proof of
Theorem~\ref{Theorem:TranscendenceCriterion} is complete.

\section{Tools from (non)commutative algebra}
\label{section:NoncommutativeAlgebra}

\subsection{The ring $\kbar[\sigma]$}

\subsubsection{Definition} Let 
$\kbar[\sigma]$ be the ring obtained by adjoining a noncommutative
variable $\sigma$ to $\kbar$ subject to the commutation relations
$$\sigma x=x^{q^{-1}}\sigma\;\;\;(x\in \kbar).$$
Every element of
$\kbar[\sigma]$ has a unique presentation of the form
$$\sum_{i=0}^\infty a_i\sigma^i\;\;\;(a_i\in
\kbar,\;\;a_i=0\;\mbox{for}\;i\gg 0),$$
and in terms of such presentations the multiplication
law in
$\kbar[\sigma]$ takes the form
$$\textstyle\left(\sum_ia_i\sigma^i\right)
\left(\sum_jb_j\sigma^j\right)=
\sum_i\sum_j a_ib_j^{q^{-i}}\sigma^{i+j}.$$
Given
$$\phi=\sum_{i=0}^\infty a_i\sigma^i\in \kbar[\sigma]\;\;\;(a_i\in
\kbar,\;\;
\mbox{$a_i=0$ for $i\gg 0$}),$$
we define
$$\deg_\sigma \phi:=\max(\{-\infty\}\cup\{i\vert a_i\neq 0\}).$$
Clearly we have
$$\deg_\sigma \phi\psi=\deg_\sigma \phi+ \deg_\sigma
\psi\;\;\;(\phi,\psi\in \kbar[\sigma]).$$
The ring $\kbar[\sigma]$ admits interpretation as the opposite of the
ring of $\FF_q$-linear endomorphisms of the additive group  over
$\kbar$. This interpretation is not actually needed in the sequel but
might serve as a guide to the intuition of the reader.
\subsubsection{Division algorithms and their uses}
\label{subsubsection:DivisionAlgorithms}
 The ring
$\kbar[\sigma]$ has a left (resp., right) division algorithm:\\
\begin{itemize}
\item For all
$\psi,\phi\in \kbar[\sigma]$ such that $\phi\neq 0$
there exist unique
$\theta,\rho\in \kbar[\sigma]$  such that
$\psi=\phi\theta+\rho$ (resp.,
$\psi=
\theta
\phi+\rho)$ and $\deg_\sigma\rho<\deg_\sigma\phi$.\\
\end{itemize} 
Some especially useful properties 
of $\kbar[\sigma]$ and of left modules over it readily deducible from the
existence of left and right division algorithms 
are as follows:\\
\begin{itemize}
\item Every left ideal of
$\kbar[\sigma]$ is principal. 
\item Every finitely generated left $\kbar[\sigma]$-module is
noetherian.
\item $\dim_\kbar \kbar[\sigma]/\kbar[\sigma]\phi=\deg_\sigma
\phi<\infty$ for all $0\neq \phi\in \kbar[\sigma]$.
\item For every
matrix
$\phi\in\Mat_{r\times s}(\kbar[\sigma])$ there exist matrices
$\alpha\in\GL_r(\kbar[\sigma])$ and $\beta\in \GL_s(\kbar[\sigma])$ such
that the product
$\alpha \phi \beta$ vanishes off the main diagonal.
\item A finitely generated free left $\kbar[\sigma]$-module has a
well-defined rank, i.~e., all
$\kbar[\sigma]$-bases  have the same cardinality.
\item A $\kbar[\sigma]$-submodule of a free left
$\kbar[\sigma]$-module of rank $s<\infty$ is free
of rank $\leq s$.
\item Every finitely generated left $\kbar[\sigma]$-module
is isomorphic to a finite direct sum of cyclic left
$\kbar[\sigma]$-modules.
\\
\end{itemize}
These facts are quite well known. The proofs run
along lines very similar to the proofs of the analogous statements
for, say,  the commutative ring $\kbar[t]$.

\subsubsection{The functors $\bmod{\,\sigma}$ and $\bmod{\,(\sigma-1)}$}
Given a homomorphism 
$$f:H_0\rightarrow H_1$$
of left $\kbar[\sigma]$-modules, let
$$f\bmod{\sigma}:\frac{H_0}{\sigma H_0}\rightarrow
\frac{H_1}{\sigma H_1},\;\;\;
f\bmod{(\sigma-1)}:\frac{H_0}{(\sigma-1)H_0}
\rightarrow
\frac{H_1}{(\sigma-1)H_1}$$
be the corresponding induced maps.

\begin{Lemma}\label{Lemma:InclusionFiniteness}
Let 
$$f:H_0\rightarrow H_1$$
be an injective homomorphism of free left
$\kbar[\sigma]$-modules of finite rank such that
$$n:=\dim_\kbar \coker(f)<\infty.$$
We have
$$\card\ker(f\bmod{(\sigma-1)})\leq
q^n$$
with equality if and only if $f\bmod{\sigma}$ is
bijective.
\end{Lemma}
\proof We may assume without loss of generality that
$$H_0=\Mat_{1\times r}(\kbar[\sigma]),\;\;\;H_1=\Mat_{1\times
s}(\kbar[\sigma]),\;\;\; f=(x\mapsto x\phi)\;\;\;(\phi\in
\Mat_{r\times s}(\kbar[\sigma])).$$ After replacing $\phi$
by
$\alpha\phi\beta$ for
suitably chosen
$\alpha\in \GL_r(\kbar[\sigma])$ and $\beta\in
\GL_s(\kbar[\sigma])$,
we
may assume without loss of generality that
$\phi$ vanishes off the main diagonal, in which case clearly
$\phi$ vanishes nowhere on the main
diagonal and
$r=s$. We might as well assume now also that
$r=s=1$. Write
$$\phi=\sum_{i=0}^na_i\sigma^i\;\;\;(a_i\in
\kbar,\;\;\;a_n\neq
0,\;\;\;a_0\neq 0\Leftrightarrow f\bmod{\sigma}\mbox{ is
bijective}).$$ 
We have
$$\kbar[\sigma]=\kbar\oplus(\sigma-1)\cdot\kbar[\sigma]$$
and
$$\;\;\;x\phi
\equiv\sum_{i=0}^na_i^{q^i}x^{q^i}
\bmod{(\sigma-1)\cdot\kbar[\sigma]}$$
for all $x\in \kbar$,
whence the result.
\qed

\begin{Lemma}\label{Lemma:YuLemmaAnalogue}
For $i=1,2$ let 
$$f_i:H_0\rightarrow H_i$$
be a homomorphism of free left $\kbar[\sigma]$-modules of
finite rank.
Assume that $H_0$, $H_1$
and $H_2$ are all of the same rank over $\kbar[\sigma]$. 
Assume further that $f_1\bmod{\sigma}$ is
bijective
 and that
$$\ker(f_1\bmod{(\sigma-1)})
\subset\ker(f_2\bmod{(\sigma-1)}).$$
Then $f_2$ factors uniquely through $f_1$, i.~e., there exists a unique
homomorphism 
$$g:H_1\rightarrow H_2$$
of left $\kbar[\sigma]$-modules such that 
$$g\circ f_1=f_2.$$
\end{Lemma}
\proof (Cf.~\cite[Lemma~1.1]{Yu}.) 
We may assume without loss of generality that
$$H_0=H_1=
H_2=
\Mat_{1\times s}(\kbar[\sigma]),$$
$$f_1=(x\mapsto x\phi),\;\;f_2=(x\mapsto x\psi),\;\;\;(\phi,\psi\in
\Mat_{s\times
s}(\kbar[\sigma])).$$ After making replacements
$$\phi\leftarrow\alpha\phi\beta,\;\;\;\psi\leftarrow\alpha\psi$$
for suitably chosen $\alpha,\beta\in \GL_s(\kbar[\sigma])$,
we may
assume without loss of generality that
$\phi$ vanishes off the main diagonal. Since $f_1\bmod{\sigma}$ is
bijective, no diagonal entry of $\phi$ vanishes. We might as well 
assume now also that $s=1$. 
Use the left division algorithm to find $\theta,\rho\in
\kbar[\sigma]$ such that
$$\psi=\phi\theta+\rho,\;\;\;
\deg_\sigma \rho<\deg_\sigma\phi.$$
Put
$$g:=(x\mapsto x\theta),\;\;\;h:=(x\mapsto x\rho).$$
Then
$$f_2=g\circ f_1+h.$$
If $h=0$  we are done. Suppose instead that $h\neq
0$. We then have
$$\ker
(f_1\bmod{(\sigma-1)})\subset\ker(h\bmod{(\sigma-1)}),\;\;\;
\dim_\kbar\coker(f_1)>\dim_\kbar\coker(h).$$
But the latter relations are
contradictory in view of
Lemma~\ref{Lemma:InclusionFiniteness} and 
 our hypothesis
that $f_1\bmod{\sigma}$ is bijective.
\qed

\subsection{The  ring $\kbar[\![\sigma]\!]$}

\subsubsection{Definition}
We define $\kbar[\![\sigma]\!]$ to be the completion
of
$\kbar[\sigma]$ with respect to the system of two-sided ideals
$$\left\{\sigma^n\kbar[\sigma]\right\}_{n=0}^\infty.$$ 
Every element of
$\kbar[\![\sigma]\!]$ has a unique presentation of the form
$$\sum_{i=0}^\infty a_i\sigma^i\;\;\;(a_i\in\kbar).$$
In terms of such presentations
the multiplication law in
$\kbar[\![\sigma]\!]$ takes the form
$$\textstyle\left(\sum_ia_i\sigma^i\right)
\left(\sum_jb_j\sigma^j\right)=
\sum_i\sum_j a_ib_j^{q^{-i}}\sigma^{i+j}.$$
The ring $\kbar[\![\sigma]\!]$ contains $\kbar[\sigma]$ as a subring.
The ring $\kbar[\![\sigma]\!]$ is a domain.

\subsubsection{The operation $\partial$}
Given
$$\phi=\sum_{i=0}^\infty a_{(i)}\sigma^i
\in \Mat_{r\times s}(\kbar[\![\sigma]\!])\;\;\;(a_{(i)}\in \Mat_{r\times
s}(\kbar)),$$
put
$$\partial\phi:=a_{(0)}.$$
The operation $\partial$ thus defined is $\kbar$-linear and satisfies
$$\partial(\phi\psi)=(\partial\phi)(\partial\psi)$$
for all matrices $\phi$ and $\psi$ with entries in $\kbar[\![\sigma]\!]$
such that the product $\phi\psi$ is defined. 

\begin{Lemma}\label{Lemma:DoubleBracketInversion}
(i) For all $\phi\in\Mat_{s\times
s}(\kbar[\![\sigma]\!])$, if $\partial\phi\in \GL_s(\kbar)$,
then $\phi\in \GL_s(\kbar[\![\sigma]\!])$. 
(ii)
 Every nonzero left ideal of
$\kbar[\![\sigma]\!]$ is generated by a power of $\sigma$.
\end{Lemma}
\proof (i) After replacing
$\phi$ by
$\alpha\phi$ for suitably chosen $\alpha\in \GL_s(\kbar)$ we may
assume 
$\partial\phi=\one_s$. Now write
$\phi=\one_s-X$.
The series
$\one_s+\sum_{n=1}^\infty X^n$
converges to a two-sided inverse to $\phi$. 
(ii) Let $I\subset\kbar[\![\sigma]\!]$ be a nonzero left ideal.
Let $\phi=\alpha\sigma^n$ be a nonzero element of $I$
where $\partial \alpha\neq 0$ and $n$ is a nonnegative integer taken
as
small as possible. Then we have $\alpha\in \kbar[\![\sigma]\!]^\times$ by
(i), hence
$\sigma^n\in I$, and hence $\sigma^n$ generates $I$.
\qed

\begin{Lemma}\label{Lemma:ExpExistence}
Let
$$\theta\in \Mat_{r\times
r}(\kbar[\![\sigma]\!]),\;\;\;a\in \Mat_{r\times
r}(\kbar),\;\;e\in\Mat_{r\times s}(\kbar),\;\; b\in \Mat_{s\times
s}(\kbar)$$ be given such that
$$
\partial\theta=a,\;\;\;(a-T\cdot\one_r)^r=0,\;\;ae=eb,\;\;(b-T\cdot\one_s)^s=0.$$
Then there exists unique
$$E\in\Mat_{r\times
s}(\kbar[\![\sigma]\!])$$
such that
$$\theta E=Eb,\;\;\;\partial E=e.$$
\end{Lemma}
\proof (Cf.~\cite[Prop.~2.1.4]{AndersonMotive}.)
Write
$$\theta=\sum_{i=0}^\infty a_{(i)}\sigma^i\in \Mat_{r\times
r}(\kbar[\![\sigma]\!])\;\;(a_{(i)}\in \Mat_{r\times r}(\kbar),\;
a_{(0)}=a)$$
and
$$E=\sum_{i=0}^\infty e_{(i)}\sigma^i
\in \Mat_{r\times
s}(\kbar[\![\sigma]\!])\;\;(e_{(i)}\in \Mat_{r\times s}(\kbar),\;e_{(0)}=e).
$$
We have
$$\theta E=\sum_{i=0}^\infty\sum_{j=0}^\infty
a_{(i)}e_{(j)}^{(-i)}\sigma^{i+j}$$
and hence $\theta E=Eb$ if  and only if 
the system of coefficients $\left\{e_{(n)}\right\}$ satisfies the
recursion
$$
e_{(n)}b^{(-n)}-ae_{(n)}=\sum_{0<i\leq
n}a_{(i)}e_{(n-i)}^{(-i)}\;\;\;\;\;\;(n=1,2,\dots).$$
Now for each
$n>0$ the $\kbar$-linear map
$$\left(z\mapsto
zb^{(-n)}-az\right):
\Mat_{r\times s}(\kbar)\rightarrow \Mat_{r\times s}(\kbar)$$
is invertible because all of its eigenvalues equal $T^{q^{-n}}-T$;
indeed, for $q^i \geq \max\{r,s\}$
the $q^i$-th 
iteration  sends
$z$  to $(T^{q^{-n}}-T)^{q^i}z$.
It follows that the recursion
satisfied by the coefficients $e_{(n)}$ has a unique solution
with $e_{(0)}=e$.  \qed

\begin{Lemma}\label{Lemma:Stiffness}
Let matrices
$$\theta,\phi,\rho\in \Mat_{s\times s}(\kbar[\sigma])$$
be given such that
$$
(\partial\theta-T\cdot \one_s)^s=0,\;\;\;\theta
\phi=\phi\rho,\;\;\; (\partial\rho-T\cdot\one_s)^s=0,$$
Assume further that the map
$$(x\mapsto x\phi):\Mat_{1\times s}(\kbar[\sigma])
\rightarrow\Mat_{1\times s}(\kbar[\sigma])$$
is injective. Then 
$$\det\partial\phi\neq 0.$$
\end{Lemma}
\proof (Cf.~\cite[Lemma 1.3]{YuEarlier}.) After making
replacements
$$\theta\leftarrow\alpha\theta\alpha^{-1},\;\;\;
\phi\leftarrow\alpha\phi\beta,\;\;\;\rho\leftarrow\beta^{-1}\rho\beta,$$
for suitably chosen $\alpha,\beta\in \GL_s(\kbar[\sigma])$ we
may assume without loss of generality that
$\phi$ vanishes off the main diagonal, in which case it is clear that no
entry of $\phi$ on the main diagonal vanishes.
By Lemma~\ref{Lemma:ExpExistence} there exist unique
matrices
$$
D,E,F\in \Mat_{s\times s}(\kbar[\![\sigma]\!])
$$ 
such that
$$\theta D=D\cdot\partial \theta,\;\;\;\theta E=E\cdot\partial
\rho,\;\;\;
\rho F=F\cdot\partial \rho,$$
and
$$\partial D=\one_r,\;\;\;\partial E=\partial \phi,\;\;\;
\partial F=\one_s.$$
By the uniqueness asserted in 
Lemma~\ref{Lemma:ExpExistence} 
we have $$\phi F=E=D\cdot\partial\phi.$$
Further, we have 
$$D,F\in \GL_s(\kbar[\![\sigma]\!])$$
by Lemma~\ref{Lemma:DoubleBracketInversion} (i).
Now consider the quotient $M$ (resp., $N$) of the free left
$\kbar[\![\sigma]\!]$-module
$\Mat_{s\times 1}(\kbar[\![\sigma]\!])$
by the $\kbar[\![\sigma]\!]$-submodule generated by the rows of $\phi$
(resp., $\partial
\phi$). Since $\phi$ is diagonal with no vanishing diagonal
entries, we have
$\dim_\kbar M<\infty$ by Lemma~\ref{Lemma:DoubleBracketInversion} (ii).
Since
$D^{-1}\phi F=\partial
\phi$, the left
$\kbar[\![\sigma]\!]$-modules
$M$ and
$N$ are isomorphic and hence we have $\dim_\kbar N<\infty$. 
Under the latter condition it is impossible for any diagonal
entry of
$\partial
\phi$ to vanish.
\qed

\subsection{The  ring $\kbar[t,\sigma]$}

\subsubsection{Definition}
Let $\kbar[t,\sigma]$ be the ring obtained by adjoining the commutative
variable
$t$ to $\kbar[\sigma]$. Every element of $\kbar[t,\sigma]$ has a
unique presentation of the form
$$\sum_{i=0}^\infty \alpha_it^i\;\;\;(\alpha_i\in
\kbar[\sigma],\;\alpha_i=0\mbox{ for
$i\gg 0$}).
$$
In terms of such presentations the multiplication law in
$\kbar[t,\sigma]$ takes the form
$$\textstyle\left(\sum_i \alpha_{i}t^i\right)
\left(\sum_j \beta_{j}t^j\right)=
\sum_i\sum_j \alpha_{i}\beta_{j}t^{i+j}.$$
Every element of $\kbar[t,\sigma]$ also has a unique presentation of
the form
$$\sum_{i=0}^\infty a_i\sigma^i\;\;\;(a_i\in
\kbar[t],\;\;a_i=0\;\mbox{for}\;i\gg 0).$$
In terms of such presentations the multiplication law in
$\kbar[t,\sigma]$ takes the form
$$\textstyle \left(\sum_i a_i\sigma^i\right)
\left(\sum_j  b_j\sigma^j\right)=
\sum_i\sum_j a_ib_j^{(-i)}\sigma^{i+j}.$$
The ring
$\kbar[t,\sigma]$ contains both the noncommutative ring
$\kbar[\sigma]$ and the commutative ring
$\kbar[t]$ as subrings.  
The
ring
$\FF_q[t]$ is contained in the center of the ring $\kbar[t,\sigma]$. 
The ring $\kbar[t,\sigma]$ is a domain.

\begin{Proposition}\label{Proposition:Seesaw}
Let $M$ be a left $\kbar[t,\sigma]$-module finitely generated over both
$\kbar[\sigma]$ and $\kbar[t]$. 
Let $M_\sigma$ (resp., $M_t$) be the sum of all
$\kbar[\sigma]$- (resp., $\kbar[t]$-) submodules
$N\subset M$ such
that $\dim_\kbar
N<\infty$. Then 
$M_\sigma=M_t$,  $\dim_\kbar M_\sigma=\dim_\kbar M_t<\infty$
and the quotient
$M/M_\sigma=M/M_t$ is free of finite rank over both $\kbar[\sigma]$ and
$\kbar[t]$.
\end{Proposition}
In particular it follows that a left
$\kbar[t,\sigma]$-module finitely generated over both $\kbar[\sigma]$
and $\kbar[t]$ is free over
$\kbar[\sigma]$ if and only if free over
$\kbar[t]$.
\proof (Cf.~\cite[Lemma 1.4.5, p.~463]{AndersonMotive}.) 
As a $\kbar[\sigma]$-module $M$ decomposes as a
finite direct sum of cyclic left
$\kbar[\sigma]$-modules. Therefore 
$\dim_\kbar M_\sigma<\infty$ and
$M/M_\sigma$ is free of finite rank over $\kbar[\sigma]$. Similarly 
$\dim_\kbar M_t<\infty$ and $M/M_t$ is free of finite rank over
$\kbar[t]$. Now for any
$\kbar[t]$-submodule $N\subset M$ of finite dimension over $\kbar$,
again $\sigma N$ is a $\kbar[t]$-submodule of $M$ of finite dimension
over $\kbar$ and hence $\sigma M_t\subset M_t$. Similarly we have
$tM_\sigma\subset M_\sigma$. Therefore each
of the modules $M_\sigma$ and
$M_t$ contains the other.
\qed

\subsubsection{Saturation}
Let 
$N\subset M$ be left $\kbar[t,\sigma]$-modules. Assume that $M$ is
finitely generated over both $\kbar[\sigma]$ and $\kbar[t]$. By
Proposition~\ref{Proposition:Seesaw} there exists a unique
$\kbar[t,\sigma]$-submodule $\tilde{N}\subset M$ 
such that
$\tilde{N}=\sum N'=\sum N''$
where $N'$ ranges over $\kbar[\sigma]$-submodules
such that $\dim_\kbar (N'+N)/N<\infty$ and
$N''$ ranges over $\kbar[t]$-submodules such that
$\dim_{\kbar}(N''+N)/N<\infty$. Clearly we have $\tilde{N}\supset N$.
 We call $\tilde{N}$
the {\em saturation} of $N$ in $M$. 
By Proposition~\ref{Proposition:Seesaw} the quotient
$M/\tilde{N}$ is free of finite rank over both $\kbar[\sigma]$ and 
$\kbar[t]$ and moreover $\dim_\kbar\tilde{N}/N<\infty$. If
$N=\tilde{N}$ we say that $N$ is {\em saturated} in $M$. 
A
necessary and sufficient condition for $N$ to be saturated in $M$ is that
$M/N$ be torsion-free over $\kbar[t]$ or torsion-free
over $\kbar[\sigma]$.

\subsection{Dual $t$-motives}
\label{subsection:DualtMotives}
\nopagebreak

\subsubsection{Definition}
A {\em dual $t$-motive} $H$ is a left $\kbar[t,\sigma]$-module with
the following three properties:
\begin{itemize}
\item $H$ is free of finite rank over $\kbar[t]$.
\item $H$ is free of finite rank over $\kbar[\sigma]$.
\item $(t-T)^nH\subset\sigma H$ for $n\gg 0$.
\end{itemize}
 A {\em morphism} of dual
$t$-motives is by definition a homomorphism of left
$\kbar[t,\sigma]$-modules. Thus dual $t$-motives form a category.
For any dual
$t$-motive
$H$ there exist
$$\ggg\in \Mat_{r\times 1}(H),\;\;\;\hh\in\Mat_{s\times
1}(H),\;\;\;\;\Phi\in
\Mat_{r\times r}(\kbar[t]),\;\;\;\;\theta\in\Mat_{s\times
s}(\kbar[\sigma])$$ such that
\begin{itemize}
\item the entries of $\ggg$ form a $\kbar[t]$-basis for $H$ and $\sigma
\ggg=\Phi\ggg$, 
\item the entries of $\hh$ form a $\kbar[\sigma]$-basis for $H$
 and
$t\hh=\theta\hh$,
\end{itemize}
in which case 
\begin{itemize}
\item $\det \Phi=c(t-T)^s$ for some nonzero $c\in \kbar$ and
\item $(\partial \theta-T\cdot \one_s)^s=0$.
\end{itemize}
The latter two assertions both follow from the assumption that the
quotient
$H/\sigma H$ is killed by a power of $t-T$. 

\subsubsection{Basic stability properties of the class of dual
$t$-motives} Let 
$$H_0\subset H_1$$ be left $\kbar[t,\sigma]$-modules.
The following statements hold:\\
\begin{itemize}
\item If
$H_1$ is a dual
$t$-motive and
$(t-T)^n H_0\subset \sigma H_0$ for
$n\gg 0$, then $H_0$ is again a dual $t$-motive.
\item If $H_1$ is a dual $t$-motive and $H_0$ is saturated in $H_1$,
then both $H_0$ and $H_1/H_0$ are again dual $t$-motives.
\item If $H_0$ and $H_1/H_0$ are both dual $t$-motives, then $H_1$ is
again a dual
$t$-motive. \\
\end{itemize}
Using the background on module theory over $\kbar[\sigma]$
and $\kbar[t,\sigma]$ provided
above the reader should have no difficulty verifying these
statements.

\begin{Proposition}
\label{Proposition:Rationale}
Let 
$$\Phi\in
\Mat_{\ell\times \ell}(\kbar[t]),\;\;\;\psi\in \Mat_{\ell\times 1}(\EE)$$
be as in Theorem~{\ref{Theorem:TranscendenceCriterion}}.
 Suppose further that there exist a dual
$t$-motive
$H$ and a vector 
$$\ggg\in \Mat_{\ell\times 1}(H)$$
with entries forming a $\kbar[t]$-basis of $H$ such that
$$\sigma \ggg=\Phi\ggg.$$
Equip $\EE$ with left $\kbar[t,\sigma]$-module
structure by the rule
$$\sigma e:=e^{(-1)}\;\;\;(e\in \EE).$$
Put
$$\begin{array}{rcl}
H_0&:=&\textup{\mbox{$\kbar[t]$-span in $\EE$ of the entries
of  the vector $\psi$}},\\
V&:=&\textup{\mbox{$\kbar$-span in
$\overline{k_\infty}$ of the entries of the vector $\psi(T)$}}.
\end{array}
$$ Then the following statements hold:
\begin{itemize}
\item
$H_0$ is a
$\kbar[t,\sigma]$-submodule of
$\EE$. 
\item $H_0$ is a dual $t$-motive admitting presentation as a quotient of 
$H$. 
\item $\rank_{\kbar[t]}H_0=\dim_\kbar V$.
\end{itemize}
\end{Proposition}
The proposition positions
Theorem~\ref{Theorem:TranscendenceCriterion} in the setting
of dual $t$-motives.  
\proof Consider the exact sequence
$$0\rightarrow H_1\subset H\rightarrow H_0\rightarrow 0$$
of left $\kbar[t]$-modules, where
$$H_1:=\{P\ggg\in H\mid P\in \Mat_{1\times \ell}(\kbar[t]),\;
P\psi=0\}$$
and the projection $H\rightarrow H_0$ is given by the rule
$$P\ggg\mapsto P\psi\;\;\;(P\in \Mat_{1\times \ell}(\kbar[t])).$$
A straightforward calculation verifies that the exact sequence
in question is in fact an exact sequence of $\kbar[t,\sigma]$-modules.
Therefore $H_0$ is a dual $t$-motive admitting presentation as a
quotient of $H$. Since every $\kbar[t]$-basis for $H_1$ can be completed to
a $\kbar[t]$-basis of $H$,
the number of
$\kbar$-linearly independent relations of
$\kbar$-linear dependence among the entries of $\psi(T)$ 
is at least as great as $\rank_{\kbar[t]}H_1$
and hence
we have
$$\rank_{\kbar[t]} H_0\geq \dim_\kbar V.$$
But we also have
$$\rank_{\kbar[t]}H_0\leq \dim_\kbar V$$
because by Theorem~\ref{Theorem:TranscendenceCriterion}
every relation of $\kbar$-linear dependence among
the entries of $\psi(T)$ lifts to a $\kbar[t]$-linear relation
among the entries of $\psi$.
\qed

\begin{Theorem}\label{Theorem:2}
For
all
dual $t$-motives $H_0$ and $H_1$ the natural map
$$\kbar\otimes_{\FF_q}\Hom_{\kbar[t,\sigma]}(H_0,H_1)\rightarrow
\Hom_{\kbar[t]}(H_0,H_1)$$
is injective. 
\end{Theorem}
\proof The proof of
\cite[Thm.~2, p.~464]{AndersonMotive} can easily be modified to prove this
result.
\qed

\subsubsection{Isogenies}
An injective morphism
$f:H_0\rightarrow H_1$ of dual $t$-motives with cokernel
finite-dimensional over $\kbar$ is called an {\em isogeny}.
We say that dual $t$-motives $H_0$ and $H_1$ are {\em isogenous}
if there exists an isogeny $f:H_0\rightarrow H_1$.

\begin{Lemma}
\label{Lemma:IsogenyTangentBehavior}
Let $f:H_0\rightarrow H_1$ be an isogeny of dual $t$-motives.
Then the induced map $f\bmod{\sigma}$ is bijective.
\end{Lemma}
\proof
Without loss of generality, we think of $f$ as an inclusion.
Let $s$ be the common rank of 
$H_0$ and $H_1$ over
$\kbar[\sigma]$. For
$i=0,1$ select 
$$\hh_{(i)}\in \Mat_{s\times 1}(H_i),\;\;\;\theta_{(i)}
\in \Mat_{s\times s}(\kbar[\sigma])$$ 
such that the entries of $\hh_{(i)}$ form a $\kbar[\sigma]$-basis of
$H_i$ and 
$$t\hh_{(i)}=\theta_{(i)}\hh_{(i)}.$$
 Let 
$$\phi\in \Mat_{s\times
s}(\kbar[\sigma])$$  be the unique solution of
the equation 
$$\hh_{(0)}=\phi\hh_{(1)},$$
noting that
$$\theta_{(0)}\phi\hh_{(1)}=\theta_{(0)}\hh_{(0)}=
t\hh_{(0)}=t\phi\hh_{(1)}=\phi t\hh_{(1)}=\phi\theta_{(1)}\hh_{(1)}.$$
By Lemma~\ref{Lemma:Stiffness}
we have $\det \partial \phi\neq 0$
and hence $f\bmod{\sigma}$ is bijective.
\qed

\begin{Theorem}\label{Theorem:IsogenySymmetry} Let
$H_0\subset H_1$ be dual $t$-motives such that
$\dim_\kbar H_1/H_0<\infty$. Then there exists $0\neq
a\in
\FF_q[t]$ such that $aH_1\subset H_0$. 
\end{Theorem}
It follows that the isogeny relation is not only reflexive and
transitive, but symmetric as well and hence an equivalence relation.
\proof  Let 
$f_1:H_0\rightarrow H_1$
be the inclusion. By Lemma~\ref{Lemma:IsogenyTangentBehavior} 
the induced map
$f_1\bmod{\sigma}$ is bijective. By
Lemma~\ref{Lemma:InclusionFiniteness} the kernel of the induced map
$f_1\bmod{(\sigma-1)}$ is finite, and, since the latter group naturally
has the structure of finite
$\FF_q[t]$-module, there exists $0\neq a\in \FF_q[t]$
killing it. 
Let $f_2:H_0\rightarrow H_0$ be the morphism of dual
$t$-motives induced by multiplication by
$a$.  Then
$$\ker(f_1\bmod{(\sigma-1)})\subset\ker(f_2\bmod{(\sigma-1)}).$$
By Lemma~\ref{Lemma:YuLemmaAnalogue} there exists a
unique
$\kbar[\sigma]$-module homomorphism 
$g:H_1\rightarrow H_0$
such that  
$f_2=g\circ f_1$. We  have
$$\ker(f_1\bmod{(\sigma-1)})\subset\ker(tf_2\bmod{(\sigma-1)})$$
and
$$tg(f_1(h))=tf_2(h)=f_2(th)=g(f_1(th))=g(tf_1(h))\;\;\;(h\in H_0).$$
By the uniqueness asserted in
Lemma~\ref{Lemma:YuLemmaAnalogue}, it follows that $g$
commutes with $t$ and hence is a morphism of dual $t$-motives. We
have
$$\ker(f_1\bmod{(\sigma-1)})\subset\ker(f_1\circ
g\circ f_1\bmod{(\sigma-1)})$$ and
$$f_1(g(f_1(h)))=f_1(f_2(h))=f_1(ah)=af_1(h)\;\;\;(h\in H_0).$$
By the uniqueness asserted 
in Lemma~\ref{Lemma:YuLemmaAnalogue}, it follows that
$f_1\circ g$ coincides with 
multiplication by $a$. We therefore have $aH_1\subset H_0$.
\qed

\begin{Corollary}\label{Corollary:RankIsogenyInvariance}
For all dual $t$-motives $H_0$ and $H_1$ the module
$$\Hom_{\kbar[t,\sigma]}(H_0,H_1)$$ is free over $\FF_q[t]$ of finite
rank and moreover its rank over $\FF_q[t]$
depends only on the isogeny classes of $H_0$ and $H_1$.
\end{Corollary}
\proof Theorem~\ref{Theorem:2} already proves that
the module in question is free of finite rank over $\FF_q[t]$.
Now let
$r(H_0,H_1)$ denote the rank over
$\FF_q[t]$ of the module in question.
For $i=0,1$ let $H_i'$
be a dual $t$-motive isogenous to $H_i$ and without loss of generality
assume that $H'_i\subset H_i$ and $\dim_\kbar H_i/H'_i<\infty$.
Choose $0\neq a\in \FF_q[t]$ such that $aH_i\subset H'_i$
for $i=0,1$. We have
$$r(H_0,H_1)\leq r(H'_0,H_1)\leq r(aH_0,H_1)=r(H_0,H_1)$$
and
$$r(H_0,H_1)=r(H_0,aH_1)\leq
r(H_0,H'_1)\leq r(H_0,H_1),$$
where each inequality is justified by the existence of a
suitably constructed 
injective $\FF_q[t]$-linear map.
\qed

\subsubsection{Simplicity}
We say that a dual $t$-motive $H$ is {\em simple}
if $H\neq \{0\}$ and there exist no saturated
$\kbar[t,\sigma]$-submodules of $H$ other than $\{0\}$ and $H$. 
\begin{Proposition}\label{Proposition:SimpleDecomposition}
 (i) A dual $t$-motive isogenous to a simple
dual
$t$-motive is again simple. 
(ii) A nonzero morphism of dual $t$-motives with
simple source and target is automatically an isogeny. 
(iii) Let $\{H_i\}$ be a family of simple dual
$t$-motives
each embedded as a $\kbar[t,\sigma]$-submodule of a dual
$t$-motive $H$. If 
$$\dim_\kbar H/\left(\sum_i H_i\right)<\infty,$$ then $H$
is isogenous to a finite direct  sum of dual
$t$-motives of the family $\{H_i\}$.
\end{Proposition}
\proof (i) Let $H_0\subset H_1$ be dual $t$-motives with
$\dim_\kbar H_1/H_0<\infty$ and $H_0$ simple. It suffices to show that
$H_1$ is simple. Let $M\subset H_1$ be a $\kbar[t,\sigma]$-submodule
saturated in $H_1$ and hence free over $\kbar[\sigma]$
and $\kbar[t]$.  Then $M\cap H_0$
is saturated in
$H_0$ and hence
$M\cap H_0=\{0\}$ or $M\cap H_0=H_0$ by the simplicity of $H_0$. In the
former case,
$M$ injects into $H_1/H_0$, so $\dim_\kbar M<\infty$ and hence $M=0$ since
$M$ is a free $\kbar[\sigma]$-module. In the latter case $\dim_\kbar
H_1/M<\infty$ and hence
$M=H_1$ since 
$M$ is saturated in $H_1$. Therefore $H_1$ is indeed simple.

(ii) Let $f:H_0\rightarrow H_1$ be a nonzero morphism of dual
$t$-motives with simple source and target. The kernel of $f$ is
a saturated $\kbar[t,\sigma]$-submodule of $H_0$ distinct from $H_0$
and hence equal to $\{0\}$. The saturation of the image of $f$
is a saturated
 $\kbar[t,\sigma]$-submodule of $H_1$
distinct from
$\{0\}$, hence equal to $H_1$, and hence the cokernel of $f$ is
of finite dimension over $\kbar$. Therefore $f$ is indeed an
isogeny.

(iii) Since $H$ is noetherian over $\kbar[t,\sigma]$, there exists a
finite set $I$ of indices such that
$$\dim_{\kbar}H/\left(\sum_{i\in I} H_i\right)<\infty.$$ Fix such a set
$I$ now with $\card I$ minimal. Consider the exact
sequence
$$0\rightarrow K\rightarrow \bigoplus_{i\in I}H_i\rightarrow
\sum_{i\in I} H_i\rightarrow 0$$ 
of left $\kbar[t,\sigma]$-modules.
It suffices to show that $K=0$; suppose instead that $K\neq 0$.
In any case
$K$ is a saturated $\kbar[t,\sigma]$-submodule of
a dual
$t$-motive, and hence $K$ is a dual $t$-motive. Let $M$ be a nonzero
saturated
$\kbar[t,\sigma]$-submodule of $K$ of minimal rank over $\kbar[t]$.
Then $M$ is a simple dual $t$-motive. For some index $i_0\in I$
the evident map
$M\rightarrow H_{i_0}$ is nonzero and hence an isogeny, in which case
$\sum_{i\in
I\setminus\{i_0\}}H_{i}$ is of finite $\kbar$-codimension in
$H$ in contradiction to the
 minimality of $\card I$. This contradiction proves that $K=0$.
\qed

\subsubsection{Rigid-analytic triviality}
\label{subsubsection:RigidAnalyticTriviality}
Given a dual $t$-motive $H$, put
$$\tilde{H}:=\CC_\infty\{t\}\otimes_{\kbar[t]}H,$$
 equip
$\tilde{H}$ with an action of $\sigma$ by the rule
$$\sigma(f\otimes h):=f^{(-1)}\otimes h,$$
and put
$$
H^{\Betti}:=\{\mbox{$\sigma$-invariant elements of $\tilde{H}$}\}.$$
We say that
$H$ is {\em rigid analytically trivial} if the natural map
$$\CC_\infty\{t\}\otimes_{\FF_q[t]}H^{\Betti}
\rightarrow\tilde{H}$$
is bijective, cf.~\cite[p.~474]{AndersonMotive}.
\begin{Lemma}\label{Lemma:RATEquation}
Let $H$ be a dual $t$-motive. Select
$$\ggg\in \Mat_{r\times 1}(H),\;\;\;\Phi\in \Mat_{r\times
r}(\kbar[t])$$ such that the entries of $\ggg$ form a
$\kbar[t]$-basis of $H$ and $\sigma\ggg=\Phi\ggg$. 
(i) A
necessary and sufficient condition for
$H$ to be rigid analytically trivial is that there exist a solution
$$\Psi\in \GL_r(\CC_\infty\{t\})$$
of the equation
$$\Psi^{(-1)}=\Phi\Psi.$$
(ii) For any such solution $\Psi$ the entries of the column vector
$\Psi^{-1}\ggg$ form an
 $\FF_q[t]$-basis for $H^{\Betti}$.
\end{Lemma}
In particular if $H$ is rigid analytically trivial, then $H^{\Betti}$
is free over
$\FF_q[t]$ of rank equal to the rank of $H$ over $\kbar[t]$.
Note that $\Psi\in\Mat_{r\times
r}(\EE)$ by Proposition~\ref{Proposition:HypothesisChecker}.
\proof  Note that by 
hypothesis 
and definition 
the diagrams
$$\begin{array}{rcl}
\Mat_{1\times r}(\kbar[t])&\xrightarrow{\times \ggg}&H\\\\
{\scriptstyle {P\mapsto P^{(-1)}\Phi}}\downarrow&&
\downarrow{\scriptstyle{\sigma \times}}\\\\
\Mat_{1\times r}(\kbar[t])&\xrightarrow{\times \ggg}&H
\end{array}\;\;\;\;\;\begin{array}{rcl}
\Mat_{1\times r}(\CC_\infty\{t\})&\xrightarrow{\times
\ggg}&\tilde{H}\\\\ {\scriptstyle {P\mapsto
P^{(-1)}\Phi}}\downarrow&&
\downarrow{\scriptstyle{\sigma \times}}\\\\
\Mat_{1\times r}(\CC_\infty\{t\})&\xrightarrow{\times \ggg}&\tilde{H}
\end{array}
$$ commute and have bijective horizontal arrows.
(i)($\Rightarrow$)
There exists by hypothesis a matrix
$\Theta\in
\Mat_{r\times r}(\CC_\infty\{t\})$ such that the entries of the
column vector
$\Theta\ggg$ are at once a $\CC_\infty\{t\}$-basis for $\tilde{H}$ and
an $\FF_q[t]$-basis for $H^{\Betti}$. Such a matrix
$\Theta$ necessarily belongs to $\GL_r(\CC_\infty\{t\})$ and satisfies
the functional equation
$\Theta^{(-1)}\Phi=\Theta$. The matrix $\Psi:=\Theta^{-1}$
has then the desired properties.
(i)($\Leftarrow$)\&(ii) The entries of the column vector
$\Psi^{-1}\ggg$ form a
$\CC_\infty\{t\}$-basis for
$\tilde{H}$ and are also an $\FF_q[t]$-linearly independent
collection of elements of
$H^{\Betti}$. Every element of $H^{\Betti}$ is of the form
$P\ggg$ for unique
$P\in
\Mat_{1\times r}(\CC_\infty\{t\})$ such that 
$P=P^{(-1)}\Phi$, and we have
$$(P\Psi)^{(-1)}=P^{(-1)}\Phi\Psi=P\Psi\in \Mat_{1\times
r}(\CC_\infty\{t\}\cap\FF_q[\![t]\!])=\Mat_{1\times r}(\FF_q[t]).$$
Therefore the entries of the column vector $\Psi^{-1}\ggg$ span
$H^{\Betti}$ over
$\FF_q[t]$ and hence form an 
$\FF_q[t]$-basis of $H^{\Betti}$. 
\qed
\begin{Lemma}\label{Lemma:DivisionJustification}
For all $0\neq a\in \FF_q[t]$ we have
$\FF_q[t]\cap a\cdot \CC_\infty\{t\}=a\cdot \FF_q[t]$.
\end{Lemma}
\proof View $\CC_\infty\{t\}$ as a subring of the Laurent series field
$\CC_\infty(\!(t)\!)$. We have
$$a^{-1}\FF_q[t]\cap \CC_\infty\{t\}\subset
\FF_q(\!(t)\!)\cap \CC_\infty\{t\}=\FF_q[t],$$
whence the result.
\qed

\begin{Theorem}\label{Theorem:RATConsequences}
For all rigid analytically trivial dual $t$-motives $H_0$ and $H_1$,
the natural map 
$$\Hom_{\kbar[t,\sigma]}(H_0,H_1)\rightarrow
\Hom_{\FF_q[t]}
\left(H_0^{\Betti},H_1^{\Betti}\right)$$
is injective and its cokernel is 
without $\FF_q[t]$-torsion.
\end{Theorem}
\proof After replacing both $H_0$ and $H_1$ by $H_0\oplus H_1$
we may assume without loss of generality that $H_0=H_1$, in which case we
might as well drop subscripts and simply write
$H=H_0=H_1$.
 Let
$\ggg$,
$\Phi$ and
$\Psi$ be as in Lemma~\ref{Lemma:RATEquation}. Fix 
$$e\in \End_{\kbar[t]}(H)$$ arbitrarily, let
$$\tilde{e}\in \End_{\CC_\infty\{t\}}(\tilde{H})$$
be the unique $\CC_\infty\{t\}$-linear extension of $e$, and
let
$$E\in \Mat_{r\times
r}(\kbar[t])$$ be the representation of $e$ with respect to the
$\kbar[t]$-basis $\ggg$, i.~e., the unique solution of the equation
$$e\ggg=E\ggg.$$
In this last and in
analogous expressions below, in order to avoid having to manipulate
indices and summations,
$e$ (resp., $\tilde{e}$) is understood to be applied entrywise to any
column vector with entries in $H$ (resp., $\tilde{H}$) that it
precedes. We have
$$\sigma e\ggg=\sigma E\ggg=E^{(-1)}\sigma
\ggg=E^{(-1)}\Phi\ggg,\;\;\;e\sigma\ggg=e\Phi\ggg=\Phi e\ggg=\Phi
E\ggg,$$ and hence
$$\begin{array}{rcl}
e\in \End_{\kbar[t,\sigma]}(H)&\Leftrightarrow&
E^{(-1)}\Phi=\Phi E\\\\
&\Leftrightarrow&
(\Psi^{-1} E\Psi)^{(-1)}=\Psi^{-1} E\Psi\\\\
&\Leftrightarrow &\Psi^{-1} E\Psi\in \Mat_{r\times
r}(\FF_q[\![t]\!]\cap\CC_\infty\{t\})=\Mat_{r\times r}(\FF_q[t]).
\end{array}$$
We have
$$
\left.\begin{array}{l}
e\in
\End_{\kbar[t,\sigma]}(H)\\\\
\tilde{e}(H^{\Betti})=0
\end{array}\right\}\Rightarrow
0=\tilde{e}\Psi^{-1} \ggg=\Psi^{-1} e\ggg\Rightarrow e\ggg=0\Rightarrow
e=0.$$
Therefore the map in
question is injective. Now fix
$0\neq a\in
\FF_q[t]$. 
We have
$$\begin{array}{rcl}
ae\in
\End_{\kbar[t,\sigma]}(H)&\Rightarrow& a\Psi^{-1} E\Psi\in
\Mat_{r\times r}(\FF_q[t]\cap a\cdot \CC_\infty\{t\})\\\\
&\Rightarrow&\Psi^{-1} E\Psi\in 
\Mat_{r\times r}(\FF_q[t])\\\\
&\Rightarrow&e\in \End_{\kbar[t,\sigma]}(H)
\end{array}$$
where the second implication is justified by
Lemma~\ref{Lemma:DivisionJustification}. Thus we rule out the possibility
of
$\FF_q[t]$-torsion in the cokernel of the map in question.
\qed

\subsection{The Dedekind-Wedderburn trick}
\label{subsection:DedekindWedderburnTrick}
\subsubsection{Dedekind domains}
Recall that a ring with unit is called a {\em Dedekind domain}
if commutative, noetherian, entire ($1\neq 0$ and no zero-divisors), 
one-dimensional (every nonzero prime is maximal), and  
integrally closed. Now let a Dedekind domain $\Kbold$ be given.
The following hold:
\begin{itemize}
\item Every finitely generated $\Kbold$-module without
$\Kbold$-torsion is projective.
\item Every projective $\Kbold$-module of
finite rank is a direct sum of projective
 $\Kbold$-modules of
rank one.
\end{itemize}
 We take these basic facts for granted.
\begin{Lemma}\label{Lemma:Ends}
For every Dedekind domain $\Kbold$ and projective $\Kbold$-module
$\Mbold$ of rank one the map
$(x\mapsto(m\mapsto xm)):\Kbold\rightarrow\End_{\Kbold}(\Mbold)$
is bijective.
\end{Lemma}
\proof Injectivity is clear. Failure of surjectivity would entail
failure of $\Kbold$ to be integrally closed. \qed
\begin{Proposition}\label{Proposition:SeesawTwo}
For $i=0,1$ let $\Lbold_i$ be a Dedekind domain.
Let $\Mbold$ be an abelian group equipped with right $\Lbold_0$-module
structure and left $\Lbold_1$-module structure in such fashion that
$(a_1m)a_0=a_1(ma_0)$
for all $a_1\in \Lbold_1$, $m\in \Mbold$ and $a_0\in \Lbold_0$.
Assume further that $\Mbold$ is projective of rank one both as an
$\Lbold_0$-module and as an $\Lbold_1$-module.
Then there exists a unique
ring isomorphism 
$\theta:\Lbold_0\rightarrow\Lbold_1$ such that
$\theta(a)m=ma$
for all $a\in \Lbold_0$ and $m\in \Mbold$.
\end{Proposition}
\proof This is a trivial (but quite useful) consequence of
Lemma~\ref{Lemma:Ends}.
\qed
\begin{Lemma}\label{Lemma:Descent}
Let $\Kbold\supset\Fbold$ be an integral extension of entire rings.
For every $0\neq b\in \Kbold$ there exists $0\neq a\in \Fbold$
such that $a/b\in \Kbold$.
\end{Lemma}
It follows that in this setting any $\Kbold$-module is
without
$\Kbold$-torsion if (and of course only if) without $\Fbold$-torsion.
\proof Since $\Kbold$ is integral over $\Fbold$, 
there exists a polynomial
$f(z)\in
\Fbold[z]$ monic in $z$ such that $f(b)=0$. Moreover, since $\Kbold$ is
entire, after factoring out a
power of
$z$, we may assume without loss of generality  that $a:=f(0)\neq 0$.  Now
write
$f(z)=a-zg(z)$ with $g(z)\in \Fbold[z]$.  Then $a/b=g(b)\in \Kbold$. \qed

\begin{Proposition}\label{Proposition:JustChecking}
Let $\Lbold\supset\Kbold\supset\Fbold$
be a tower of rings where
$\Lbold$
and $\Fbold$ are Dedekind domains and $\Lbold$ is a finitely
generated projective $\Fbold$-module. Then $\Kbold$ is a
Dedekind domain if and only if the quotient 
$\Lbold/\Kbold$
is a projective $\Fbold$-module.
\end{Proposition}
\proof 

($\Rightarrow$) For all $x\in \Lbold$ and $0\neq a\in
\Fbold$ such that $ax\in \Kbold$ we must have $x\in \Kbold$
since $\Lbold$ is integral over $\Kbold$
and $\Kbold$ is integrally closed.
Therefore the quotient $\Lbold/\Kbold$ is without
$\Fbold$-torsion and since finitely generated
over $\Fbold$ must be projective over $\Fbold$.

($\Leftarrow$) Clearly $\Kbold$ is commutative,
noetherian and entire. Let
$P\subset \Kbold$ be a nonzero prime of $\Kbold$. By
Lemma~\ref{Lemma:Descent} the prime
$P\cap\Fbold$ is nonzero and hence a maximal ideal.
Since $\Kbold$ is integral over $\Fbold$
the prime $P$ has to be maximal, too. Therefore $\Kbold$ is
one-dimensional. 
Since
$\Lbold$ is integrally closed, every element of the integral
closure of $\Kbold$ belongs to $\Lbold$
and hence gives rise to a $\Kbold$-torsion
element of the quotient
$\Lbold/\Kbold$. But by
Lemma~\ref{Lemma:Descent} and hypothesis, the quotient
$\Lbold/\Kbold$ is without
$\Kbold$-torsion.  Therefore
$\Kbold$ is integrally closed. 
Therefore $\Kbold$ is a Dedekind domain.
\qed

\subsubsection{Wedderburn's theorem}
Let $R$ be a (possibly noncommutative) ring with unit.
Let
$M$ be a simple faithful left $R$-module and put
$K:=\End_R(M)$. View $M$ as a left
$K$-module. Assume that
$M$ is finitely generated over
$K$. Then according to {\em Wedderburn's theorem} (a special case of the
{\em Jacobson density theorem}) every element of $\End_K(M)$ is of the
form
$m\mapsto rm$ for unique $r\in R$.

\begin{Proposition}[``The Dedekind-Wedderburn Trick'']
Let $\Lbold\supset\Fbold$ be an extension of Dedekind domains
such that $\Lbold$ is a projective $\Fbold$-module of finite rank.
Let $\Mbold$ be an $\Lbold$-module finitely generated and projective over
$\Fbold$. Assume that $\Lbold$ and $\Mbold$ are of the same rank over
$\Fbold$. Let
$\Rbold$ be a subring
of $\End_{\Fbold}(\Mbold)$
such that
$$\Rbold\supset\{(m\mapsto xm)\in
\End_{\Fbold}(\Mbold)\mid x\in
\Lbold\}$$ and
the quotient
$\End_{\Fbold}(\Mbold)/\Rbold$ is without
$\Fbold$-torsion.  Put
$$\Kbold:=\{x\in \Lbold\mid
(m\mapsto xm)\in\End_{\Rbold}(\Mbold)\}.$$
Then we have 
$$\Rbold=\End_{\Kbold}(\Mbold).$$
\end{Proposition}
The following trivial consequences of the definition
of $\Kbold$ bear emphasis since they are crucial in
applications:
\begin{itemize}
\item 
$\Kbold$ is a Dedekind domain equipped with
$\Fbold$-algebra structure.
\item $\Kbold$ is a projective $\Fbold$-module of finite rank.
\item $\Mbold$ is a finite direct sum of projective
$\Kbold$-modules of rank one.  
\end{itemize}

\proof By Lemma~\ref{Lemma:Descent} and hypothesis the module
$\Mbold$ is projective over $\Lbold$ of rank one and 
by Lemma~\ref{Lemma:Ends} the
evident map $\Lbold\rightarrow\End_{\Lbold}(\Mbold)$ is bijective.
To simplify notation we now identify $\Lbold$
with the subring $\End_{\Lbold}(\Mbold)$ of $\End_{\Fbold}(\Mbold)$.
Note that
$$\Lbold=\End_{\Lbold}(\Mbold)\supset
\End_{\Rbold}(\Mbold)=\Kbold.$$ 
Now put
$$F:=\mbox{fraction field of
$\Fbold$},\;\;\;
M:=F\otimes_{\Fbold}\Mbold,$$
$$\End_F(M)\supset R:=F\otimes_{\Fbold}{\Rbold}
\supset L:=F\otimes_{\Fbold}\Lbold=
\mbox{fraction field of
$\Lbold$},$$ 
$$K:=F\otimes_\Fbold\Kbold=\mbox{fraction field of
$\Kbold$}=\End_R(M)\subset \End_L(M)=L.$$
By Wedderburn's
theorem 
$R=\End_K(M)$ and hence
$${\Rbold}\subset R\cap \End_{\Fbold}(\Mbold)=
\End_K(M)\cap\End_{\Fbold}(\Mbold)=
\End_{\Kbold}(\Mbold).$$
Finally, we have 
$\Rbold=\End_{\Kbold}(\Mbold)$ because the
quotient
$\End_{\Kbold}(\Mbold)/\Rbold$ is a torsion $\Fbold$-submodule
of $\End_{\Fbold}(\Mbold)/\Rbold$
and hence by hypothesis vanishes.
\qed

\subsection{Geometric complex multiplication (GCM)}
\label{subsection:GCM}

\subsubsection{GCM $\FF_q[t]$-algebras}
Let $\Lbold$ be an $\FF_q[t]$-algebra satisfying the following
conditions: 
\begin{itemize}
\item $\Lbold$ is a free $\FF_q[t]$-module of finite rank.
\item $\Lbold$ is a Dedekind domain.
\item $\kbar\otimes_{\FF_q}\Lbold$ is also a Dedekind domain.
\end{itemize}
Under these conditions we say  that $\Lbold$ is a {\em GCM}
$\FF_q[t]$-algebra. ({\em GCM} is short for {\em
geometric complex multiplications}.) By
Proposition~\ref{Proposition:JustChecking} every
$\FF_q[t]$-subalgebra $\Kbold\subset\Lbold$ such that $\Kbold$ is a
Dedekind domain is again a GCM
$\FF_q[t]$-algebra.

\subsubsection{Constructions functorial in GCM $\FF_q[t]$-algebras}
For any GCM $\FF_q[t]$-algebra $\Lbold$ put
$$\Lboldbar:=\kbar\otimes_{\FF_q}\Lbold.$$
For any integer $n\in \ZZ$ let
$$\left(a\mapsto a^{(n)}\right):\Lboldbar\rightarrow\Lboldbar$$
be the unique $\Lbold$-linear extension of the
automorphism 
$$(x\mapsto x^{q^n}):\kbar\rightarrow\kbar.$$
Let $\Lboldbar[\sigma]$  be the ring obtained by adjoining
a noncommutative variable $\sigma$ to $\Lboldbar$ subject to the relations
$$\sigma a=a^{(-1)}\sigma\;\;\;(a\in \Lboldbar).$$
Every element of $\Lboldbar[\sigma]$ has a unique presentation of the form
$$\sum_{i=0}^\infty a_i\sigma^i\;\;\;(a_i\in \Lboldbar,\;\;\;
a_i=0\;\mbox{for}\; i\gg 0)$$
and in terms of such presentations the multiplication law 
in $\Lboldbar[\sigma]$ takes the form
$$\textstyle\left(\sum_i a_i\sigma^i\right)
\left(\sum_j b_j\sigma^j\right)=
\sum_i\sum_j a_ib_j^{(-i)}\sigma^{i+j}.$$
Note that $\Lbold$ is 
contained in the center of $\Lboldbar[\sigma]$. 
If $\Lbold=\FF_q[t]$, then
$\Lboldbar[\sigma]=\kbar[t,\sigma]$.

\subsubsection{Geometric complex multiplications}
Let $\Lbold$ be a GCM $\FF_q[t]$-algebra.
Let $H$ be a left $\Lboldbar[\sigma]$-module such that
the $\Lboldbar$-module underlying $H$ is projective of rank one 
and the
$\kbar[t,\sigma]$-module underlying $H$ is a dual
$t$-motive.
We call $H$ a dual $t$-motive with {\em geometric
complex multiplications} by
$\Lbold$ (for short: {\em GCM} by $\Lbold$).
We call the $\kbar[t,\sigma]$-module underlying $H$
the {\em bare} dual $t$-motive underlying $H$. We define the {\em GCM type}
of
$H$ with respect to $\Lbold$ to be the ideal of $\Lboldbar$ annihilating
the quotient
$H/\sigma H$. If the natural map
$\Lbold\rightarrow\End_{\kbar[t,\sigma]}(H)$
is bijective, we say that $H$ has {\em tight} GCM by $\Lbold$.
We
say that dual
$t$-motives
$H_0$ and
$H_1$  both with GCM by $\Lbold$
are {\em $\Lbold$-linearly isogenous} if there exists
an injective homomorphism
 $H_0\rightarrow H_1$
of  $\Lboldbar[\sigma]$-modules
with cokernel of finite dimension over $\kbar$.
By Theorem~\ref{Theorem:IsogenySymmetry} the relation of $\Lbold$-linear
isogeneity is an equivalence relation
and by
 Lemma~\ref{Lemma:IsogenyTangentBehavior}
the GCM type of a dual $t$-motive with GCM by $\Lbold$
depends only on its $\Lbold$-linear
isogeny class.

\begin{Theorem}\label{Theorem:SimpleGCM}
Let $H$ be a dual $t$-motive with GCM by $\Lbold$.
Let $H_0$ be any simple dual $t$-motive embedded in the bare dual
$t$-motive underlying $H$, e.g., any nonzero saturated
$\kbar[t,\sigma]$-submodule of minimum possible rank over $\kbar[t]$. (i)
The bare dual $t$-motive underlying $H$ is isogenous to a finite
direct sum of copies of the dual $t$-motive $H_0$.
(ii) If $H$ has tight GCM by $\Lbold$ then the bare dual $t$-motive
underlying $H$ is simple.
\end{Theorem}
\proof 
(i) Consider the family $$\{aH_0\}_{0\neq a\in \Lbold}$$
of isomorphic copies of $H_0$ embedded
$\kbar[t,\sigma]$-linearly in $H$. 
The sum $\sum_a aH_0$ is a nonzero
$\Lboldbar[\sigma]$-submodule of $H$ 
and {\em a fortiori} a $\kbar[t,\sigma]$-submodule of $H$ of finite
codimension over $\kbar$. Therefore the bare dual $t$-motive underlying
$H$ is isogenous to a direct sum of simple dual $t$-motives
isomorphic to $H_0$ by
Proposition~\ref{Proposition:SimpleDecomposition}. 

(ii) Making use of
Proposition~\ref{Proposition:SimpleDecomposition}  in a more precise way,
we obtain a positive integer $n$ and
$a_1,\dots,a_n\in \Lbold$ such that the map
$$\phi:=((h_1,\dots,h_n)\mapsto
a_1h_1+\cdots+a_nh_n):H_0^n\rightarrow H$$ is an isogeny of (bare) dual
$t$-motives.  By
Theorem~\ref{Theorem:2} there exists 
$0\neq a\in
\FF_q[t]$ such that $aH$ is contained in the image of $\phi$.
Put\\
$$\left.\begin{array}{rcl}
\pi_i&:=&((h_1,\dots,h_n)\mapsto
h_i):H^n_0\rightarrow H_0\\\\
e_i&:=&\left(h\mapsto
a_i\pi_i(\phi^{-1}(ah))\right)\in
\End_{\kbar[t,\sigma]}(H)
\end{array}\right\}\;\;\;(i=1,\dots,n).$$\\ The endomorphisms of 
the bare dual $t$-motive underlying $H$
thus constructed satisfy the relations
$$e_i\neq 0,\;\;\;\;e_ie_j=a\delta_{ij}e_i\;\;\;(i,j=1,\dots,n).$$
But unless $n=1$,  such a system of relations is forbidden
because the ring
 $\End_{\kbar[t,\sigma]}(H)$ is isomorphic
as an $\FF_q[t]$-algebra to the domain $\Lbold$.
\qed

\begin{Theorem}\label{Theorem:LiningUp}
For $i=0,1$ let $H_i$ be a dual $t$-motive with tight GCM by
$\Lbold_i$ and let $I_i\subset\Lboldbar_i$ be the GCM type of $H_i$
with respect to $\Lbold_i$. Assume
that the bare dual $t$-motives underlying
$H_0$ and
$H_1$ are isogenous. 
(i) Then
there exists a unique
$\FF_q[t]$-algebra isomorphism
$$\theta:\Lbold_0\iso\Lbold_1$$
such that all the diagrams
$$\begin{array}{rcl}
H_0&\xrightarrow{\phi}&H_1\\
{\scriptstyle a\times}\downarrow&&\downarrow{\scriptstyle
\theta(a)\times}\\ H_0&\xrightarrow{\phi}&H_1
\end{array}\;\;\;
\left(\phi\in \Hom_{\kbar[t,\sigma]}(H_0,H_1),\;\;a\in
\Lbold_0\right)
$$
commute. (ii) Under the unique $\kbar$-linear isomorphism
$\Lboldbar_0\rightarrow\Lboldbar_1$ induced by $\theta$
the ideal $I_0$ maps bijectively to ideal $I_1$.
\end{Theorem} 
\proof
(i) Put
$\Mbold:=\Hom_{\kbar[t,\sigma]}(H_0,H_1)$,
regarding $\Mbold$ as a right $\Lbold_0$-module and a left
$\Lbold_1$-module in the evident fashion.
Since $H_0$, $H_1$ have tight GCM by $\Lbold_0$, $\Lbold_1$,
Corollary~\ref{Corollary:RankIsogenyInvariance} shows that the
$\FF_q[t]$-modules $\Mbold$, $\Lbold_0$ and $\Lbold_1$ are free of the
same finite rank over $\FF_q[t]$.
It follows by Lemma~\ref{Lemma:Descent}
that $\Mbold$ is projective of rank one both over $\Lbold_0$
and over $\Lbold_1$.  Existence and uniqueness of $\theta$
now follow by Proposition~\ref{Proposition:SeesawTwo}.
(ii) For any isogeny $\phi:H_0\rightarrow H_1$ the induced
map $\phi\bmod{\sigma}$ is bijective by
Lemma~\ref{Lemma:IsogenyTangentBehavior}, whence the result.
\qed

\begin{Theorem}\label{Theorem:Decomposition}
Let $H$ be a dual $t$-motive with GCM by $\Lbold$.
Assume that the bare dual $t$-motive underlying $H$ is rigid analytically
trivial. Put
$$ \Rbold:=\End_{\kbar[t,\sigma]}(H),\;\;\;
\Kbold:=\{x\in \Lbold\mid (x\mapsto xh)\in
\End_\Rbold(H)\}.$$
Then $\Kbold$ is a GCM $\FF_q[t]$-subalgebra of $\Lbold$
and there exists a
$\Kboldbar[\sigma]$-submodule 
$$H_0\subset H$$ with the
following properties:
\begin{itemize}
\item  $H_0$ is a dual $t$-motive with tight GCM by $\Kbold$.
\item  
The bare dual
$t$-motive underlying
$H_0$ is simple and rigid analytically trivial. 
\item The GCM
type of
$H$ with respect to
$\Lbold$ is generated as an ideal of $\Lboldbar$ by the
GCM type of
$H_0$ with respect to
$\Kbold$.
\end{itemize}
\end{Theorem}
\proof We claim that the natural map
$$\Rbold=\End_{\kbar[t,\sigma]}(H)
\rightarrow\End_\Kbold\left(H^{\Betti}\right)$$
is bijective.  In any case the
natural map
$$\End_{\kbar[t,\sigma]}(H)\rightarrow\End_{\FF_q[t]}\left(H^{\Betti}
\right)$$
is
injective and has cokernel without $\FF_q[t]$-torsion by
Theorem~\ref{Theorem:RATConsequences}. Moreover, 
we have
$$\rank_{\FF_q[t]}\Lbold=\rank_{\kbar[t]}\Lboldbar=
\rank_{\kbar[t]}H=\rank_{\FF_q[t]}H^{\Betti}$$
where the first equality is trivial, the second
holds by
definition of GCM, and third
holds by
Lemma~\ref{Lemma:RATEquation} and hypothesis. 
Therefore the Dedekind-Wedderburn trick proves the claim.
It follows that $\Kbold$ is a GCM $\FF_q[t]$-subalgebra of $\Lbold$.

Let $n$ be the rank of $H^{\Betti}$ as a projective
$\Kbold$-module and let
$$H^{\Betti}=\bigoplus_{i=0}^{n-1} \Mbold_i$$
be a decomposition of $H^{\Betti}$ as a direct sum of projective
$\Kbold$-modules of rank one.
Let $e_0\in \Rbold$ be the idempotent endomorphism of $H$
inducing the projection of $H^{\Betti}$ to the direct summand $\Mbold_0$,
and put 
$$H_0:=e_0H.$$
Then $H_0$ is a rigid analytically trivial dual $t$-motive since it is a
$\kbar[t,\sigma]$-linear direct summand of a rigid analytically trivial
dual $t$-motive. Moreover, we have
$$H_0^{\Betti}=\Mbold_0,\;\;\;
\rank_{\kbar[t]}\Kboldbar=\rank_{\FF_q[t]}\Kbold=\rank_{\FF_q[t]}\Mbold_0=\rank_{\kbar[t]}H_0,$$
the last equality above by Lemma~\ref{Lemma:RATEquation}.
The inclusion $\Kbold \hookrightarrow \End_{\Kbold}(H^{\Betti})$
induces a $\Kboldbar$-module structure on $H_0$, and therefore $H_0$
has GCM by $\Kbold$.  By 
a repetition of
the argument of the first paragraph of the
proof
we have a natural bijective map
$$\End_{\kbar[t,\sigma]}(H_0)\rightarrow
\End_{\Kbold}(\Mbold_0)=\Kbold.$$
Therefore $H_0$ is naturally equipped with tight GCM by
$\Kbold$.
Moreover $H_0$ is automatically simple as a bare dual $t$-motive by 
Theorem~\ref{Theorem:SimpleGCM}.  The natural map
$$(a\otimes h\mapsto ah):\Lboldbar\otimes_{\Kboldbar} H_0\rightarrow H$$
is injective with  $\Lboldbar[\sigma]$-stable image. Call the image
$H'$. The $\Lboldbar$-module underlying $H'$ is projective of
rank one and {\em a fortiori} of finite codimension over
$\kbar$ in $H$. Clearly the $\kbar[t]$- and $\kbar[\sigma]$-modules
underlying $H'$ are free of finite rank. Since
$H'$ is generated over
$\Lbold$ by
$H_0$ and $\Lbold$ is central in $\Lboldbar[\sigma]$,
we have 
$$(t-T)^nH_0\subset \sigma H_0\Rightarrow (t-T)^nH'\subset \sigma H'$$
for all integers $n\geq 0$ and hence the $\kbar[t,\sigma]$-module
underlying
$H'$ is a dual $t$-motive. The upshot is that $H'$ is a dual
$t$-motive with GCM by $\Lbold$ and that $H'$ is $\Lbold$-linearly
isogenous to
$H$ via the inclusion. Necessarily the GCM types of
$H'$ and $H$ with respect to $\Lbold$ are equal. Finally, we have
at our disposal an isomorphism
$$\Lboldbar\otimes_{\Kboldbar}\frac{H_0}{\sigma H_0}
=\frac{H'}{\sigma H'}$$
of torsion $\Lboldbar$-modules, whence it follows that the GCM type
of $H'$ with respect to $\Lbold$ is generated as an ideal of $\Lboldbar$
by the GCM type of $H_0$ with respect to $\Kbold$. 
\qed

\begin{Corollary}\label{Corollary:SimplicityCriterion}
Let $\Lbold$ be a GCM $\FF_q[t]$-algebra.
Assume that the fraction field of $\Lbold$ is Galois over $\FF_q(t)$
and let
$G$ be the associated Galois group.  Extend the action of $G$ on $\Lbold$
to
$\Lboldbar$ by $\kbar$-linearity. Let
$H$ be a  dual
$t$-motive with GCM by
$\Lbold$  
and 
with rigid analytically trivial underlying bare dual
$t$-motive. Let
$H_0$ be a simple quotient of the bare dual
$t$-motive underlying $H$. Let
$I\subset\Lboldbar$ be the GCM type of
$H$ with respect to $\Lbold$. Let $r$ be the cardinality of the set of
ideals of $\Lboldbar$ that are $G$-conjugate to $I$. Then $H_0$ is of
rank
$\geq r$ over
$\kbar[t]$.
\end{Corollary}
\proof  By
Theorem~\ref{Theorem:Decomposition}  there exist a
GCM $\FF_q[t]$-subalgebra
$\Kbold\subset\Lbold$ and a
 $\Kboldbar[\sigma]$-submodule
$H'_0\subset H$
with the following properties:
\begin{itemize}
\item $H'_0$ has tight GCM by $\Kbold$. 
\item The bare dual $t$-motive underlying $H'_0$ is simple. 
\item The GCM
type
$I_0\subset\Kboldbar$ of $H'_0$ with respect to $\Kbold$ generates
$I\subset\Lboldbar$.
\end{itemize}
Further, by Theorem~\ref{Theorem:SimpleGCM} 
the given quotient $H_0$
of the bare dual
$t$-motive
underlying $H$
and the bare dual $t$-motive underlying $H'_0$ are
isogenous. Finally, we have
$$\rank_{\kbar[t]}H_0=\rank_{\kbar[t]}H'_0=\rank_{\FF_q[t]}\Kbold\geq r$$
because
$\gamma I=I$ for every $\gamma\in G$ acting as the identity on $\Kbold$.
\qed

\begin{Corollary}\label{Corollary:NonisogenyCriterion}
Let $\Lbold$ be a GCM $\FF_q[t]$-algebra.
Assume that the fraction field of $\Lbold$ is Galois over $\FF_q(t)$
and let
$G$ be the associated Galois group.  Extend the action of $G$ on $\Lbold$
to
$\Lboldbar$ by $\kbar$-linearity.  For $i=0,1$
let
$H_i$ be a dual
$t$-motive with GCM by
$\Lbold$ with 
rigid analytically trivial underlying bare dual $t$-motive,  let
$I_i\subset\Lboldbar$ be the GCM type of
$H_i$ with respect to
$\Lbold$, and let $H_{i0}$ be a simple quotient of the bare dual
$t$-motive underlying $H_i$. If $H_{00}$ and $H_{10}$ are isogenous, then
the ideals
$I_0$ and
$I_1$ are
$G$-conjugate.
\end{Corollary}
\proof For $i=0,1$,   there exist by
Theorem~\ref{Theorem:Decomposition} a
GCM $\FF_q[t]$-subalgebra
$\Kbold_i\subset\Lbold$ and a $\Kboldbar_i[\sigma]$-submodule
$H'_{i0}\subset H_i$
with the following properties:
\begin{itemize}
\item $H'_{i0}$ has tight GCM by $\Kbold_i$. 
\item The bare dual $t$-motive underlying $H'_{i0}$ is simple. 
\item The GCM
type
$I_{i0}\subset\Kboldbar_i$ of $H'_{i0}$ with respect to $\Kbold_i$
generates
$I_i\subset\Lboldbar$.
\end{itemize}
Further, for $i=0,1$ by Theorem~\ref{Theorem:SimpleGCM} 
the bare dual $t$-motive underlying $H'_{i0}$ is isogenous to $H_{i0}$,
and hence the bare dual $t$-motives underlying $H'_{00}$ and $H'_{10}$
are isogenous.
Finally, by Theorem~\ref{Theorem:LiningUp} there exists some $\gamma\in G$
such that $\gamma
\Kbold_0=\Kbold_1$ and
$\gamma I_{00}=I_{10}$, in which case necessarily $\gamma I_0=I_1$.
\qed
\subsubsection{Remark}
The theory of rigid analytically trivial GCM dual $t$-motives
worked out above is essentially just the dual
$t$-motivic translation of a method introduced in
\cite{BrownPapa} for analyzing the $t$-submodule structure of
the geometric
$t$-modules of \cite{Sinha}.

\section{Special functions}
\label{section:Ragbag}

\subsection{The Carlitz exponential and its fundamental period}   
\label{subsection:CarlitzExponential}
 The reference \cite[Chap.~3]{Goss} is a
good source of background material on this topic.

\subsubsection{The Carlitz
exponential}\label{subsubsection:ExpDef} Put 
$$\exp_C z:=\sum_{n=0}^\infty
\frac{z^{q^n}}{D_n}\;\;\;\left(D_n:=\prod_{i=0}^{n-1}
\left(T^{q^n}-T^{q^i}\right)\right)$$
thereby defining an
$\FF_q$-linear power series in $z$ with coefficients in $k$ called
the {\em Carlitz exponential}. 
   The Carlitz exponential satisfies the functional equation
$$ \exp_C (Tz)=T\exp_C z+(\exp_C z)^q$$
by \cite[Prop.~3.3.1]{Goss}.

\subsubsection{The fundamental period}
The power series
$\exp_C z$ has an infinite radius of convergence with respect to the
valuation
$|\cdot|_\infty$ and has a Weierstrass product expansion of
the form 
$$\exp_C z=z\prod_{
\begin{subarray}{c}
a\in A\\
a\neq 0\end{subarray}}\left(1-\frac{z}{\varpi
a}\right)$$
for unique 
$\varpi\in k_\infty(\Tbar)$ such that
$$|\varpi-\Tbar T|_\infty<|\Tbar T|_\infty$$ 
by \cite[Cor.~3.2.9 and Remarks~3.2.10]{Goss}. It follows that the
sequences
$$0\rightarrow \varpi\cdot
A\subset\overline{k_\infty}\xrightarrow{\exp_C}
\overline{k_\infty}\rightarrow 0,\;\;\;\;0\rightarrow \varpi\cdot
A\subset\CC_\infty
\xrightarrow{\exp_C}\CC_\infty\rightarrow 0$$
are exact. 
We
call
$\varpi$ the {\em fundamental period} of the Carlitz exponential. 
Transcendence of
$\varpi$ over
$k$ was first proved in
\cite{Wade}. 

The next two results relate the Carlitz exponential and
its fundamental period to the power series $\Omega(t)$ and
transcendental number $\Omega(T)$ discussed in
\S\ref{subsubsection:TranscendenceIllustration}.

\begin{Proposition}\label{Proposition:EboldOmegaConnection}
We
have
$$
1/\Omega^{(-1)}(t)=\sum_{i=0}^\infty
\exp_C\left(\varpi/T^{i+1}\right)t^i.$$
\end{Proposition}
\proof Recall that
$$\Omega(t)=\Tbar^{-q}\prod_{i=1}^\infty\left(1-\frac{t}{T^{q^i}}\right),\;\;\;
(t-T)\cdot\Omega=\Omega^{(-1)}.$$
Temporarily denote the power series on the right side of the
identity to be proved by
$\Theta(t)$.  
Consider
the Maclaurin expansion
$$\Omega^{(-1)}(t)\Theta(t)=\sum_{i=0}^\infty
c_it^i\;\;\;\;\;(c_i\in k_\infty(\Tbar)).$$
The functional equation noted in
\S\ref{subsubsection:ExpDef} implies that
$$t\Theta=T\Theta+\Theta^{(1)},$$
and hence
$$
\Omega\Theta^{(1)}=
\Omega\cdot
(t-T)\cdot \Theta=
\Omega^{(-1)}\Theta=\left(\Omega\Theta^{(1)}\right)^{(-1)}.$$  By
this last we have
$c_i=\sqrt[q]{c_i}$ and hence
$c_i\in \FF_q$ for all $i\geq 0$.
By plugging into the
Weierstrass product
expansion of 
$\exp_Cz$ we find that 
$|c_0-1|<1$ and hence 
$$c_0=\Tbar^{-1}\cdot\exp_C\left(\varpi/T\right)=1.$$
Now write
$$\Omega^{(-1)}(t)=\sum_{i=0}^\infty
a_it^i,\;\;\;\Theta(t)=\sum_{i=0}^\infty
b_it^i\;\;\;(a_i,b_i\in k_\infty(\Tbar)).$$ We have, so we claim,
bounds
$$|a_i|_\infty\leq |\Tbar|^{-1},\;\;\;|b_i|_\infty\leq
|\Tbar|$$
with strict  inequality for $i>0$. The bound for
$|a_i|_\infty$ is clear. The bound for $|b_i|_\infty$ is
verified by plugging into the Weierstrass product 
expansion of $\exp_C z$. Thus the claim is proved. It follows that
$|c_i|_\infty<1$ for all $i>0$, hence 
$c_i=0$ for all $i>0$, and hence $\Omega^{(-1)}\Theta=1$.
\qed

\begin{Corollary} \label{Corollary:CarlitzPeriodProduct}
We have
$$\varpi=T\Tbar\prod_{i=1}^\infty
\left(1-T^{1-q^i}\right)^{-1}=-1/\Omega(T).$$
\end{Corollary}
\proof (Cf.~\cite[Corollary~2.5.8]{AndersonThakur}.)  By
Proposition~\ref{Proposition:EboldOmegaConnection} and a little high
school algebra (summation of geometric series), we have
$$\frac{\frac{1}{\Omega(t)}+\varpi}{t-T}=\frac{1}{\Omega^{(-1)}(t)}-
\sum_{n=0}^\infty\left(
\frac{\varpi }{T^{n+1}}\right)t^n
=\sum_{n=0}^\infty\left(\sum_{i=1}^\infty\frac{1}{D_i}
\left(\frac{\varpi}{T^{n+1}}\right)^{q^i}
\right)t^n.$$ 
The power series on the right is convergent in the disc
$|t|_\infty<|T|^q_\infty$ and so is the power series
$1/\Omega(t)$.  Therefore we have
$1/\Omega(T)+\varpi=0$.
\qed
 
\subsection{The functions $\ebold$ and $\ebold^*$ and their division
polynomials}
\subsubsection{Definition of
$\ebold$}\label{subsubsection:eboldDef}
 Put
$$\ebold:=(x\mapsto \exp_C \varpi
x):k_\infty\rightarrow k_\infty(\Tbar),$$
where $\varpi$ is the fundamental period of the Carlitz
exponential. The function $\ebold$ maps the quotient
$k_\infty/A$ isomorphically to a compact additive
subgroup of $k_\infty(\Tbar)$. The function
$\ebold$ is in many respects analogous to the function
$(x\mapsto e^{2\pi i x}):\RR\rightarrow\CC$
mapping the quotient $\RR/\ZZ$ 
isomorphically to the unit circle in $\CC$.

\subsubsection{Division polynomials for $\ebold$}
Given $a\in A$,
write
$$a=Tb+\epsilon\;\;\;(b\in A,\;\;\epsilon\in \FF_q)$$
and put
$$C_a(t,z):=\left\{\begin{array}{cl}
0&\mbox{if $a=0$,}\\
C_b(t,tz+z^q)+\epsilon z&
\mbox{if $a\neq 0$,}
\end{array}\right.
$$
thereby recursively defining a polynomial
$$C_a(t,z)\in \FF_q[t,z].$$ 
For example we have
$$C_1(t,z)=z,\;\;\;C_T(t,z)=tz+z^q,
\;\;\;C_{T^2}(t,z)=t^2z+(t+t^q)z^q+z^{q^2}.$$
For all $f=f(T)\in A_+$ of positive degree we have\\
$$
C_f(t,z)=f(t)z+ \left(\begin{array}{l}
\textup{\mbox{$\FF_q[t]$-linear
combination of}}\\
\textup{\mbox{terms
$z^{q^i}$ with
$0<i<\deg f$}}
\end{array}\right)+z^{q^{\deg
f}}
$$\\
as can be verified by an evident induction. 
{From} the functional equation noted in \S\ref{subsubsection:ExpDef}
it follows that
$$T\ebold(x)+\ebold(x)^q=\ebold(Tx)$$
for all $x\in k_\infty$ and hence by an evident induction that
$$C_a\left(T,\ebold(x)\right)=\ebold(ax)$$
for all $a\in A$ and $x\in k_\infty$.
{From} the latter identity it follows  that $C_a(t,z)$
depends
$\FF_q$-linearly on
$a$ and that the {\em composition law}
$$C_a(t,C_b(t,z))=C_{ab}(t,z)$$
holds for all $a,b\in A$. 
We call
$C_a(t,z)$ the {\em division polynomial} for $\ebold$ indexed by $a$.

\subsubsection{Torsion values of
$\ebold$}\label{subsubsection:DelicateEbold} A
number of the form
$$\ebold(x)\;\;\;(x\in k)$$
will be called a {\em torsion value} of
$\ebold$. For every $f\in A_+$ we have
$$C_f(T,z)=\prod_{\begin{subarray}{c}
a\in A\\
\deg a<\deg f
\end{subarray}}(z-\ebold(a/f))$$
and hence every torsion value of $\ebold$
is separably algebraic over
$k$.  
\subsubsection{Remark}
The ring homomorphism
$$(a\mapsto (x\mapsto
C_a(T,x))):A\rightarrow\End_{\mbox{\scriptsize alg.~gp.}/A}({\mathbf
G}_a)$$ 
is called the {\em Carlitz
module} according to \cite[Def.~3.3.5]{Goss}.
The torsion values of $\ebold$ admit interpretation as the torsion
points of the Carlitz module defined over $\kbar$. The latter
interpretation makes it clear that torsion values of $\ebold$ generate
abelian extensions of
$k$. The Carlitz module plays in explicit class field theory over
$k$ a role analogous to that played
by the multiplicative group in explicit class field theory over
$\QQ$. See \cite{Hayes} for a treatment of explicit class field
theory over $k$. We consider some of the more delicate properties of
torsion values of $\ebold$ below in
\S\ref{subsubsection:TorsionCloserLook}.

\subsubsection{Definition of
$\ebold^*$} \label{subsubsection:eboldstarDef}
Let  
$$\Res:k_\infty\rightarrow \FF_q$$ be the unique $\FF_q$-linear
functional such that
$$\ker\Res=\FF_q[T]+(1/T^2)\FF_q[\![1/T]\!],\;\;\;\;
\Res(1/T)=1.$$
Write
$$\Omega^{(-1)}(t)=
\Tbar^{-1}\prod_{i=0}^\infty\left(1-\frac{t}{T^{q^i}}\right)=\sum_{i=0}^\infty a_it^i\;\;\;(a_i\in
k_\infty(\Tbar)).$$ 
For each
$x\in k_\infty$ put
$$\ebold^*(x):=\sum_{i=0}^\infty \Res(T^ix)a_i,$$
e.~g.,
$$\ebold^*(1/T^{n+1})=\left\{\begin{array}{cl}
a_n&\mbox{if $n\geq 0$}\\
0&\mbox{if $n<0$}
\end{array}\right.
$$
for any $n\in \ZZ$,
thereby defining a
function
$$\ebold^*:k_\infty\rightarrow k_\infty(\Tbar).$$
A number of the form
$$\ebold^*(x)\;\;\;(x\in k)$$
will be called a {\em torsion value} of $\ebold^*$.

\begin{Lemma}\label{Lemma:EboldStarEstimate}
We have $\ebold^*(x)\neq 0$ for 
all $x\in k_\infty$ such that $0<|x|_\infty<1$.
\end{Lemma}
\proof With the
coefficients
$a_i\in
k_\infty(\Tbar)$ as defined in \S\ref{subsubsection:eboldstarDef},
we have
$$\left|a_i-(-T)^{-\left(q^0+\dots+q^{i-1}\right)}\Tbar^{-1}\right|_\infty=
|a_i-\Tbar^{-q^i}|_\infty<|\Tbar^{-q^i}|_\infty,$$
for all integers $i\geq 0$ and hence
$$|\ebold^*(x)|_\infty=\left|\Tbar^{-q^{\min\{i\vert\Res(T^ix)\neq
0\}}}\right|_\infty>0$$
for all $x\in k_\infty$ such that $0<|x|_\infty<1$. 
\qed

\subsubsection{Division polynomials for $\ebold^*$}
Put 
$$\overline{C}_1(t,z):=z.$$
For each $f\in A_+$ of positive degree, write
$$f=Tg+\epsilon\;\;\;(g\in A_+,\;\;\epsilon\in\FF_q)$$
and put
$$\overline{C}_f(t,z):=
\overline{C}_g(t^q,tz^q+z)+\epsilon z^{q^{\deg f}}.
$$
In this way we recursively  define polynomials
$$\overline{C}_f(t,z)\in \FF_q[t,z]$$
for all $f\in A_+$.
For
example we have
$$\overline{C}_1(T,z)=z,\;\;\;\overline{C}_T(t,z)=tz^q+z,
\;\;\;\overline{C}_{T^2}(t,z)=t^{2q}z^{q^2}+(t+t^q)z^q+z.$$
For all $f=f(T)\in A_+$ of positive degree we have\\
$$\overline{C}_f(t,z)=z+ \left(\begin{array}{l}
\textup{\mbox{$\FF_q[t]$-linear
combination of}}\\
\textup{\mbox{terms
$z^{q^i}$ with
$0<i<\deg f$}}
\end{array}\right)+f(t)^{q^{-1+\deg f}}z^{q^{\deg
f}}
$$
as can be verified by an evident induction.
With the
coefficients
$a_i\in
k_\infty(\Tbar)$ as defined in \S\ref{subsubsection:eboldstarDef}
and in view
of the functional equation 
$\Omega^{(-1)}=(t-T)\cdot\Omega$
we have
$$
Ta_i^q+a_i=\left\{\begin{array}{cl}
a_{i-1}^q&\mbox{if $i>0$}\\
0&\mbox{if $i=0$}
\end{array}\right.$$
for all nonnegative integers $i$ and hence
$$T\ebold^*(x)^q+\ebold^*(x)=\ebold^*(Tx)^q$$ 
for all $x\in k_\infty$. It follows by an evident induction that
$$\overline{C}_f(T,\ebold(x))=\ebold(fx)^{q^{\deg f}}$$
for all $f\in A_+$ and $x\in k_\infty$. This last is the end to
which the definition of
$\overline{C}_f(T,z)$ is contrived.
We call $\overline{C}_f(t,z)$ the {\em division polynomial} for
$\ebold^*$ indexed by $f$.

\subsubsection{Key properties of $\ebold^*$}
\label{subsubsection:EboldStarKey}
The function
$$\ebold^*:k_\infty\rightarrow k_\infty(\Tbar)$$ has, so we claim,
the following properties:\\
\begin{itemize}
\item
$\sum_{i=0}^\infty
\ebold^*\left(1/T^{i+1}\right)t^i=
\Tbar^{-1}\prod_{i=0}^\infty\left(1-t/T^{q^i}\right)=
\Omega^{(-1)}(t).$
\item $\ebold^*$ is $\FF_q$-linear and $|\cdot|_\infty$-continuous.
\item $\ker\ebold^*=A$.
\item $\ebold^*(x)\Tbar\in \FF_q[\![1/T]\!]$ for all $x\in
k_\infty$.
\item
The torsion values of
$\ebold^*$ are separably algebraic over $k$.\\
\end{itemize}
The function $\ebold^*$ has the first two of the
claimed properties by definition.
We have $A\subset\ker\ebold^*$ by definition
and
$\ker\ebold^*\subset A$ by Lemma~\ref{Lemma:EboldStarEstimate}.
Therefore the third property holds. The first three properties imply
the fourth.
 Each torsion value of $\ebold^*$ is a solution of an equation of the
form $\overline{C}_f(T,z)=0$ for some $f\in A_+$ and hence
separably algebraic over $k$. Therefore the fifth property holds.
Thus all claims concerning
$\ebold^*$ are proved. Note that the first three properties above
already determine the function $\ebold^*$ uniquely.

\subsubsection{Remark}
Given $f\in A_+$, write
$$C_f(t,z)=\sum_{i=0}^{n} f_i
z^{q^i}\;\;\;(f_i\in
\FF_q[t],\;\;\;n=\deg f).$$
It can be shown that
$$\overline{C}_f(t,z)=
\sum_{i=0}^{n}f_{i}^{q^{n-i-1}}z^{q^{n-i}}.$$
See \cite[\S3.7 and \S4.14]{Goss} for an attractive
interpretation of  this phenomenon.  It
turns out that the
torsion values of the function
$\ebold^*$ admit interpretation as the torsion points defined over
$\kbar$ of a gadget called the {\em
adjoint Carlitz module} and that the latter is dual to the Carlitz
module in a sense made precise by a duality theory due to Elkies and
Poonen. 

\subsubsection{Key properties of $\ebold$}
\label{subsubsection:EboldKey}
The properties of the function $\ebold^*$ are strikingly parallel to
those of the function
$$\ebold:k_\infty\rightarrow k_\infty(\Tbar)$$ as is made plain by
the following summary of facts already proved:
\\
\begin{itemize}
\item $\sum_{i=0}^\infty
\ebold\left(1/T^{i+1}\right)t^i=\Tbar\prod_{i=0}^\infty
\left(1-t/T^{q^i}\right)^{-1}=1/\Omega^{(-1)}(t)$.
\item $\ebold$ is
$\FF_q$-linear and $|\cdot|_\infty$-continuous.
\item $\ker\ebold=A$.
\item $\ebold(x)/\Tbar\in\FF_q[\![1/T]\!]$ for all $x\in k_\infty$.
\item The torsion values of $\ebold$ are separably algebraic over
$k$.\\
\end{itemize}
Note that the first three
of the properties listed above determine $\ebold$ uniquely 
without any reference to the Carlitz exponential or its fundamental
period.

\subsubsection{Remark}\label{subsubsection:DigitPatterns}
For each 
$n\in
\ZZ_{\geq 0}$
let
$$n=\sum_i n_iq^i\;\;\;(0\leq n_i<q)$$
be the base $q$ representation of $n$
and put
$$\alpha(n):=
\left\{\begin{array}{cl}
\sum_i n_i&\mbox{if $n_i\in \{0,1\}$ for all $i$,}\\
-\infty&\mbox{otherwise.}
\end{array}\right.
$$
By definition of $\ebold^*$ we have
$$\Tbar\ebold^*(x)=\sum_{n=0}^\infty
\Res\left(\left(-T\right)^{\alpha(n)}x\right)T^{-n}\;\;\;\;\;\left(\mbox{by
convention:}\;(-T)^{-\infty}=0\right)$$ for all $x\in k_\infty$. Note
that the very simple ``program''\\
$$\alpha(n)=\left\{\begin{array}{cl}
0&\mbox{if $n=0$, else}\\\\
\displaystyle \alpha(n/q)&\mbox{if
$n\equiv 0\bmod{q}$, else}\\\\
\displaystyle \alpha((n-1)/q)+1&
\mbox{if $n\equiv 1\bmod{q}$, else}\\\\
-\infty&
\end{array}\right.
$$\\
computes $\alpha(n)$ recursively.  Now fix $x\in k$. Then the
sequence
$$\left\{\Res\left(\left(-T\right)^n x\right)\right\}_{n=0}^\infty$$
is eventually periodic, hence  the sequence
$$\left\{\Res\left((-T)^{\alpha(n)}x\right)\right\}_{n=0}^\infty$$ 
of coefficients is computable by
a ``$q$-automaton'' and hence the number $\Tbar\ebold^*(x)$ is
algebraic over $k$ by Christol's algebraicity
criterion \cite{Christol} (see also
 \cite{CKMR}). Thus we can give by
``computer science'' a proof of the algebraicity of the torsion
values of $\ebold^*$ over $k$ completely independent of the {\em ad
hoc} theory of division polynomials worked out above. The simplicity
of the ``digit-patterns'' of the torsion points of the adjoint
Carlitz module seems not previously to have been noticed.  This
remark will not be used in the sequel but does we think deserve
further investigation.

\subsection{The geometric $\Gamma$- and $\Pi$-functions and their
special values}
\subsubsection{The Moore determinant identity}
Given an $\FF_q$-algebra $R$ and $x_1,\dots,x_N\in R$,
 put
$$\Moore_q(x_1,\dots,x_N):=\det_{i,j=1}^N x_j^{q^{N-i}}=\left|\begin{array}{lcl}
x_1^{q^{N-1}}&\dots&x_N^{q^{N-1}}\\
\vdots&&\vdots\\
x_1^{q^0}&\dots&x_N^{q^0}\end{array}\right|\in R$$
thereby defining the 
{\em Moore determinant} of $x_1,\dots,x_N$.
The {\em Moore determinant identity} reads
$$\Moore_q(x_1,\dots,x_N)=
\prod_{
\begin{subarray}{c}
c\in \FF_q^N\\
c:\mbox{\scriptsize monic}
\end{subarray}}\left(\sum_{i=1}^N
c_ix_i\right)$$ where {\em ad hoc} we say
that a vector $c=(c_1,\dots,c_N)\in
\FF_q^N$ is {\em monic} if $c\neq 0$ and the leftmost
nonzero entry of $c$ equals
$1$.  See
\cite[Chap.~1, Sec.~3]{Goss} for  proof and further discussion of
the identity.

\subsubsection{Definition of $\Psi_N(z)$}\label{subsubsection:PsiDef}
For
each integer $N\geq 0$, put
$$\Psi_N(z):=\left(\prod_{\begin{subarray}{c}
a\in A\\
\deg a<N\end{subarray}}(z-a)\right){\bigg /}
\left(\prod_{\begin{subarray}{c}
a\in A_+\\
\deg a=N\end{subarray}}a\right)$$
thereby defining an $\FF_q$-linear polynomial in $z$ with coefficients in $k$.
We have
$$\Psi_N(z)=
\frac{\Moore_q\left(z,T^{N-1},\dots,1\right)}{\Moore_q\left(T^N,\dots,1\right)}
= \frac{z^{q^N}}{D_N}+\left(\mbox{terms of degree $<q^N$ in $z$}\right)
$$
by combining the Moore and Vandermonde determinant identities. 
The definition of
$\Psi_N(z)$ goes back to the paper \cite{CarlitzAncient}.  
A contemporary reference for this material is \cite[\S3.5]{Goss}.

\subsubsection{Key relations satisfied by $\Psi_N(z)$}
\label{subsubsection:KeyPsi} Fix an integer $N\geq 0$.
We have
$$1+ \Psi_N(z)=\prod_{\begin{subarray}{c}
a\in A_+\\
\deg a=N
\end{subarray}}\left(1+\frac{z}{a}\right)
$$
as can be verified by comparing zeroes and constant terms.
Clearly we have
$$\deg a<N\Rightarrow \Psi_N(z+a)=\Psi_N(z)$$
for all $a\in A$.
We have
$$N>0\Rightarrow\Psi_{N}(z)=\frac{\Psi_{N-1}(z)^q-\Psi_{N-1}(z)}{T^{q^{N}}-T}$$
as can be verified by comparing zeroes
 and leading terms; for convenient reference we dub
this important fact the {\em fundamental recursion} for
$\Psi_N(z)$. For any 
$f\in A_+$ we have
$$\prod_{\begin{subarray}{c} a\in
A\\
\deg a<\deg f\end{subarray}}
\frac{1+\Psi_N\left(\frac{z+a}{f}\right)}
{1+\Psi_N\left(\frac{a}{f}\right)}=1+\Psi_{N+\deg
f}(z)$$
as can be verified by comparing zeroes and constant terms.

\subsubsection{Definitions of  $\Pi(z)$ and $\Gamma(z)$}
\label{Subsubsection:DefinitionsOfGammaAndPi}
Put
$$
\Pi(z):=\prod_{a\in
A_+}\left(1+\frac{z}{a}\right)^{-1},\;\;\;
\Gamma(z):=z^{-1}\Pi(z)=z^{-1}\prod_{a\in
A_+}\left(1+\frac{z}{a}\right)^{-1},$$
thereby defining the (one-variable) {\em geometric factorial}
and {\em
geometric
$\Gamma$-function}  attached to
$A$, respectively. The basic references concerning these functions
are
\cite[Sec.~9.9]{Goss} and
\cite{ThakurGamma}. Note that
$\Pi(z)^{-1}$ is the unique entire function of $z$ 
taking the value $1$ at the origin
with no zeroes other than simple zeroes at each point of the set
$$-A_+:=\{-a\vert a\in A_+\}.$$
Similarly  $\Gamma(z)^{-1}$ is the unique entire
function of $z$ normalized by the condition
$$\Gamma(z)^{-1}=z+O(z^2)$$
with no zeroes other than simple zeroes at each point of the set
$$-A_+\cup\{0\}.$$
 We obtain identities
$$\Pi(z)=\prod_{N=0}^\infty(1+\Psi_N(z))^{-1},\;\;\;\;
\Gamma(z)=z^{-1}\prod_{N=0}^\infty(1+\Psi_N(z))^{-1}
$$
of crucial importance in the
sequel by grouping factors in the natural way by degree.

\subsubsection{The standard functional equations}
\label{subsubsection:KeyPi}
 We have
a {\em translation identity}
$$\frac{\Pi(z+a_0)}{\Pi(z)}=\prod_{i=0}^{\deg a_0}
\frac{1+\Psi_i(z)}{1+\Psi_i(z+a_0)}$$
for each $0\neq a_0\in A$. The
translation identity follows in evident fashion from the corresponding
translation-invariance property of $\Psi_N(z)$ noted in
\S\ref{subsubsection:KeyPsi}.  We have a {\em reflection identity}
$$\prod_{\epsilon\in \FF_q^\times}\Pi(\epsilon z)
=\frac{\varpi z}{\exp_C\varpi z}$$
which  is essentially \cite[Theorem 6.1.1]{ThakurGamma}.  To
prove the reflection identity in the form stated here you have only to
compare the Weierstrass product expansion of the reciprocal of the left
side to that of the Carlitz exponential
$\exp_C z$. 
For each
$f\in A_+$ we have a {\em Gauss multiplication identity}
$$\prod_{\begin{subarray}{c}
a\in A\\
\deg a<\deg f
\end{subarray}}\Pi\left(\frac{z+a}{f}\right)
=\Pi(z)\cdot\prod_{0\leq i<\deg f} (1+\Psi_i(z))
\cdot \prod_{\begin{subarray}{c}
a\in A_+\\
\deg a<\deg f
\end{subarray}}\prod_{\epsilon\in \FF_q^\times}\Pi\left(\epsilon
\frac{a}{f}\right)$$  which is essentially
\cite[Theorem 6.2.1]{ThakurGamma}. The Gauss multiplication
identity in the form stated here can easily be recovered from the
analogous identity satisfied by $\Psi_N(z)$ noted at 
the end of \S\ref{subsubsection:KeyPsi}. We refer to the identities
above as the {\em standard functional equations} satisfied by
$\Pi(z)$. Of course we have analogous standard functional equations
for $\Gamma(z)$. The latter we
do not write out, referring the interested reader to \cite{ThakurGamma}.

\subsubsection{Special $\Gamma$- and $\Pi$-values}
Numbers of the form
$$\Pi(x)\;\;\;(x\in k\setminus -A_+)$$ 
will be called {\em special $\Pi$-values}. 
Numbers of the form
$$\Gamma(x)\;\;\;(x\in k\setminus(-A_+\cup\{0\}))$$
will be called {\em special $\Gamma$-values}. 
The goal of the paper is to determine all Laurent polynomial relations
among $\varpi$ and special $\Gamma$-values with coefficients in $\kbar$.
{From} the point of view of transcendence theory over $k$ it is clearly
all the same whether we study special
$\Gamma$-values or special $\Pi$-values. But the
standard functional equations satisfied by
$\Pi(z)$ take a slightly simpler form than do those satisfied by
$\Gamma(z)$ and it is convenient also that $\Pi(0)=1$. 
Accordingly, for the sake of convenience and for no deeper reason,
we stress
$\Pi(z)$ over
$\Gamma(z)$ and special $\Pi$-values over special 
$\Gamma$-values in the sequel. 
Now in view of the translation identity
of \S\ref{subsubsection:KeyPi}, a special  $\Pi$-value $\Pi(x)$
up  to factors in $k^\times$ depends only on $x\bmod{A}$.
Accordingly, without any loss of generality, we stress special
$\Pi$-values  of the form $\Pi(x)$
with $|x|_\infty<1$ in the sequel.

\subsubsection{Remark}
All $\kbar$-linear relations among $1$, $\varpi$
and the special $\Pi$-values were determined in
\cite{BrownPapa}, building on the work of \cite{Sinha}.
In particular, it is known that $\Pi(z)$ is transcendental
for all $z\in k\setminus A$.

\subsection{Interpolation formulas}
\label{subsection:InterpolationFormulas}
\begin{Lemma}\label{Lemma:OmegaPsiRelation}
For all integers $N$ the identity
$$\Omega^{(N)}/\Omega^{(-1)}=
\left\{\begin{array}{cl}
\displaystyle-\sum_{n=0}^\infty
\Psi_N\left(1/T^{n+1}\right)t^n&\mbox{if $N\geq 0$}\\\\
\displaystyle\prod_{i=1}^{|N|-1}\left(t-T^{q^{-i}}\right)&\mbox{if
$N<0$}
\end{array}
\right.
$$
holds in the power series ring $\kbar[\![t]\!]$.
\end{Lemma}
\proof  {From} the functional equation 
$\Omega^{(-1)}=\Omega\cdot(t-T)$ it follows that
$$\Omega^{(N)}/\Omega^{(-1)}=
\left\{\begin{array}{cl}
\displaystyle\prod_{i=0}^N\left(t-T^{q^i}\right)^{-1}&\mbox{if $N\geq
0$,}\\\\
\displaystyle\prod_{i=1}^{|N|-1}\left(t-T^{q^{-i}}\right)&\mbox{if
$N<0$.}
\end{array}
\right.
$$
Accordingly, we may assume without loss of generality for
the rest of the proof that $N\geq 0$. The identity
$$(t-T)^{-1}=-\sum_{n=0}^\infty
t^n/T^{n+1}= -\sum_{n=0}^\infty
\Psi_0\left(1/T^{n+1}\right)t^n$$ dispatches the case $N=0$. Assume
now that
$N>0$. The recursion
$$
\left(\prod_{i=0}^{N-1}
\left(t-T^{q^{i}}\right)^{-1}\right)^{(1)}-\prod_{i=0}^{N-1}
\left(t-T^{q^{i}}\right)^{-1}=
(T^{q^{N}}-T)\prod_{i=0}^{N}\left(t-T^{q^i}\right)^{-1}$$
matches the fundamental recursion 
$$\Psi_{N-1}(z)^q-\Psi_{N-1}(z)=\left(T^{q^{N}}-T\right)\Psi_{N}(z)$$
noted in \S\ref{subsubsection:KeyPsi}. We are done by induction on
$N$.
\qed

\begin{Lemma}\label{Lemma:InterpolationIntermediary1}
For all $x\in k_\infty$ and integers $N$ such that
$$|x|_\infty
\leq |T|_\infty^{\min(-1,N)}$$
we have
$$\sum_{i=0}^\infty
\ebold^*(T^{-i-1})^{q^{N+1}}\ebold(T^i
x)=\left\{\begin{array}{rl}
-\Psi_N(x)&\mbox{if $N\geq 0$,}\\\\
\Res\left(T^{-N-1}x\right)&\mbox{if $N<0$.}
\end{array}\right.
$$\\
\end{Lemma}
\proof 
Both sides of the identity to be proved depend
$\FF_q$-linearly and
 $|\cdot|_\infty$-continuously on $x$.
We may therefore assume without loss of generality that $x=T^{-n-1}$
for some integer $n$ such that 
$-n-1\leq \min(-1,N)$.
By Lemma~\ref{Lemma:OmegaPsiRelation},
the formal properties of $\ebold^*$
noted in
\S\ref{subsubsection:EboldStarKey}, and the formal properties of $\ebold$
noted in \S\ref{subsubsection:EboldKey}, both sides of the claimed
identity admit interpretation as  the coefficient with which
$t^n$ appears in  the Maclaurin expansion of $\Omega^{(N)}/\Omega^{(-1)}$
in powers of
$t$.
\qed

\subsubsection{$f$-dual families}
Fix $f\in A_+$. We say that families
$$\left\{a_i\right\}_{i=1}^{\deg
f},\;\;\;\left\{b_j\right\}_{j=1}^{\deg f}\;\;\; (a_i,b_j\in A)$$
are {\em $f$-dual} if
$$\Res(a_ib_j/f)=\delta_{ij}\;\;\;\;(i,j=1,\dots,\deg f).$$
Since the square matrix
$$\left\{\Res(T^{i+\deg f-j-1}/f)\right\}_{i,j=0}^{\deg f-1}$$
is lower triangular with $1$'s along the diagonal,
the pairing
$$\left((a\bmod{f},b\bmod{f})\mapsto \Res(ab/f)\right):A/f\times
A/f\rightarrow \FF_q$$
is perfect and hence $f$-dual
families exist.

\begin{Theorem}\label{Theorem:Interpolation}
Fix $f\in A_+$ and $f$-dual
families 
$$\left\{a_i\right\}_{i=1}^{\deg
f},\;\;\;\left\{b_j\right\}_{j=1}^{\deg f}\;\;\; (a_i,b_j\in A).$$
Fix $a\in A$ such that
$$\deg a<\deg f.$$
(i) We have
$$
\sum_{i=1}^{\deg f}\ebold^*(a_i/f)^{q^{N+1}}\ebold(b_ia/f)
=
-\Psi_N(a/f)
$$ for all integers $N\geq 0$. (ii) Moreover, if $a\in A_+$ we have
$$\sum_{i=1}^{\deg f}\ebold^*(a_i/f)
\ebold(b_ia/f)^{q^{\deg f-\deg a-1}}=1.
$$
\end{Theorem}
\proof (Cf.~\cite[Thm.~2, p.~58]{AndersonStick}.)
For every $N\in \ZZ$ we have
$$\begin{array}{cl}
&\displaystyle \sum_{i=1}^{\deg f}
\ebold^*(a_i/f)^{q^{N+1}}\ebold(b_i a/f)\\\\
=&\displaystyle
\sum_{i=1}^{\deg f}\sum_{n=0}^\infty\ebold^*(1/T^{n+1})^{q^{N+1}}
\Res(T^na_i/f)\ebold(b_i
a/f)\\\\ =&\displaystyle\sum_{n=0}^\infty
\ebold^*(1/T^{n+1})^{q^{N+1}}\ebold(T^n a/f)\\\\
=&\left\{\begin{array}{cl}
\Psi_N(a/f)&\mbox{if $N\geq 0$}\\
1&\mbox{if $a\in A_+$ and $N=\deg a-\deg f$}
\end{array}\right.
\end{array}
$$
by definition of the special function $\ebold^*$ and
Lemma~\ref{Lemma:InterpolationIntermediary1}.
\qed
\subsubsection{Remark}
Theorem~\ref{Theorem:Interpolation}
``interpolates'' $\Psi_N(a/f)$ by an algebraic expression in which $N$
figures as the power to which a Frobenius endomorphism is raised.  We
learned this seemingly strange but in fact fundamental notion of
interpolation from examples in
\cite[\S9.3]{ThakurGamma}; see \cite[\S2]{AndersonStick}
for an appreciation of Thakur's work.
The possibility of such an interpolation was proved in
\cite[Thm.~2, p.~58]{AndersonStick} without the
interpolating expression being made explicit.

\subsection{Diamond brackets and $L$-functions}
\subsubsection{Definitions}\label{subsubsection:DiamondBrackets}
For all $x\in k_\infty$
and integers $N\geq 0$, put\\
$$\diamondbracket{x}_N:=\left\{\begin{array}{rl} 1&\mbox{if
$\displaystyle 0<\inf_{a\in
\FF_q[T]}\left|x-a-T^{-N-1}\right|_\infty<|T|_\infty^{-N-1}$,}\\\\
0&\mbox{otherwise,}
\end{array}\right.
$$
and also set
$$\diamondbracket{x}:=\sum_{N=0}^\infty \diamondbracket{x}_N.$$
The sum on the right makes sense because at most one of its terms is
nonzero. We call the functions
$$\diamondbracket{\cdot},\diamondbracket{\cdot}_N:k_\infty\rightarrow\{0,1\}$$
thus defined the {\em
diamond bracket} and the {\em generalized diamond bracket},
respectively. 
\subsubsection{Remark}
The value assigned to the expression $\diamondbracket{x}$ here
coincides with that assigned in
\cite{SinhaDR}, but with that assigned to
$\diamondbracket{-x}$ 
in \cite[Def.~7.6.1]{ThakurGamma}. The choice
of sign is just a normalization not affecting the utility of the
definition.

\subsubsection{Evaluation of $L$-functions}
\label{subsubsection:Lvalues}
Fix $f\in A_+$ of positive degree and a character $\chi:(A/f)^\times
\rightarrow
\CC^\times$ not factoring through $(A/g)^\times$ for any proper divisor
$g\in A_+$ of $f$. Note that under these hypotheses
$\chi$ is not identically equal to $1$. Put
$$\begin{array}{rcl}
L(s,{\chi})&=&\displaystyle \sum_{\begin{subarray}{c}
a\in A_+\\
(a,f)=1
\end{subarray}}{\chi(a\bmod{f})} q^{-s\deg a}\\\\
&=&\displaystyle \sum_{\begin{subarray}{c}
a\in A_+\\
(a,f)=1\\
\deg a<\deg f
\end{subarray}}{\chi(a\bmod{f})} q^{-s\deg a}\\\\
&=&\displaystyle\sum_{\begin{subarray}{c}
a\in A\\
\deg a<\deg f\\
(a,f)=1
\end{subarray}}
{\chi(a\bmod{f})}\left(\sum_{N=0}^\infty
\diamondbracket{\frac{a}{f}}_{N}q^{(N+1-\deg
f)s}\right),
\end{array}$$
noting  in particular that
$$L(0,\chi)=\sum_{\begin{subarray}{c}
a\in A\\
\deg a<\deg f\\
(a,f)=1
\end{subarray}}\chi(a\bmod{f})\diamondbracket{\frac{a}{f}}.$$
The $L$-function $L(s,\chi)$ is the $\FF_q[T]$-analogue of a Dirichlet
$L$-function. 
It is well known that $L(0,\chi)=0$ if and only if $\chi(\epsilon
\bmod{f})=1$ for all
$\epsilon\in
\FF_q^\times$. 
\subsubsection{Remark}
The connection between $L$-functions and diamond brackets
recalled above is the motivation for the definition of diamond brackets.

\subsubsection{Diamond bracket relations}
\label{subsubsection:DiamondBracketRelations}
We make the following
claims:\\
\begin{itemize}
\item $\diamondbracket{x+a}=\diamondbracket{x}$
for all $x\in k_\infty$ and $a\in A$. \\\\
\item $\displaystyle\sum_{\epsilon\in
\FF_q^\times}\diamondbracket{\epsilon x}=\left\{\begin{array}{cl}
1&\mbox{if $x\not\in A$}\\
0&\mbox{if $x\in A$}\end{array}\right.$
for all $x\in k_\infty$. \\\\
\item $\displaystyle\sum_{\begin{subarray}{c}
a\in A\\
\deg a<\deg f
\end{subarray}}\left(
\diamondbracket{\frac{x+a}{f}}-\diamondbracket{\frac{a}{f}}\right)=
\diamondbracket{x}$
for all $x\in k_\infty$ and $f\in A_+$. \\\\
\end{itemize}
Clearly the first and second claims hold. 
To prove the third claim we may assume without loss of
generality that $|x|_\infty<1$. Then the only summand on the
left possibly differing from $0$ is the one indexed by
$a=0$, and we have 
$\diamondbracket{x/f}=\diamondbracket{x}$. Therefore
the third claim  holds.

\section{Analysis of the algebraic relations among special
$\Pi$-values}
\label{section:Finale}
Throughout \S\ref{section:Finale} we fix $f\in A_+$ of positive
degree.

\subsection{The diamond bracket criterion}
\label{subsection:DiamondBracketYoga}

\subsubsection{The free abelian group $\AAA_f$}
\label{subsubsection:AAAFormalism}

Let
$\AAA_f$ be the free abelian group on symbols of the form
$$[x]\;\;\;\left(x\in f^{-1}A\right),$$
where symbols $[x]$ and $[x']$ are identified
if $x\equiv x'\bmod{A}$.
The group $\AAA_f$ is free abelian
of rank $q^{\deg f}$. 
Let $\DDD_f$ be the subgroup of $\AAA_f$ generated by all elements of the form
$$\left[x\right]-\sum_{\begin{subarray}{c}
a\in A\\
\deg a<\deg g
\end{subarray}}\left[\frac{x+a}{g}\right]\;\;\;\left(g\in A_+
\mbox{ dividing}\;f,\;\;\;x\in \frac{g}{f}\cdot A\right).$$
The quotient $\AAA_f/\DDD_f$ is the analogue over $\FF_q[T]$ of the {\em
universal ordinary distribution} of level $f$, cf.~\cite{Kubert} or
\cite[\S2]{LangCyclotomic}.
  Let $\RRR_f$
be the subgroup of
$\AAA_f$ generated by $\DDD_f$ along with all elements of the form
$$\sum_{\epsilon\in\FF_q^\times}\left[\epsilon x\right]\;\;\;\left(x\in
f^{-1}A\right).$$
Let $\widetilde{\RRR}_f$ be the subgroup of $\AAA_f$ 
consisting of $\abold$ such that $\weight\abold=0$ and
$N\abold\in \RRR_f$ for some positive integer $N$.

\subsubsection{The star action}
For each $a\in A$ prime to $f$ there exists a unique automorphism
$$(\abold\mapsto a\star \abold):\AAA_f\rightarrow \AAA_f$$
of free abelian groups
such that
$$a\star[x]=[ax]$$
for all $x\in \frac{1}{f}A$.   
Note that this automorphism stabilizes $\DDD_f$,
$\RRR_f$ and
$\widetilde{\RRR}_f$. 
Note that
$$a\star(b\star \abold)=(ab)\star \abold$$
for all $a,b\in A$ prime to $f$ and $\abold\in \AAA_f$.
Thus via the star operation $\AAA_f$
is equipped with an action of
$(A/f)^\times$ passing to the quotients $\AAA_f/\DDD_f$, $\AAA_f/\RRR_f$
and
$\AAA_f/\widetilde{\RRR}_f$. 

\begin{Theorem}\label{Theorem:Kubert}
The rational vector space
$$\QQ\otimes(\AAA_f/\DDD_f)$$ is
$(A/f)^\times$-equivariantly isomorphic to the
rational group ring 
$$\QQ\left[(A/f)^\times\right].$$
\end{Theorem}
\proof The classical model for the theorem
is proved in \cite{Kubert}; alternatively,
see\cite[\S2]{LangCyclotomic}.  The methods of Kubert carry over to our
function field situation without any difficulty. We omit the details.
\qed

\begin{Corollary}\label{Corollary:TranscendenceDegreeDetermination}
Fix a character
$\chi:(A/f)^\times\rightarrow\CC^\times$.
The
dimension over
$\CC$ of the $\chi$-isotypical component of 
$\CC\otimes (\AAA_f/\widetilde{\RRR}_f) $  is $\leq 1$,
with strict inequality if and only if $\chi$ is not identically equal to
$1$ but $\chi\left(\epsilon\bmod{f}\right)=1$ for all $\epsilon\in
\FF_q^\times$. 
\end{Corollary}
It follows that the free abelian group
$\AAA_f/\widetilde{\RRR}_f$ is
of rank $1+\frac{q-2}{q-1}\cdot\card(A/f)^\times$.
\proof We have
$$\AAA_f/\RRR_f= 
\mbox{module of $\FF_q^\times$-coinvariants
in $\AAA_f/\DDD_f$}.$$
The sequence of
$(A/f)^\times$-modules
$$0\rightarrow \CC\cdot \left(\sum_{\begin{subarray}{c}
 0\neq a\in A\\
\deg a<\deg f
\end{subarray}}\left[\frac{a}{f}\right]\bmod{\widetilde{\RRR}_f}
\right)
\rightarrow  \CC\otimes (\AAA_f/\widetilde{\RRR}_f)
\rightarrow \CC\otimes(\AAA_f/\RRR_f)\rightarrow 0$$
is exact. The result follows.
\qed
\subsubsection{Extension of the definition of diamond brackets}
 For all integers $N\geq 0$ let
$$(\abold\mapsto\diamondbracket{\abold}_N):\AAA_f\rightarrow \ZZ$$
be the unique homomorphism such that
$$\diamondbracket{[x]}_N=\diamondbracket{x}_N$$
for all $x\in \frac{1}{f}A$
 and put
$$\diamondbracket{\abold}=\sum_{N=0}^\infty \diamondbracket{\abold}_N$$
for all $\abold\in \AAA_f$. The sum on the right makes sense because only
finitely many of its terms are nonzero. For all $\abold,\bbold\in \AAA_f$
we write 
$\abold\sim_f\bbold$
if and only if 
$\diamondbracket{a\star\abold}=
\diamondbracket{a\star\bbold}$ for all $a\in A$ 
such that $(a,f)=1$ and $\deg a<\deg f$,
thereby defining an equivalence relation $\sim_f$
in $\AAA_f$.
When $\sim_f$ fails to hold between $\abold$
and $\bbold$ we write $\abold\not\sim_f\bbold$.

\begin{Theorem}\label{Theorem:KoblitzOgus}
For all $\abold,\bbold\in
\AAA_f$ we have 
$\abold\equiv\bbold\bmod{\widetilde{\RRR}_f}$
if and only if $\abold\sim_f\bbold$.
\end{Theorem}
\proof 
 ($\Rightarrow$) Since $\widetilde{\RRR}_f$ is stable under
the action of
$(A/f)^\times$ it is enough to prove
$$\abold\in \widetilde{\RRR}_f\Rightarrow
\diamondbracket{\abold}=0.$$
In turn it is enough to prove that
$$\abold\in \RRR_f\Rightarrow \diamondbracket{\abold}=\weight
\abold.$$
In order to prove the latter implication we may assume
without loss of generality that $\abold$ is one of the generators of
$\RRR_f$ exhibited in \S\ref{subsubsection:AAAFormalism}.  Then the
diamond bracket relations of
\S\ref{subsubsection:DiamondBracketRelations} do the job.

($\Leftarrow$)  
This is proved by an evident modification of the
proof of the Koblitz-Ogus criterion \cite{DeligneCorvallis}.
The connection between diamond brackets and 
values of $L$-functions at $s=0$ noted in
\S\ref{subsubsection:Lvalues} is the essential point of the proof. 
 We
omit the details.
\qed

\subsubsection{$\Pi$-monomials and their relationship with
$\Gamma$-monomials} Let
$$(\abold\mapsto \Pi(\abold)):\AAA_f\rightarrow \CC_\infty^\times$$
be the unique homomorphism such that
$$\Pi([x])=\Pi(x)$$
for all $x\in \frac{1}{f}A$ such that $|x|_\infty<1$.
Numbers in the image of the
homomorphism $\AAA_f\xrightarrow{\Pi}\CC_\infty^\times$ 
we call {\em $\Pi$-monomials} of level $f$.
The reflection
identity satisfied by the $\Pi$-function and the hypothesis $\deg
f>0$ imply that up to a factor in
$\kbar^\times$, the number $\varpi$ is a $\Pi$-monomial of level $f$.
Taking the translation identity satisfied by the $\Pi$-function also into
account, as well as the simple relationship between $\Gamma$-
and $\Pi$-functions, it is clear
that  up to a factor in
$\kbar^\times$ every
$\Gamma$-monomial belonging to the group of such generated
by the set 
$$\{\varpi\}\cup\left\{\Gamma(x)\left|
x\in \frac{1}{f}A\setminus(\{0\}\cup-A_{+})\right.\right\}$$
is a $\Pi$-monomial of level $f$.

The following result is the direct analogue in our setting of the
Koblitz-Ogus criterion stated in \cite{DeligneCorvallis}. It allows us to
decide in  more or less mechanical fashion whether between a given pair
of 
$\Gamma$-monomials there exists a relation of $\kbar$-linear dependence
explained by the standard functional equations.

\begin{Corollary}[Diamond bracket criterion]
\label{Corollary:KoblitzOgus} For all $\abold,\bbold\in \AAA_f$ we have
$$\abold\sim_f\bbold\Rightarrow\Pi(\abold)/\Pi(\bbold)\in
\kbar^\times.$$
\end{Corollary}

\proof It is enough to prove that 
$$\abold\in \RRR_f\Rightarrow\varpi^{-\weight
\abold}\;
\Pi(\abold)\in \kbar^\times.$$ In order to do so we may assume 
without loss of generality that $\abold$ is one of the generators of 
$\RRR_f$
specified in
\S\ref{subsubsection:AAAFormalism}. Then the standard
functional equations stated in
\S\ref{subsubsection:KeyPi} do the job.
\qed

\subsubsection{Remark}
Deligne's reciprocity law
\cite[Thm.~7.15,  p.~91]{DMOS} refines the
Koblitz-Ogus criterion by giving information concerning
the field to which an algebraic $\Gamma$-monomial belongs.
The analogous refinement of Corollary~\ref{Corollary:KoblitzOgus}
was proved in \cite{SinhaDR}.

\subsubsection{Remark}
The converse to Corollary~\ref{Corollary:KoblitzOgus}
is the $\FF_q[T]$-analogue ``at level $f$'' of the
conjecture of Rohrlich discussed in 
\S\ref{subsubsection:RohrlichConjecture}.

\subsection{Formulation and discussion of the main result}
\label{subsection:MainResultFormulation}
The following theorem is the main result of this paper.
It is the $\FF_q[T]$-analogue  ``at level $f$''
of the conjecture of
Lang discussed in \S\ref{subsubsection:LangConjecture} above. It
restates Theorem~\ref{Theorem:LinearRelationsIntro} in a precise way and
allows us to determine all
$\kbar$-linear relations among $\Gamma$-monomials.

\begin{Theorem}\label{Theorem:LinearRelations}
Let 
\[
\abold_1,\dots,\abold_N\in \AAA_f
\]
 be given.  A necessary and sufficient condition for the corresponding
$\Pi$-monomials
\[
\Pi(\abold_1),\dots,\Pi(\abold_N)
\]
of level $f$ to be $\kbar$-linearly independent
is that 
\[
  \abold_i \not\sim_f \abold_j 
\]
for all $1 \le i < j \le N$.
\end{Theorem}
The proof of the theorem takes up almost all of the rest of the
paper, concluding in \S\ref{subsection:LooseEnds}.

\begin{Corollary}\label{Corollary:TranscendenceDegree}
Put
$$\nu_f:=1+\frac{q-2}{q-1}\cdot\card(A/f)^\times.$$
Let $E_f$ be the subfield of
$\CC_\infty$ generated over $\kbar$ by the set
$$\left\{\Pi(a/f)\left|a\in A,\;\deg a<\deg f\right.\right\}.$$
The transcendence degree of $E_f$ over $\kbar$ equals $\nu_f$.
\end{Corollary}
The following proposition establishes equivalence of theorem and
corollary. 
\begin{Proposition}\label{Proposition:Strategy}
The transcendence degree of $E_f/\kbar$ is bounded above by $\nu_f$.
Moreover, a necessary and sufficient condition
for $\nu_f$ strictly to  exceed the transcendence degree of $E_f/\kbar$ is
that for some integer 
$N\geq 2$ there exist
$\abold_1,\dots,\abold_N\in\AAA_f$ and
nonzero $c_1,\dots,c_N\in\kbar$
such that
$\abold_i\not\sim_f\abold_j$
for all $1\leq i<j\leq N$
and 
$c_1\Pi(\abold_1)+\cdots+c_N\Pi(\abold_N)=0$.
\end{Proposition}
\proof
Let $\kbar[\AAA_f/\widetilde{\RRR}_f]$ be the group
ring with coefficients in $\kbar$ of  the finitely generated free
abelian group $\AAA_f/\widetilde{\RRR}_f$.  Let 
$$\lambda:\AAA_f\rightarrow \kbar^\times$$ be any group homomorphism
agreeing with $\Pi$ on $\widetilde{\RRR}_f$
and let
$$\kbar[\AAA_f/\widetilde{\RRR}_f]\xrightarrow{\Lambda}
E_f$$ be the unique $\kbar$-algebra homomorphism
such that
$$\Lambda
(\abold\bmod{\tilde{\RRR}_f})=\frac{\Pi(\abold)}{\lambda(\abold)}$$
for all $\abold\in \AAA_f$. Put
$$I_f:=
\ker\left(\kbar[\AAA_f/\widetilde{\RRR}_f]\xrightarrow{\Lambda}E_f\right).$$
By
Corollary~\ref{Corollary:TranscendenceDegreeDetermination} 
the free abelian group $\AAA_f/\widetilde{\RRR}_f$ is of rank
$\nu_f$ and hence the
ring
$\kbar[\AAA_f/\widetilde{\RRR}_f]$ is isomorphic to the ring
of Laurent polynomials in $\nu_f$ independent variables with
coefficients in
$\kbar$. By construction the ideal $I_f$ is prime.
Clearly the field
$E_f$ is isomorphic as a $\kbar$-algebra to the field of fractions of the
ring
$\kbar[\AAA_f/\widetilde{\RRR}_f]/I_f$.
Therefore $\nu_f$ bounds the transcendence degree of $E_f$
over $\kbar$. Moreover, a necessary and sufficient condition for
$\nu_f$ strictly to  exceed the transcendence degree of $E_f/\kbar$ is
that $I_f\neq 0$. Note that every nonzero element of $I_f$
has to be a formal $\kbar$-linear combination of at least two
elements of
$\AAA_f/\widetilde{\RRR}_f$ since $\Pi(\abold)$ is non-zero for all
$\abold \in \AAA_f$. 
In view of Theorem~\ref{Theorem:KoblitzOgus}
and Corollary~\ref{Corollary:KoblitzOgus}, it is clear that 
from
any nonzero element of
$I_f$ we can  produce a
relation of
$\kbar$-linear dependence among
$\Pi$-monomials of the indicated form.  \qed

\subsection{Coleman functions} 
\label{subsection:ColemanConstruction}
\subsubsection{A closer look at torsion values of $\ebold$}
\label{subsubsection:TorsionCloserLook}
In the course of our
discussion of the special function $\ebold$ the following facts
were verified:
\begin{itemize}
\item $k(\ebold(1/f))$ is
separable over
$k$.
\item $k(\ebold(1/f))$ is a splitting field over $k$ for $C_f(T,z)$.
\item $\ebold(a/f)/\Tbar\in \FF_q[\![1/T]\!]$ for all $a\in A$.
\end{itemize} 
The following more delicate facts are well known:
\begin{itemize}
\item There exists for each $\gamma\in \Gal(k(\ebold(1/f))/k)$ 
unique $a\in A$ such that $(a,f)=1$, $\deg a<\deg f$ and
$\gamma\ebold(b/f)=\ebold(ab/f)=C_a(T,\ebold(b/f))$ for all $b\in A$.
\item The construction $\gamma\mapsto a$
induces an isomorphism
$\Gal(k(\ebold(1/f))/k)\iso (A/f)^\times$.
\item $\Gal(k(\ebold(1/f))/k)$ is generated by
its inertia subgroups.
\item The integral closure of $A$ in $k(\ebold(1/f))$
is $A[\ebold(1/f)]$.
\end{itemize}
See \cite{Hayes} for a treatment of the latter material. 

\subsubsection{Definition of $C_f^\star(t,z)$}
There exists a unique factor
$$C^\star_f(t,z)\in \FF_q[t,z]$$
of the division polynomial $C_f(t,z)$ such that
$$C^\star_f(T,z)=\prod_{\begin{subarray}{c}
a\in A\\
\deg a<\deg f\\
(a,f)=1
\end{subarray}}(z-\ebold(a/f)).$$
The polynomial $C_f^\star(t,z)$ is the Carlitz analogue of a cyclotomic
polynomial and has, so we claim, the following properties:
\begin{itemize}
\item The discriminant of $C_f^\star(t,z)$ with respect to $z$ does not
vanish identically.
\item For all $a\in A$, $C_f^\star(t,z)$ divides $C_a(t,z)$
if and only if $f$ divides $a$.
\item
$C_f^\star(t,z)$ is irreducible in $\FF_q[t,z]$ and remains
so in $\kbar[t,z]$.
\item 
$C_f^\star$, $\partial C_f^\star/\partial t$ and $\partial
C_f^\star/\partial z$ generate the unit ideal of $\FF_q[t,z]$.
\end{itemize} 
The first two properties are clear, as is irreducibility
 over $\FF_q$. Irreducibility over 
any finite algebraic extension of $\FF_q$  follows from the fact that
$\Gal(k(\ebold(1/f))/k)$ is generated by inertia.
Were $C_f^\star(t,z)$ to  be reducible over
$\kbar$, then $C_f^\star(t,z)$ would be reducible over some
field
$L$ finite algebraic over $k$ and hence reducible over the finite residue
field of some discrete valuation of $L$, a contradiction.  Therefore
$C_f^\star(t,z)$ has the third property. 
Failure of the fourth
property would imply failure of the ring
$\FF_q[t,z]/(C_f^\star(t,z))$ to be integrally closed:
but the isomorphic ring
$A[\ebold(1/f)]$ is known to be integrally closed.
Therefore $C_f^\star(t,z)$ has the fourth property. 
Thus our claim is proved.

\subsubsection{The nonsingular projective curve $X/\FF_q$}
Let $U/\FF_q$ be the irreducible nonsingular  plane algebraic
curve in the affine $(t,z)$-plane$/\FF_q$ defined by the equation
$C_f^\star(t,z)=0$ and let $X/\FF_q$ be the nonsingular projective model of
$U/\FF_q$. We regard $t$ and $z$ as regular functions on $U$ and
meromorphic functions on $X$. We regard  $X/\FF_q$ as a covering of
the projective
$t$-line$/\FF_q$. We say that the closed points of
$X$ in the complement of $U$ are {\em at infinity}.
For all
$a\in A$ prime to
$f$ put
$$\xi_a:=\left(T,\ebold(a/f)\right),$$
thereby defining a $\kbar$-valued point of $U$.
A trivial but important remark to make here is that the set $\{\xi_a\}$ is
the collection of
$\kbar$-valued points  of $U$ above the $\kbar$-valued point $t=T$ of
the affine
$t$-line. We make the following claims:
\begin{itemize}
\item $X$ is a Galois covering of the $t$-line
of degree $\card(A/f)^\times$.
\item For each automorphism $\gamma$ of $X$ over the
$t$-line there exists unique
$a\in A$ such that
$(a,f)=1$, $\deg a<\deg f$ and
$\gamma^*C_b(t,z)=C_{ab}(t,z)$ for all $b\in A$.
\item Moreover, with $\gamma$ and $a$ as above,
we have $\gamma\xi_b=\xi_{ab}$ for all $b\in A$ prime to $f$.
\item Further, the construction $\gamma\mapsto a$ induces an isomorphism
{from} the group of automorphisms of $X$ over the $t$-line to
$(A/f)^\times$.
\item
Each closed point of $X$ at infinity has residue field $\FF_q$
and is ramified of order $q-1$ over the point at infinity on the $t$-line.
\item There are $\card(A/f)^\times/(q-1)$ closed points at infinity,
and these are permuted transitively by the group of automorphisms of $X$
over the $t$-line.
\item
For all
$a\in A$ the function
$C_a(t,z)$ has at each of the points of
$X$ at infinity no singularity worse than a simple pole. 
\end{itemize}
These claims are verified by making a routine translation
{from} arithmetical language to geometrical language. We can safely omit
the details.

\subsubsection{The base-change $\Xbar/\kbar$,
deck transformations and 
$n$-fold twisting}
Put
$$\Ubar:=\kbar\otimes_{\FF_q}U,\;\;\;\;\Xbar:=\kbar\otimes_{\FF_q}X,$$
thereby defining nonsingular irreducible curves over $\kbar$,
the former being the curve in the affine $(t,z)$-plane$/\kbar$
defined by the equation $C_f^\star(t,z)=0$, and the latter
being the nonsingular projective model of the former.
Closed points of $\Xbar$ in the complement of
$\Ubar$ as before are said to be {\em at infinity}. No new closed
points at infinity appear in the base-change $\Xbar$ because all the closed
points in $X$ at infinity are already $\FF_q$-rational.

By construction
$\Xbar/\kbar$ is a Galois covering of the projective $t$-line$/\kbar$.
Just so as to have a convenient short turn of phrase at our disposal, we
call an automorphism of
$\Xbar/\kbar$ over the projective $t$-line$/\kbar$ a {\em deck
transformation}. By construction every deck transformation is the
base-change of a unique automorphism of
$X/\FF_q$ over the projective $t$-line$/\FF_q$, and hence the group of
deck transformations is canonically isomorphic to $(A/f)^\times$.

We define the {\em
$n$-fold twisting operation}
on the function
field of $\Xbar$ to be the unique automorphism
extending the $(q^n)^{th}$ power automorphism of $\kbar$
and fixing every element of the function field of $X$.
We denote the result of applying the $n$-fold twisting operation to a
function $h$ by $h^{(n)}$. The $n$-fold twisting automorphism of the
function field of
$\Xbar$ commutes with all pull-backs via deck transformations. For each
$x\in X(\kbar)$ we define the {\em
$n$-fold twist} $x^{(n)}\in X(\kbar)$ to be the point obtained by applying
the
$(q^n)^{th}$ power automorphism of $\kbar$ to the coordinates of $x$. 
Since each $\kbar$-rational point of $X$ at infinity is already defined
over $\FF_q$, each such point is fixed by the $n$-fold twisting operation.
Identifying closed points of
$\Xbar$ with
$\kbar$-valued points of
$X$ in evident fashion, we extend the $n$-fold twisting operation
to the group of divisors of $\Xbar$ by $\ZZ$-linearity. 
The $n$-fold twisting operation on divisors commutes with the action
of deck transformations.
The operation of forming the divisor of a nonzero meromorphic function on
$\Xbar$ commutes with
$n$-fold twisting, i.e., we have $(h^{(n)})=(h)^{(n)}$
for all nonzero meromorphic functions $h$ on $\Xbar$.

 \subsubsection{Definition of
Coleman functions}
\label{subsubsection:ColemanDefinition}
Fix
$$x\in f^{-1}A\setminus A.$$
Fix $f$-dual
families
$$\{a_i\}_{i=1}^{\deg f},\{b_j\}_{j=1}^{\deg f}\;\;\;(a_i,b_j\in A)$$
and write 
$$x=a_0/f\;\;\;(a_0\in A,\;\;\mbox{$a_0$ not divisible by $f$}).$$
Put
$$g_x:=1-\sum_{i=1}^{\deg f}
\ebold^*(a_i/f)C_{a_0b_i}(t,z),$$ thereby defining a
meromorphic function on $\Xbar$ regular on $\Ubar$ with singularities
at infinity no worse than simple poles. Now for
$a$ ranging over
$A$, both the function
$C_a(t,z)$ and the number $\ebold^*(a/f)$ depend
$\FF_q$-linearly on $a$ and moreover depend
only on $a\bmod{f}$. Therefore the
function
$g_x$ depends only on $x$, not on the intervening choice of
$f$-dual families $\{a_i\}$ and $\{b_j\}$. Moreover it is clear that
$g_x$ depends only on $x\bmod{A}$. We call $g_x$ a {\em Coleman
function}.

\subsubsection{The divisors of Coleman functions}
\label{subsubsection:ColemanDivisor}
Let $\infty_X$ be the formal sum of the
 $\kbar$-valued points of
$X$ at infinity,  multiplied by $(q-1)$ and viewed as a divisor of
$\Xbar$.  We have
$$\deg \infty_X=\card(A/f)^\times,\;\;\;(\infty_X)^{(1)}=\infty_X.$$
For every
$x\in f^{-1}A\setminus A$ we have, so we claim, an equality
$$(g_x)=-\frac{1}{q-1}\cdot\infty_X+\sum_{\begin{subarray}{c} a\in A\\
(a,f)=1\\
\deg a<\deg f
\end{subarray}}\sum_{N=0}^\infty\diamondbracket{ax}_{N}\cdot
\xi^{(N)}_a$$
of divisors of $\Xbar$. To see this, call the divisor on the right $D$.
There appear only finitely many nonzero terms
in the sum defining $D$ 
and hence $D$ is well-defined. Moreover we have
$$\deg D=-\frac{\card(A/f)^\times}{q-1}+
\sum_{\begin{subarray}{c}
a\in A\\ (a,f)=1\\
\deg a<\deg f
\end{subarray}}\diamondbracket{ax}=0.$$
Now let $\{a_i\}$, $\{b_j\}$ and $a_0$ be as in
\S\ref{subsubsection:ColemanDefinition}, and fix $a\in A$ prime to $f$ such
that
$$\diamondbracket{ax}=\sum_{N=0}^\infty\diamondbracket{ax}_N=1.$$
Let $b$ be the unique element of $A_+$ such that
$$aa_0\equiv b\mod{f},\;\;\;\;\deg b<\deg f$$
and put
$$N:=\deg f-\deg b-1,$$
noting that we have
$$\diamondbracket{b/f}_{N}=1.$$ We have
$$\begin{array}{rcl}
g_x\left(\xi^{(N)}_a\right)&=&\displaystyle
1-\sum_{i=1}^{\deg f}
\ebold^*(a_i/f)C_{a_0b_i}\left(T^{q^N},\ebold(ax)^{q^{N}}\right)\\\\
&=&\displaystyle
1-\sum_{i=1}^{\deg f}
\ebold^*(a_i/f)\ebold(b_ib/f)^{q^{N}}\\\\
&=&0,
\end{array}
$$\\
the last equality by part (ii) of Theorem~\ref{Theorem:Interpolation}.
Therefore $g_x$ has at least as many zeroes in $\Ubar$ as we claim for it.
In any case $g_x$ has no singularities at infinity worse than simple poles.
Therefore the divisor $(g_x)-D$ is effective
and of degree $0$, so it vanishes identically.
Thus our claim is proved.

\subsubsection{Interpolation properties of Coleman functions}
\label{subsubsection:ColemanInterpolation}
Fix 
$$x\in f^{-1}A\setminus A,\;\;\;a\in A,\;\;\;
N\in \ZZ,\;\;\;y\in k$$
such that
$$(a,f)=1,\;\;\;N\geq 0,\;\;\;ax\equiv y\bmod{A},\;\;\;|y|_\infty<1.$$
We claim that 
$$g_x^{(N+1)}(\xi_a)=1+\Psi_N(y).$$
Let $\{a_i\}$, $\{b_j\}$ and $a_0$
be as in
\S\ref{subsubsection:ColemanDefinition}. We have
$$\begin{array}{rcl}
g_x^{(N+1)}(\xi_a)&=&\displaystyle
1-\sum_{i=1}^{\deg f}
\ebold(a_i/f)^{q^{N+1}}C_{a_0b_j}\left(T,\ebold(a/f)\right)
\\\\
&=&\displaystyle
1-\sum_{i=1}^{\deg f}
\ebold(a_i/f)^{q^{N+1}}\ebold(b_jy)\\\\
&=&1+\Psi_N(y),
\end{array}
$$\\
the last equality by part (i) of Theorem~\ref{Theorem:Interpolation}.
The claim is proved.

\begin{Remark}  
The notion of Coleman function
was introduced in \cite{Sinha}, building on the foundation of
\cite{AndersonStick}. The  notion of Coleman function
was in large part inspired by  beautiful examples of
\cite{Coleman}; see \cite[\S2]{AndersonStick} for an
appreciation of Coleman's work. 
The approach to the theory of Coleman functions presented
here is quite a bit simpler than previous approaches and was developed in
an attempt more closely to approximate Coleman's own simple and very
attractive point of view.
\end{Remark}

\subsubsection{Generalized Coleman functions}
\label{subsubsection:GeneralizedColeman} 
Given effective $\abold\in \AAA_f$ such that $\weight \abold>0$, write
$$\abold=\sum_{\begin{subarray}{c}
a\in A\\
\deg a<\deg f
\end{subarray}}m_a\left[\frac{a}{f}\right]\;\;\;\left(m_a\in
\ZZ,\;\;\;m_a\geq 0\right)$$ and put
$$g_{\abold}:=\prod_{\begin{subarray}{c}
0\neq a\in A\\
\deg a<\deg f
\end{subarray}} g_{\frac{a}{f}}^{m_a},
$$
thereby defining a meromorphic function on $\Xbar$ regular on $\Ubar$.
We call $g_{\abold}$ a {\em generalized Coleman function}.
We define divisors of $\Xbar$ by the formulas
$$
\xi_\abold:=\sum_{\begin{subarray}{c} a\in A\\
(a,f)=1\\
\deg a<\deg f
\end{subarray}}\diamondbracket{a\star\abold}\cdot
\xi_{a},\;\;\;\;W_\abold:=\sum_{\begin{subarray}{c} a\in A\\
(a,f)=1\\
\deg a<\deg f
\end{subarray}}\sum_{N=1}^\infty
\diamondbracket{a\star\abold}_N\cdot \left(\sum_{i=0}^{N-1}
\xi_{a}^{(i)}\right).$$
The definition of $W_\abold$ makes sense because only finitely many
nonzero terms appear on the right side.
By the divisor calculation of \S\ref{subsubsection:ColemanDivisor}
we have
$$\begin{array}{rcl}
(g_{\abold})&=&\displaystyle-(\weight
\abold)\cdot
\infty_X+\sum_{\begin{subarray}{c}
a\in A\\
\deg a<\deg f\\
(a,f)=1
\end{subarray}}\sum_{N=0}^\infty
\diamondbracket{a\star \abold}_N\xi_a^{(N)}
\\\\
&=&\displaystyle-(\weight
\abold)\cdot
\infty_X
+\xi_\abold+W^{(1)}_{\abold}-W_{\abold}.
\end{array}$$
By the interpolation
formula of \S\ref{subsubsection:ColemanInterpolation} we have
$$\Pi(a\star\abold)^{-1}=\prod_{N=1}^\infty
g_\abold^{(N)}(\xi_a).$$\\

\begin{Proposition}\label{Proposition:SimplicityCriterionSupplement}
Fix effective $\abold,\bbold\in \AAA_f$ such that $\weight
\abold,\weight\bbold>0$. 
For any deck transformation $\gamma$  and $a\in A$ prime to $f$
corresponding canonically one to the other
in the sense that
$\gamma^*z=C_a(t,z)$, we have
$$\gamma^{-1}\xi_{\abold}=\xi_{\bbold}\Leftrightarrow a\star\abold\sim_f
\bbold.
$$
\end{Proposition}
\proof Clearly
we have $\xi_{\abold}=\xi_{\bbold}\Leftrightarrow
\abold\sim_f\bbold$
and further
$$\gamma^{-1}\xi_\abold=\sum_{\begin{subarray}{c} b\in A\\
(b,f)=1\\
\deg b<\deg f
\end{subarray}}\diamondbracket{b\star\abold}\cdot
\xi_{\bar{a}b}
=\sum_{\begin{subarray}{c} b\in A\\
(b,f)=1\\
\deg b<\deg f
\end{subarray}}\diamondbracket{ab\star\abold}\cdot
\xi_{b}=\xi_{a\star\abold}$$
where $\bar{a}\in A$ satisfies the congruence $a\bar{a}\equiv 1\bmod{f}$,
whence the result.
\qed

\subsection{A construction of 
rigid analytically trivial
 GCM dual $t$-motives}
\label{subsection:ColemanGCM}  For convenience we again put
$\ell:=\card(A/f)^\times$ and we also
fix effective $\abold\in \AAA_f$ such that $\weight\abold>0$.
Put
$$\Lbold:=\FF_q[t,z]/(C_f^\star(t,z)),\;\;\;
\Lboldbar:=\kbar[t,z]/(C_f^\star(t,z)).$$
The rings $\Lbold$ and $\Lboldbar$ are the coordinate rings of
the nonsingular irreducible affine curves curves $U/\FF_q$ and
$\Ubar/\kbar$, respectively. Clearly
$\Lbold$ qualifies as a GCM
$\FF_q[t]$-algebra.  We are going to use the
generalized
 Coleman function $g_\abold$
to create a nice dual $t$-motive $H(\abold)$ with GCM by $\Lbold$.
To do so we
translate to the dual $t$-motivic setting 
a basic construction  that
in the $t$-motivic setting was originally given in \cite{Sinha}.

\begin{Lemma}\label{Lemma:TildeFiniteGeneration}
Let
$\tilde{H}(\abold)$ be the left $\Lboldbar[\sigma]$-module obtained by
equipping $\Lboldbar$ with an action of
$\sigma$ by the rule 
$$\sigma h:=g_\abold h^{(-1)}.$$
Then the $\kbar[\sigma]$-module underlying $\tilde{H}(\abold)$ 
is free of finite rank.
\end{Lemma}
\proof By Proposition~\ref{Proposition:Seesaw}, because the
$\kbar[t]$-module underlying $\tilde{H}(\abold)$ is free of finite rank, 
we have only to prove that the $\kbar[\sigma]$-module underlying
$\tilde{H}(\abold)$ is finitely generated. Temporarily put
$$D:=\weight \abold\cdot \infty_X.$$
Since the multiplicity of $D$ at each point of $\Xbar$ at infinity
is the order of the pole of the generalized Coleman function $g_{\abold}$
at that point according to the divisor formula of
\S\ref{subsubsection:GeneralizedColeman}, the induced maps
$$\frac{\OO_\Xbar\left(nD\right)}
{\OO_\Xbar\left((n-1)D\right)}
\xrightarrow{g_\abold\times}
\frac{\OO_\Xbar\left((n+1)D\right)}{
\OO_\Xbar\left(nD\right)}
$$
of skyscraper sheaves are bijective for all $n\in \ZZ$.
By the Riemann-Roch theorem there exists
an integer
$n_0$ such that the natural maps
$$\frac{H^0\left(\Xbar,\OO_\Xbar\left(nD\right)\right)}
{H^0\left(\Xbar,\OO_\Xbar\left((n-1)D\right)\right)}
\rightarrow
H^0
\left(\Xbar,\frac{\OO_\Xbar\left(nD\right)}
{\OO_\Xbar\left((n-1)D\right)}\right)
$$ are bijective for all $n\geq n_0$.  Clearly we have
$$H^0\left(\Xbar,\OO_\Xbar\left(nD\right)\right)^{(-1)}=
H^0\left(\Xbar,\OO_\Xbar\left(nD^{(-1)}\right)\right)=
H^0\left(\Xbar,\OO_\Xbar\left(nD\right)\right)$$
for all $n\in \ZZ$.
Therefore
$$H^0\left(\Xbar,\OO_\Xbar\left(n D\right)\right)
+g_\abold \cdot H^0\left(\Xbar,\OO_\Xbar\left(n D\right)\right)^{(-1)}
=H^0\left(\Xbar,\OO_\Xbar\left((n+1) D\right)\right)
$$
for all $n\geq n_0$
and hence the vector space
$$H^0\left(\Xbar,\OO_\Xbar\left(n_0 D\right)\right)$$
is finite-dimensional over $\kbar$ and generates
$$\tilde{H}(\abold)=H^0(\Ubar,\OO_{\Xbar})=\bigcup_nH^0\left(\Xbar,\OO_\Xbar\left(n
D\right)\right)$$
over $\kbar[\sigma]$. 
\qed

\subsubsection{Construction of the dual $t$-motive $H(\abold)$}
Recall now the divisor formula
$$(g_\abold)=-(\weight
\abold)\cdot
\infty_X
+\xi_\abold+W^{(1)}_{\abold}-W_{\abold}$$
 of
\S\ref{subsubsection:GeneralizedColeman}
and recall also that the divisors
$\xi_\abold$ and $W_{\abold}$ figuring in this formula are effective.
 Put
$$H(\abold):=H^0\left(\Ubar,\OO_\Xbar\left(-W_\abold^{(1)}\right)\right)
\subset
H^0\left(\Ubar,\OO_{\Xbar}\right)=\tilde{H}(\abold)$$
thereby defining an $\Lboldbar$-submodule of $\tilde{H}(\abold)$. It is
easy to verify that $H(\abold)$ is $\sigma$-stable and hence an
$\Lboldbar[\sigma]$-submodule of $\tilde{H}(\abold)$.
It is clear that $H(\abold)$ is projective over $\Lboldbar$ of rank
one and free of finite rank 
over
$\kbar[t]$. Moreover
$H(\abold)$ is a $\kbar[\sigma]$-submodule of a
$\kbar[\sigma]$-module free of finite rank by
Lemma~\ref{Lemma:TildeFiniteGeneration} and hence a free
$\kbar[\sigma]$-module of finite rank. Since we have
$$\frac{H(\abold)}{\sigma
H(\abold)}=
\frac{H^0(\Ubar,\OO_\Xbar(-W_\abold^{(1)}))}{g_\abold\cdot
H^0(\Ubar,\OO_\Xbar(-W_\abold^{(1)}))^{(-1)}}=\frac{H^0(\Ubar,\OO_\Xbar(-W_\abold^{(1)}))}{
H^0(\Ubar,\OO_\Xbar(-\xi_{\abold}-W_\abold))},$$ 
and all the points in the support of the divisor $\xi_\abold$
lie above the point $t=T$ on the $t$-line,
it follows that $H(\abold)/\sigma H(\abold)$
is annihilated by a sufficiently high power of $t-T$.
Therefore $H(\abold)$ is a dual $t$-motive
with GCM by $\Lbold$. Note that the ideal
$$I_\abold:=H^0\left(\Ubar,\OO_\Xbar(-\xi_{\abold})\right)\subset\Lboldbar$$
is the GCM type of $H(\abold)$ with respect to $\Lbold$.

\begin{Lemma}\label{Lemma:CriticalMatrixCalculation}
(i) Let
$$\Phi_\abold\in \Mat_{\ell\times
\ell}(\kbar[t])$$ be the unique solution
of the congruence
$$g_{\abold}\left[\begin{array}{c}
1\\
z\\
\vdots\\
z^{\ell-1}
\end{array}\right]\equiv\Phi_\abold \left[\begin{array}{c}
1\\
z\\
\vdots\\
z^{\ell-1}
\end{array}\right]\bmod{C_f^\star(t,z)}.$$
Then the infinite product
$$\Psi_\abold:=\prod_{N=1}^\infty \Phi^{(N)}_\abold$$
converges with respect to Banach norm 
$$\left\Vert\sum_{i=0}^\infty a_it^i\right\Vert_\infty:=
\sup_{i=0}^\infty|a_i|_\infty$$ on $\CC_\infty\{t\}$
to an element of
$$\GL_\ell(\CC_\infty\{t\})\cap\Mat_{\ell\times \ell}(\EE)$$
satisfying the functional equation
$$\Psi^{(-1)}_\abold=\Phi_\abold\Psi_\abold.$$
(ii) Consider now the matrix 
$$\Psi_\abold(T)\in\Mat_{\ell\times \ell}(\overline{k_\infty})$$
obtained by
evaluating
$\Psi_\abold$ at $t=T$.
The sets
$$\{\Psi_\abold(T)_{ij}\mid
i,j=1,\dots,\ell\},\;\;\;\{\Pi(a\star\abold)^{-1}\mid a\in A,\;
(a,f)=1\}$$ span the same $\kbar$-subspace of $\overline{k_\infty}$.
\end{Lemma}
\proof (i) 
{From} the construction of the
generalized Coleman function $g_\abold$ it is clear that we have
$$g_\abold\equiv 1+\sum_i\sum_j c_{ij}t^iz^j\bmod{C_f^\star(t,z)}$$
for some constants $c_{ij}\in \kbar$, all but finitely many of which vanish
and all of which satisfy the bound 
$$|c_{ij}|_\infty\leq \sup_{x\in k_\infty}|\ebold^*(x)|_\infty=
|1/\Tbar|_\infty<1.$$ 
Now let
$$Z=Z(t)\in\Mat_{\ell\times \ell}(\FF_q[t])$$
be the unique solution of the congruence
$$z\left[\begin{array}{c}
1\\
z\\
\vdots\\
z^{\ell-1}
\end{array}\right]\equiv Z\left[\begin{array}{c}
1\\
z\\
\vdots\\
z^{\ell-1}
\end{array}\right]\bmod{C_f^\star(t,z)}.$$
Clearly we have
$$\Phi_{\abold}^{(N)}=\one_\ell+\sum_i\sum_j c_{ij}^{q^N}t^iZ^j$$
for all integers $N$.
It follows that the infinite product defining $\Psi_\abold$ converges 
to an element of $\GL_\ell(\CC_\infty\{t\})$ satisfying
the desired functional equation. It follows also
that $\det\Phi_{\abold}(0)\neq 0$
and hence by  Proposition~\ref{Proposition:HypothesisChecker}
that the matrix $\Psi_\abold$ has entries
in $\EE$.

(ii) Note that by construction the roots of the equation $C_f^\star(T,z)=0$ give the
$\ell$ eigenvalues of the matrix
$$Z(T)\in\GL_\ell(\kbar).$$
Now choose any matrix 
$$M\in \GL_\ell(\kbar)$$
such that
$$MZ(T)M^{-1}=\left[\begin{array}{ccc}
\ebold(a_1/f)\\
&\ddots\\
&&\ebold(a_\ell/f)
\end{array}\right]$$
where
$$\{a_1,\dots,a_\ell\}=\{a\in A\mid \deg a<\deg f,\;(a,f)=1\}.$$
Then
$$\left(M\Phi_\abold^{(N)}(T)M^{-1}\right)_{ij}=g^{(N)}(\xi_{a_i})\cdot
\delta_{ij},$$
and hence
$$(M\Psi_\abold(T)M^{-1})_{ij}=\Pi(a_i\star
\abold)^{-1}\cdot\delta_{ij},$$ which proves the result.
\qed

 \begin{Proposition}\label{Proposition:Last}
Let
$$\ggg\in \Mat_{\ell\times
1}(H(\abold)),\;\;\;\;\Phi\in \Mat_{\ell\times \ell}(
\kbar[t])$$
be given such that the entries of $\ggg$ form a $\kbar[t]$-basis
of $H(\abold)$ and
$$\sigma \ggg=\Phi\ggg.$$
There exists
$$\Psi\in \GL_\ell(\CC_\infty\{t\})\cap \Mat_{\ell\times
\ell}(\EE)$$
satisfying the functional equation
$$\Psi^{(-1)}=\Phi\Psi$$
and with the further property that the sets
$$\{\Psi(T)_{ij}\mid i,j=1,\dots,\ell\},\;\;\;
\{\Pi(a\star\abold)^{-1}\mid a\in A,\;(a,f)=1\}$$
span the same $\kbar$-subspace of $\overline{k_\infty}$.
\end{Proposition}
In particular
$H(\abold)$ is rigid analytically trivial
by Lemma~\ref{Lemma:RATEquation}.
\proof
Let
$$
Q\in \Mat_{\ell\times \ell}(\kbar[t])$$
be the unique solution of the congruence
$$\ggg\equiv Q\left[\begin{array}{c}
1\\
z\\
\vdots\\
z^{\ell-1}
\end{array}\right]\bmod{C_f^\star(t,z)}.$$
Now since the effective divisor $W_\abold^{(1)}$
is supported in the set of $\kbar$-valued points of $X$
lying above the points
$$t=T^q,T^{q^2},\dots$$
on the $t$-line, the module
$$\frac{\tilde{H}(\abold)}{H(\abold)}
=\frac{H^0(\Ubar,\OO_{\Xbar})}
{H^0(\Ubar,\OO_\Xbar(-W_\abold^{(1)}))}$$
is annihilated by
$$\prod_{i=1}^N(t-T^{q^i})^N$$
for $N\gg 0$
and hence we have
$$Q\in\GL_\ell(\CC_\infty\{t\}),\;\;\;\det Q(T)\neq 0.$$
Notation as in Lemma~\ref{Lemma:CriticalMatrixCalculation},
we
have
$$\Phi Q\left[\begin{array}{c}
1\\
z\\
\vdots\\
z^{\ell-1}
\end{array}\right]\equiv\sigma\ggg\equiv 
Q^{(-1)}\Phi_\abold\left[\begin{array}{c} 1\\
z\\
\vdots\\
z^{\ell-1}
\end{array}\right]\bmod{C_f^\star(t,z)},$$
hence
$$Q^{(-1)}\Phi_\abold=\Phi Q,$$
and hence the matrix
$$\Psi:=Q\Psi_\abold\in
\GL_\ell(\CC_\infty\{t\})\cap\Mat_{\ell\times
\ell}(\EE)$$
has all the desired properties.
\qed
\begin{Proposition}\label{Proposition:SimplicityClincher}
Let $r$ be the cardinality of the orbit of the coset
$\abold\bmod{\widetilde{\RRR}_f}$ under the action of $(A/f)^\times$.
Then any simple quotient of the bare dual $t$-motive underlying
$H(\abold)$ is of rank
$\geq r$ over
$\kbar[t]$.
\end{Proposition}
\proof We temporarily denote the group of deck transformations by $G$. We
have
$$\card\{a\star \abold\bmod{\widetilde{\RRR}_f}\mid a\in A,\;(a,f)=1\}
=\card\{\gamma\xi_\abold\mid \gamma\in G\}
=\card\{\gamma^*I_\abold\mid \gamma\in G\},$$
where the first equality holds by
Proposition~\ref{Proposition:SimplicityCriterionSupplement}
and the second equality is trivial. The quantity on the right by
Corollary~\ref{Corollary:SimplicityCriterion} is a lower bound for the
rank over $\kbar[t]$ of any simple quotient of the bare dual $t$-motive
underlying $H(\abold)$.
\qed
\begin{Proposition}\label{Proposition:NonisogenyClincher}
For all $\bbold\in \AAA_f$ such that
$\weight\bbold>0$, if the bare dual $t$-motives underlying $H(\abold)$
and $H(\bbold)$ have isogenous simple quotients,
then 
$a\star\abold\sim_f \bbold$
for some $a\in A$ prime to $f$.
\end{Proposition}

\proof Once again let $G$ denote the group of deck transformations.
Consider the following sets:
\begin{gather*}
\{\gamma^*I_\abold=I_\bbold\mid \gamma\in G\}, \\
\{\gamma \xi_\abold=\xi_\bbold\mid \gamma\in G\}, \\
\{a\star \abold\equiv
\bbold\bmod{\widetilde{\RRR}_f}\mid a\in A,\;(a,f)=1\}.
\end{gather*}
The first set is nonempty by
Corollary~\ref{Corollary:NonisogenyCriterion}.
Nonemptiness of the first set trivially
implies that of the second.
Nonemptiness of the second set implies that of the third
by Proposition~\ref{Proposition:SimplicityCriterionSupplement}.
\qed

\subsection{Proof of
Theorem~\ref{Theorem:LinearRelations} and 
Corollary~\ref{Corollary:TranscendenceDegree}}\label{subsection:LooseEnds}
We need only prove the sufficiency asserted in the theorem
because the diamond bracket criterion takes care of necessity
and Proposition~\ref{Proposition:Strategy} takes care of the corollary.
\subsubsection{Easy
reductions} 
By hypothesis:
\begin{itemize}
\item The cosets $\abold_1\bmod{\widetilde{\RRR}_f},\dots,
\abold_N\bmod{\widetilde{\RRR}_f}$
are distinct.
\end{itemize}
After enlarging the set
$\{\abold_1,\dots,\abold_N\}$ suitably we may assume that:
\begin{itemize}
\item
The finite set $\{\abold_i\bmod{\widetilde{\RRR}_f}\mid i=1,\dots,N\}
\subset \AAA_f/\widetilde{\RRR}_f$
is $(A/f)^\times$-stable.
\end{itemize}
After relabeling the $\abold_i$ we may assume that
for some integer $1\leq n\leq N$:
\begin{itemize}
\item The set 
$\{\abold_i\bmod{\widetilde{\RRR}_f}\mid i=1,\dots,n\}$
forms a set of representatives for the $(A/f)^\times$-orbits
in $\{\abold_i\bmod{\widetilde{\RRR}_f}\mid i=1,\dots,N\}$.
\end{itemize}
For $i=1,\dots,n$ put\\
$$\begin{array}{rcl}
r_i&:=&\card\{a\star \abold_i\bmod{\widetilde{\RRR}_f}\mid
a\in A,\;(a,f)=1\},\\\\
V_i&:=&\left(\mbox{$\kbar$-span of $\{\Pi(a\star \abold_i)^{-1}\mid a\in
A,\;(a,f)=1\}$}\right)\subset\overline{k_\infty}.
\end{array}
$$ \\
We have $N=\sum_i r_i$. Moreover the $\kbar$-span of
$\{\Pi(-\abold_i)\mid i=1,\dots,N\}$ equals $\sum_i V_i$
by the diamond bracket criterion. It therefore suffices to
show that
$$\sum_{i=1}^n r_i\leq\dim_\kbar \sum_{i=1}^n
V_i.$$ 
After adding a fixed positive integral
multiple of
$$\sum_{\begin{subarray}{c}
a\in A\\
\deg a<\deg f
\end{subarray}}\left[\frac{a}{f}\right]\in \AAA_f$$
to $\abold_i$ for all $i$, we may further assume that\\
\begin{itemize}
\item $\abold_i$ is
effective and 
 $\weight \abold_i>0$ for $i=1,\dots,n$.
\end{itemize}

\subsubsection{Further reductions}
For $i=1,\dots,n$ we make the following constructions.
Put $$H_i:=H(\abold_i).$$
Choose 
$$\ggg_{(i)}\in \Mat_{\ell\times
1}(H_i),\;\;\;
\Phi_{(i)}\in \Mat_{\ell\times
\ell}(\kbar[t])\;\;\;(\ell:=\card(A/f)^\times)$$
such that the entries of $\ggg_{(i)}$ form a $\kbar[t]$-basis of $H_i$
and
$$\sigma \ggg_{(i)}=\Phi_{(i)}\ggg_{(i)}.$$
By Proposition~\ref{Proposition:Last} there exists
$$\Psi_{(i)}\in \GL_\ell(\CC_\infty\{t\})\cap\Mat_{\ell\times \ell}(\EE)$$
such that the functional equation
$$\Psi_{(i)}^{(-1)}=\Phi_{(i)}\Psi_{(i)}$$
holds. Let $\psi_{(i)}$ be
the first column of $\Psi_{(i)}$. Put\\
$$
\begin{array}{rcl}
H_{i0}&:=&\mbox{$\kbar[t]$-span in $\EE$ of the entries of
$\psi_{(i)}$,}\\\\
V_{i0}&:=&\mbox{$\kbar$-span in $\overline{k_\infty}$ of the entries of
$\psi_{(i)}(T)$.}
\end{array}
$$\\
By Proposition~\ref{Proposition:Last}
 we have
\begin{itemize}
\item $V_{i0}\subset V_i$ for $i=1,\dots,n$.
\end{itemize}
By Proposition~\ref{Proposition:Rationale} we have
\begin{itemize}
\item $\displaystyle \rank_{\kbar[t]}\sum_{i=1}^n H_{i0}
=\dim_\kbar \sum_{i=1}^n V_{i0}$.
\end{itemize}
It therefore suffices to prove that
$$\sum_{i=1}^n r_i\leq \rank_{\kbar[t]}\sum_{i=1}^n H_{i0}.$$

\subsubsection{Endgame}
In a sense made precise by Proposition~\ref{Proposition:Rationale}:
\begin{itemize}
\item $H_{i0}$ is a 
nonzero
 dual $t$-motive admitting presentation as
a quotient of  the bare dual $t$-motive underlying $H_i$ for
$i=1,\dots,n$.
\end{itemize}
We  have
\begin{itemize}
\item $r_i\leq \rank_{\kbar[t]}H_{i0}= \dim_\kbar V_{i0}\leq
\dim_\kbar V_i\leq
r_i$, for $i=1,\dots,n$,
\end{itemize}
by Proposition~\ref{Proposition:SimplicityClincher}
(inequality at the extreme left), Proposition~\ref{Proposition:Rationale}
(equality at second juncture)
and the diamond bracket criterion (inequality at the extreme
right). It follows that
\begin{itemize}
\item $H_{i0}$ is simple and of rank $r_i$ over $\kbar[t]$ for
$i=1,\dots,n$.
\end{itemize}
The simple dual $t$-motives $H_{i0}$ belong to distinct
isogeny classes by Proposition~\ref{Proposition:NonisogenyClincher}
and
hence:
\begin{itemize}
\item The natural map $\displaystyle\bigoplus_{i=1}^n
H_{i0}\rightarrow\sum_{i=1}^n H_{i0}$ is bijective.
\end{itemize}
 Therefore we have
$$\sum_{i=1}^n r_i=
\sum_{i=1}^n\rank_{\kbar[t]}H_{i0}=\rank_{\kbar[t]}\sum_{i=1}^n
H_{i0}$$ and the proof of sufficiency in
Theorem~\ref{Theorem:LinearRelations} is finished.
As noted above, with sufficiency proved,
the proofs of Theorem~\ref{Theorem:LinearRelations}
and Corollary~\ref{Corollary:TranscendenceDegree}
are complete.

\subsection{Remarks concerning transcendence bases}
Giving an explicit transcendence basis for the field $E_f$ over $\kbar$ is
in general not as straightforward as one might suspect from the
statements of Theorem~\ref{Theorem:LinearRelations} and
Corollary~\ref{Corollary:TranscendenceDegree}. 
In principle transcendence bases and systems of relations
could be constructed explicitly
by a translation to our context of the methods of \cite{Kucera},
but we do not attempt such a construction here.
We just work out the special case in which $f$ is a power of an
irreducible polynomial and then give a cautionary example.

\begin{Proposition} \label{Proposition:IrreduciblePowerBasis}
Let $f_1 \in A_+$ be irreducible, and suppose $f = f_1^s$ for some
$s$.  Let
\[
  B_f = \{ \Pi(a/f) \mid a \in A\setminus A_+, \deg a < \deg f, (a,f) =
  1 \} \cup \{ \varpi \}.
\]
Then $E_f = \kbar(B_f)$, and the numbers in $B_f$ are algebraically
independent over $\kbar$.
\end{Proposition}
\proof
Consider $\kappa f_1^e/f$, where $(\kappa,f_1) = 1$
and $\deg \kappa < (s-e)\deg f_1$.  As an element of $\DDD_f$,
\begin{align*}
\left[\frac{\kappa f_1^e}{f}\right] -
\sum_{\deg a < e\deg f_1} \left[\frac{\kappa f_1^e/f + a}{f_1^e}\right]
&=\left[\frac{\kappa f_1^e}{f}\right] -
\sum_{\deg a < e\deg f_1} \left[\frac{\kappa f_1^e + af}{f_1^ef}\right] \\
&=\left[\frac{\kappa f_1^e}{f}\right] -
\sum_{\deg a < e\deg f_1} \left[\frac{\kappa + af_1^{s-e}}{f}\right].
\end{align*}
Notice now that the terms on the right have numerators congruent to
$\kappa$ modulo a power of $f_1$.  Thus the numerators are relatively
prime to $f$.  According to the diamond bracket criterion we can therefore
express every special
$\Pi$-value $\Pi(b/f)$ as a
$\kbar$-multiple of a product of $\Pi(a/f)$, with $a$ relatively prime
to $f$, divided by a power of $\varpi$.  Finally, for every $a \in
A_+$, the reflection identity dictates that
\[
  \varpi^{-1}\prod_{c \in \FF_q^{\times}} \Pi(ca/f) \in \kbar^{\times},
\]
and so we conclude that $E_f = \kbar(B_f)$ as claimed.  Since $\# B_f =
\nu_f$ is the transcendence degree of $E_f$ over $\kbar$, the rest
follows.
\qed

\subsubsection{Cautionary example}
\label{subsubsection:ParthianShot}
In light of Proposition~\ref{Proposition:IrreduciblePowerBasis}, it would
not be far-fetched to imagine that 
\[
  B_f = \{ \Pi(a/f) \mid a \in A\setminus A_+, \deg a < \deg f, (a,f) =
  1 \} \cup \{ \varpi \}.
\]
would provide a trancendence basis for $E_f$ over $\kbar$ for
\emph{all} $f \in A_+$.  That however is not the always the case.

Consider the example of $q = 3$ and $f = T^2 - T$.  By
Corollary~\ref{Corollary:TranscendenceDegree} the transcendence degree of
$E_f$ over $\kbar$ is $3$.  In \S 4.2 of \cite{SinhaDR} it is shown
that
\[
  \Pi\left( \frac{1}{T^2 - T} \right) \bigg/ \Pi\left( \frac{1}{T}
  \right) \in \kbar^{\times},
\]
and 
\[
  \Pi\left( \frac{T+1}{T^2 - T} \right) \bigg/ \Pi\left( \frac{1}{T}
  \right) \in \kbar^{\times},
\]
by applying the diamond bracket criterion. 
Consequently, in view of the reflection identity satisfied by the
$\Pi$-function,
$\kbar(B_f)$ has transcendence degree at most $2$ over $\kbar$. But by
Corollary~\ref{Corollary:TranscendenceDegree} and the diamond bracket
criterion, we know that $E_f$ is the rational function field
\[
  E_f = \kbar\left( \varpi, \Pi\left( \frac{1}{T} \right), \Pi\left(
  \frac{1}{T-1} \right) \right).
\]
\subsubsection{Remark}
The $\Pi$-monomial 
$\Pi\left( \frac{1}{T^2 - T} \right) 
/ \Pi\left( \frac{1}{T} \right)$ and its companion 
happen to be examples of the  sort considered
above in \S\ref{subsubsection:DangerousBend}.

\end{document}